\documentclass[11pt,letterpaper]{amsart}

\usepackage[centertags]{amsmath}
\usepackage{amsthm}
\usepackage{amssymb}
\usepackage{amscd}
\usepackage{eucal}
\usepackage{epsfig}
\usepackage{verbatim}

\usepackage{latexsym,enumerate,amsxtra}

\newtheorem{thm}{Theorem}[section]
\newtheorem{theorem}[thm]{Theorem}
\newtheorem{corollary}[thm]{Corollary}
\newtheorem{lemma}[thm]{Lemma}
\newtheorem{proposition}[thm]{Proposition}
\newtheorem{definition}[thm]{Definition}

\newtheorem{remarks}[thm]{Remarks}
\newtheorem{remark}[thm]{Remark}

\newtheorem{numbering}[thm]{}
\newtheorem{defrmk}[thm]{Definition and Remark}

\newcommand{\CaA}{\mathcal A}
\newcommand{\CaB}{\mathcal B}

\newcommand{\CaE}{\mathcal E}

\newcommand{\CaO}{\mathcal O}
\newcommand{\CaH}{\mathcal H}

\newcommand{\CaL}{\mathcal L}

\newcommand{\frt}{\mathfrak t}
\newcommand{\frS}{\mathfrak S}
\newcommand{\frB}{\mathfrak B}

\newcommand{\bbC}{\mathbb C}
\newcommand{\bbF}{\mathbb F}
\newcommand{\bbH}{\mathbb H}
\newcommand{\bbQ}{\mathbb Q}
\newcommand{\bbR}{\mathbb R}
\newcommand{\bbW}{\mathbb W}
\newcommand{\bbZ}{\mathbb Z}
\newcommand{\tbR}{\tilde\bbR}

\newcommand{\rA}{A}
\newcommand{\rT}{\mathrm T}

\newcommand{\rH}{H}

\newcommand{\rd}{\mathrm d}

\newcommand{\xo}{y}
\newcommand{\yo}{y}
\newcommand{\xx}{x}
\newcommand{\Xn}{X}
\newcommand{\Set}{S}

\newcommand{\sfC}{\mathsf C}

\newcommand{\sfM}{\mathsf M}
\newcommand{\sfP}{\mathsf P}

\newcommand{\sfS}{\mathsf S}
\newcommand{\sfK}{\mathsf K}

\newcommand{\val}{\nu}
\newcommand{\Ad}{\textup{Ad}}
\newcommand{\ad}{\textup{ad}}
\newcommand{\pid}{\mathfrak p}
\newcommand{\orb}{\mathcal O}

\newcommand{\Tr}{\textup{Tr}}

\newcommand{\bG}{\mathbf G}
\newcommand{\bH}{\mathbf H}

\newcommand{\bC}{\mathbf C}
\newcommand{\bT}{\mathbf T}
\newcommand{\bM}{\mathbf M}
\newcommand{\bP}{\mathbf P}
\newcommand{\bN}{\mathbf N}
\newcommand{\bU}{\mathbf U}
\newcommand{\bS}{\mathbf S}

\newcommand{\bK}{\sfK}
\newcommand{\bZ}{\mathbf Z}
\newcommand{\bX}{\mathbf X}
\newcommand{\bGL}{\mathbf{GL}}
\newcommand{\bSL}{\mathbf{SL}}
\newcommand{\T}{T}

\newcommand{\C}{C}
\newcommand{\D}{D}
\newcommand{\G}{G}
\newcommand{\rM}{M}
\newcommand{\rU}{U}
\newcommand{\rN}{N}
\newcommand{\rP}{P}
\newcommand{\rL}{L}
\newcommand{\rK}{K}
\newcommand{\rJ}{J}
\newcommand{\heish}{h}
\newcommand{\rh}{{\mathrm{h}}}
\newcommand{\tx}{w_1}
\newcommand{\txx}{w}
\newcommand{\ty}{w_2}

\newcommand{\Jv}{K_{\vdash}}

\newcommand{\Jvy}{K_{y\vdash}}
\newcommand{\tJ}{\tilde J}

\newcommand{\lieG}{\mathfrak g}
\newcommand{\lieH}{\mathfrak h}
\newcommand{\lieZ}{\mathfrak z}
\newcommand{\lieM}{\mathfrak m}
\newcommand{\lieN}{\mathfrak n}
\newcommand{\lieU}{\mathfrak u}
\newcommand{\lieT}{\mathfrak t}

\newcommand{\blieG}{\boldsymbol{\lieG}}
\newcommand{\blieH}{\boldsymbol{\lieH}}
\newcommand{\blieT}{\boldsymbol{\lieT}}
\newcommand{\blieZ}{\boldsymbol{\lieZ}}

\newcommand{\der}{\mathrm{der}}
\newcommand{\gss}{\Gamma}
\newcommand{\sgss}{\Gamma}
\newcommand{\gsc}{\mathcal S}
\newcommand{\sgsc}{\frS}
\newcommand{\datum}{\Sigma}
\newcommand{\Index}{\Upsilon_\datum}
\newcommand{\good}{good}

\newcommand{\cuspd}{\rho}
\newcommand{\cuspf}{\tau}

\newcommand{\conv}{\ast}
\newcommand{\addch}{\Lambda}
\newcommand{\bilinear}{\mathrm B}
\newcommand{\Bd}{\CaB}

\newcommand{\Apt}{\CaA}

\newcommand{\chamberD}{D}

\newcommand{\dth}{\rd}

\newcommand{\dpi}{\varrho}
\newcommand{\rtm}{r}

\newcommand{\utm}{u}
\newcommand{\atm}{a}
\newcommand{\btm}{b}
\newcommand{\tm}{t}
\newcommand{\stm}{s}

\newcommand{\charpi}{\Theta_\pi}
\newcommand{\Fcharpi}{\widehat{\Theta}_\pi}
\newcommand{\Ff}{\widehat f}
\newcommand{\resf}{f_{0^+}}

\newcommand{\Ccs}{C_c^\infty}
\newcommand{\distr}{\mathcal J}

\newcommand{\supp}{\textup{Supp}}

\newcommand{\nil}{\mathcal N}

\newcommand{\Ind}{\textup{Ind}}
\newcommand{\Gal}{\textup{Gal}}

\newcommand{\admG}{\CaE}
\newcommand{\dualG}{\CaE^u}
\newcommand{\temG}{\CaE^t}
\newcommand{\dualg}{\widehat{\lieG}}
\newcommand{\adualg}{\lieG^\ast}
\newcommand{\eps}{\epsilon}
\newcommand{\slf}{{F}}
\newcommand{\vol}{\textit{vol}}
\newcommand{\mult}{m}
\newcommand{\type}{{\mathfrak s}}

\newcommand{\Hypk}{\textup{(H$k$)}}
\newcommand{\HypB}{\textup{(HB)}}
\newcommand{\HypGT}{\textup{(HGT)}}
\newcommand{\HypN}{\textup{(H$\nil$)}}

\newcommand{\itt}{\alpha}

\newcommand{\slft}{\slf_\itt}
\newcommand{\lp}{\left(}
\newcommand{\rp}{\right)}

\newcommand{\midvsp}{\vspace{7pt}}

\numberwithin{equation}{section}

\begin{document}

{ }

\title{supercuspidal representations: an exhaustion theorem}
\author{Ju-Lee Kim}

\begin{abstract}
Let $\G$ be a reductive $p$-adic group.
We prove that all supercuspidal representations of $\G$ arise 
through Yu's construction subject to certain hypotheses on $k$ 
(depending on $\G$).
As a corollary, under the same hypotheses, we see that any
supercuspidal representation is compactly induced from a representation of
an open subgroup which is compact modulo the center.
\end{abstract}

\thanks{}
\subjclass{Primary 22E50; Secondary 22E35, 20G25.}
\keywords{supercuspidal representations, $\bK$-types, $p$-adic groups}

\address
{Department of Mathematics\\University of Illinois at Chicago
\\Chicago, IL 60607}
\email{julee@math.uic.edu}

\vspace{.7cm}

\maketitle

\section*{Introduction}

Let $k$ be a $p$-adic field of characteristic zero and residue characteristic
$p$. Let $\G$ be the group of $k$-points of a connected reductive group 
$\bG$ defined over $k$.
In \cite{yu}, Yu gives a fairly general construction of
supercuspidal representations of $\G$ in a certain tame situation.
In this paper, subject to some hypotheses on $\bG$ and $k$,
we prove that all supercuspidal representations arise 
through his construction.

While there have been numerous constructions of supercuspidal 
representations, the question of whether they are exhaustive is resolved 
only for depth zero representations \cite{MP2, Mr}
and for groups of type $A_n$ 
such as $\bGL_n$ \cite{BK, HM, Moy}, $\bSL_n$ \cite{BKsl1, BKsl2}.
In \cite{Moy}, Moy proves the exhaustiveness for Howe's construction 
of supercuspidal representations via the generalized Jacquet-Langlands 
correspondence in the tame case. 
In \cite{HM}, Howe and Moy prove it by analyzing Hecke algebras when $p>n$.
In \cite{BK}, Bushnell and Kutzko construct supercuspidal representations 
and prove their exhaustiveness by analyzing simple types and split types 
with no assumption on $k$. 
Recently, Stevens showed that any supercuspidal
representation of classical groups of positive depth
contains a certain semisimple characters \cite{Stevens}.
However, since no analogue of the Jacquet-Langlands correspondence
for general groups has been developed yet, and 
since types, or Hecke algebras for general groups,
are far less understood than 
for $\bGL_n$, it is not easy to extend their methods to other groups.

In this paper, we approach this problem via harmonic analysis on $\G$. 
We now briefly describe the main idea of the proof.
From now on, we assume that the residue characteristic $p$ of $k$ is
sufficiently large (see (\ref{num: hypotheses}) for the precise condition).
We first prove that any supercuspidal representation is either of depth zero,
or otherwise contains a $\bK$-type constructed in \cite{Asym2}. 
This we do by relating the Plancherel formulas on $\G$ and 
on its Lie algebra $\lieG$ 
and by using some results on asymptotic expansions \cite{Asym2}.
We relate the $\bK$-type further to a supercuspidal type 
constructed in \cite{yu} by analyzing appropriate Hecke algebras 
and Jacquet modules.
Before expanding our account of the main strategy of the proof, 
we first recall some results on $\gss$-asymptotic expansions.

\medskip

\subsection{\sc Results on $\gss$-asymptotic expansions}
Let $\admG:=\admG(\G)$ denote the set of all equivalence classes of
irreducible admissible representations of $\G$. We use the same notation
for a representation $\pi$ and its equivalence class. 
For $(\pi,V_\pi)\in\admG$, let $\charpi$ be the character of $\pi$. 
Let $\Bd(\bG,k)$ be the extended Building of $\bG$ over $k$.
In \cite{Asym2}, we found a certain character expansion
of $\charpi$ depending on $\bK$-types contained 
in $(\pi,V_\pi)$. The construction of $\bK$-types is 
based on a (strongly) good positive $\G$-datum $\datum$ which is a quadruple
(see (\ref{def: posGdatum}) and (\ref{def: goodatum})).
However, if $\datum$ is strongly good, $\datum$ can be 
alternatively described as a pair $\datum=(\gss,\xo)$ of a semisimple element 
$\gss\in\lieG$ and $\xo\in\Bd(\bC_\bG(\gss),k)$ satisfying the following 
(see (\ref{def: posGdatum}), (\ref{def: goodatum}) and (\ref{rmk: xtoy})):

$(i)$ $\gss=\gss_d+\gss_{d-1}+\cdots+\gss_0$ where 
$\gss_i$, $0\le i\le d-1$, is a $\G$-good element of depth $-\rtm_i$
and $\gss_d$ is either zero or an element in the center of $\lieG$ 
of depth $-\rtm_d$. Set $\rtm_d:=\rtm_{d-1}$ if $\gss_d=0$.

$(ii)$ $0<\rtm_0<\rtm_1<\cdots<\rtm_{d-1}$. 
If $\gss_d\ne0$, $\rtm_{d-1}<\rtm_d$.

$(iii)$ $\bC_\bG(\gss)=\bG^0\subsetneq\bG^1
\subsetneq\cdots\subsetneq\bG^{d-1}\subsetneq\bG$
where $\bG^i=\bC_\bG(\gss_d+\gss_{d-1}+\cdots\gss_i)$.

\noindent
Set $\stm_i:=\frac{\rtm_i}2$ and $\vec\stm^+
:=(0^+,\stm_0^+,\stm_1^+,\cdots,\stm_{d-1}^+)$.
Then the associated $\bK$-type (which we denote by 
$(\rK^+_\datum,\phi_\datum)$ in the text) 
is $(\vec\G_{\xo,\vec\stm^+},\chi_\gss)$, 
where $\vec\G_{\xo,\vec\stm^+}$ is an open compact subgroup
defined in \cite{yu} and $\chi_\gss$ is 
the character on $\vec\G_{\xo,\vec\stm^+}$
represented by $\gss$ via a logarithmic map.
Any irreducible admissible representation containing $\chi_\gss$
when restricted to $\vec\G_{\xo,\vec\stm^+}$ has depth $\dpi=\rtm_{d}>0$.
Moreover, such a representation is not supercuspidal in general, that is,
$(\vec\G_{\xo,\vec\stm^+},\chi_\gss)$ is not necessarily a supercuspidal type.

The main result of \cite{Asym2} states that 
if $\pi$ contains $(\vec\G_{\xo,\vec\stm^+},\chi_\gss)$,
for any $f\in\Ccs(\lieG_{\stm_{d-1}^+})$, then
$\charpi(f\circ\log)
=\sum_{\orb\in\orb(\gss)}c_\orb(\pi)\widehat{\mu_\orb}(f).$
Here, $\lieG_{\rtm^+}=\cup_{x\in\Bd(\bG,k)}\lieG_{x,\rtm^+}$, $\rtm\in\bbR$,
$\orb(\gss)$ is the set of $\G$-orbits whose closures contain $\gss$,
and $\widehat{\mu_\orb}$ is the Fourier transform of the orbital integral
$\mu_\orb$. In \cite{Asym2}, we also define a certain subspace 
$\distr^\gss$ of the space of $\G$-invariant distributions on $\lieG$ 
having the property that $\distr^\gss$, when restricted to the image of  
the Fourier transform of $\Ccs(\lieG_{\stm_{d-1}^+})$,
coincides with the finite dimensional space spanned by $\mu_\orb$, 
$\orb\in\orb(\gss)$.
An important property of this expansion is that 
there are test functions $f_\orb$ in $\Ccs(\lieG)$ 
indexed by $\orb\in\orb(\gss)$ such that, for two $\G$-invariant distributions 
$\rT_1$ and $\rT_2$ on $\lieG$ with 
their Fourier transforms $\widehat\rT_1$, $\widehat\rT_2$ in $\distr^\gss$,
if $\widehat\rT_1(f_\orb)=\widehat\rT_2(f_\orb)$ 
for all $\orb\in\orb(\gss)$, then $\rT_1\equiv\rT_2$ on 
$\Ccs(\lieG_{\stm_{d-1}^+})$ (see (\ref{thm: test ftns})).

\subsection{\sc First step}\label{subsec: first step}
Let $\temG$ be the set of equivalence classes of 
irreducible tempered representations. We first show that 
{\em almost every irreducible tempered representation $(\pi,V_\pi)$ is either
of depth zero, or otherwise contains $(\vec\G_{\xo,\vec\stm^+},\chi_\gss)$
for some $(\gss,\xo)$.}

We begin by observing that, thanks to  
the Plancherel formulas on $\lieG$ and $\G$ (\cite{Dx}), we have the equality
\begin{equation}\label{eq: Pl-intro}
\int_{\lieG}\widehat f(X)\,dX=f(0)=\int_{\temG}\charpi(f\circ\log)\,d\pi
\end{equation}
for $f\in\Ccs(\lieG)$ supported in a small neighborhood of $0$. Here, note 
that we have identified $\lieG$ with the unitary dual $\dualg$ of $\lieG$.
Now, refining this equality,
we will find a matching between spectral decomposition factors of each side
of (\ref{eq: Pl-intro}), parameterized by some equivalence 
classes on the union of $\{0\}$ and semisimple elements satisfying $(i)-(iii)$.

Let $\gss$ be a semisimple element as above or $\gss=0$.
Then we define $\lieG_\gss:=\,^\G(\gss+\lieG^0_0)$ where 
$\lieG^0$ is the Lie algebra of $\G^0=\C_\G(\gss)$ 
(recall that $\lieG^0_0=\cup_{\xx\in\Bd(\bG^0,k)}\lieG^0_{\xx,0}$). 
Each $\lieG_\gss$ is a $\G$-domain, 
an open and closed $\G$-invariant subset of $\lieG$.
We say that two such semisimple elements $\gss$ and $\gss'$ are equivalent 
if $\lieG_\gss=\lieG_{\gss'}$ (see (\ref{def: ss equivalence})).
Let $\sgsc$ be the set of equivalence classes of $\gss$'s.
Then, $\lieG$ is the disjoint union of $\G$-domains
$\lieG_\gss$, $\gss\in\sgsc$ (see (\ref{prop: dec dualg})), and we have 
\[
\int_{\lieG}\widehat f(X)\,dX
=\sum_{\gss\in\sgsc}\int_{\lieG_\gss}\widehat f(X)\,dX\ .
\]
On the other hand, each $\gss\in\sgsc$ also parameterizes a subset 
$\admG_\gss$ of $\admG$. Roughly speaking, $\admG_\gss$ is the subset 
of $\admG$ which consists of $(\pi,V_\pi)\in\admG$ containing 
$(\vec\G_{\xo,\vec\stm^+},\chi_{\gss'})$ for some $\gss'\sim\gss$ and 
$\xo\in\Bd(\bC_{\bG}(\gss'),k)$, and
$\admG_0$ is the set of depth zero representations
(see (\ref{def: admGamma}) and (\ref{rmk: indep epsilon}) for details).
Moreover, $\admG_\gss=\admG_{\gss'}$ if $\gss$ and $\gss'$
are equivalent, and $\admG_\gss\cap\admG_{\gss'}=\emptyset$ otherwise
(see (\ref{lem: dec admG})).
Setting $\temG_{\sgsc}
:=\left(\overset\circ\cup_{\gss\in\sgsc}\temG_\gss\right)$ where
$\temG_\gss=\temG\cap\admG_\gss$, we have
\[
\int_{\temG}\charpi(f\circ\log)\,d\pi
=\sum_{\gss\in\sgsc}\int_{\temG_\gss}\charpi(f\circ\log)\,d\pi\ 
+\int_{\temG\setminus\temG_{\sgsc}}\charpi(f\circ\log)\,d\pi.
\]
It is obvious that $\temG_{\sgsc}\,\subset\,\temG$. Our claim above 
is that $\temG\setminus\temG_{\sgsc}$ has Plancherel measure zero. 

Now, we match terms parameterized by $\gss\in\sgsc$ (see \S\ref{sec: match}):
\begin{equation}
\label{eq: main}
\int_{\lieG_\gss}\Ff(X)\,dX
=\int_{\temG_\gss} \charpi(f\circ\log)\,d\pi\ 
\end{equation}
for $f\in\Ccs(\lieG_{0^+})$ if $\gss\sim0$,
and for $f\in\Ccs(\lieG_{\stm_{d-1}^+})$ if $\gss\not\sim0$.
A similar equality was considered in \cite{Howe2} when $\gss$ is regular, and
in \cite{Dpl1} when $\gss$ is a good element. 
If $\gss\sim0$, this is already proven in \cite{Dpl1}.
If $\gss\not\sim0$, we regard both sides of (\ref{eq: main}) 
as distributions on $\Ccs(\lieG_{\stm_{d-1}^+})$, and denote the distributions
on the left and the right side of (\ref{eq: main}) by $\rT_\ell$ and $\rT_r$
respectively.
We first need to prove that the $\distr^{-\gss}$ from \cite{Asym2}
(see also the previous section) 
contain both distributions $\widehat\rT_\ell$ and $\widehat\rT_r$. Then, 
by matching $\widehat\rT_\ell(f_\orb)=\widehat\rT_r(f_\orb)$ for each
test function $f_{\orb}$, $\orb\in\orb(-\gss)$, found in \cite{Asym2},
we verify that they are equal on $\Ccs(\lieG_{\stm_{d-1}}^+)$. 
Since we prove this equality by matching them only on test functions, 
and these test functions have the property that
$\{\pi\in\temG\mid\charpi(\Ff_{\orb})\neq0,
\textrm{ for some }\orb\in\orb(-\gss)\}$ is a subset of $\temG_\gss$,
we do not need any explicit knowledge of the Plancherel measure $d\pi$.

Using this, we can also prove that the equality in (\ref{eq: main})
holds for any characteristic function $f_{x,\stm}$
of a lattice $\lieG_{x,\stm}$ with $x\in\Bd(\bG,k)$ and $\stm>0$.
Summing over all $\gss\in\sgsc$, we have 
\begin{equation*}
\int_\lieG\widehat f_{x,\stm}(X)\,dX
=\sum_{\gss\in\sgsc}\int_{\lieG_\gss}\widehat f_{x,\stm}(X)\,dX
=\sum_{\gss\in\sgsc}\int_{\temG_\gss}\charpi(f_{x,\stm}\circ\log)\,d\pi
=\int_{\temG}\charpi(f_{x,\stm}\circ\log)\,d\pi\ ,
\end{equation*}
which will lead to a proof that almost every irreducible 
tempered representation (in particular, supercuspidal representation)
is an element of $\temG_\gss$ for some $\gss\in\sgsc$ 
(Theorems \ref{thm: tem exhaustion} and \ref{thm: tempered}).

\subsection{\sc Generic $\G$-datum}

The construction in \cite{yu} is based on 
a generic $\G$-datum which consists of a quintuple
$\datum_Y=(\vec\bG,\xo,\vec\rtm,\vec\phi,\rho)$ 
satisfying the following five conditions 
(see \cite[\S3]{yu} or \S\ref{sec: sc}):

$\mathbf D1.$ $\vec\bG=(\bG^0,\bG^1,\cdots,\bG^d=\bG)$ is a tamely ramified
Levi sequence and $\bZ_{\bG^0}/\bZ_{\bG}$ is $k$-anisotropic, where
$\bZ_{\bG^0}$ (resp. $\bZ_{\bG}$) is the center of $\bG^0$ (resp. $\bG$).

$\mathbf D2.$ $\xo\in\Bd(\bG^0,k)$.

$\mathbf D3.$ $\vec\rtm=(\rtm_0,\rtm_1,\cdots,\rtm_{d-1},\rtm_d)$ is
a sequence of positive real numbers 
with $0<\rtm_0<\cdots<\rtm_{d-2}< \rtm_{d-1}\le\rtm_d$ if $d>0$,
$0\le\rtm_0$ if $d=0$.

$\mathbf D4.$ $\vec\phi=(\phi_0,\phi_1,\cdots,\phi_d)$ is a sequence of 
quasi-characters; $\phi_i$, $0\le i\le d-1$, 
is a $\G^{i+1}$- generic character 
of $\G^i$ of depth $\rtm_i$ at $\xo$. If $\rtm_{d-1}<\rtm_d$, 
$\phi_d$ is a generic character of $\G^d$ of depth $\rtm_d$ at $\xo$, and
$\phi_d$ is trivial otherwise.

$\mathbf D5.$ $\rho$ is an irreducible representation of $\G^0_{[\xo]}$,
the stabilizer in $\G^0$ of the image $[\xo]$ of $\xo$ 
in the reduced building of $\bG^0$, such that $\rho\mid\G^0_{\xo,0^+}$ 
is a multiple of the trivial representation and
$c\textrm{-Ind}_{\G^0_{[\xo]}}^{\G^0}(\rho)$ is 
irreducible and supercuspidal.

Based on the above data, Yu constructs a pair $(K_{\datum_Y},\rho_{\datum_Y})$ 
of an open compact modulo center subgroup $K_{\datum_Y}$ and 
its irreducible representation $\rho_{\datum_Y}$ such that 
$c\textrm{-Ind}_{K_{\datum_Y}}^\G\rho_{\datum_Y}$ is supercuspidal. 
Denote the resulting supercuspidal representation 
$c\textrm{-Ind}_{K_{\datum_Y}}^\G\rho_{\datum_Y}$ by $\pi_{\datum_Y}$.

\subsection{\sc Second step}
Now, let $(\pi,V_\pi)$ be a supercuspidal representation of $\G$. We want to 
prove that there is a generic $\G$-datum $\datum_Y$ such that 
$\pi\simeq\pi_{\datum_Y}$.
Since the case of depth zero supercuspidal representations is
already known (see \cite{Mr}, \cite[(6.6), (6.8)]{MP2}), we may assume that 
$\pi$ is of positive depth.
Then from the first step, $\pi$ contains a $\bK$-type 
$(\vec\G_{\xo,\vec\stm^+},\chi_\gss)$ for some $(\gss,\xo)$.
Let $\gss_i$, $\bG^i$ and $\rtm_i$ be as in $(i)$-$(iii)$.
Let $\phi_i$ be a quasi-character of $\G^i$
extending the character $\chi_{\gss_i}$ of $\G^i_{\xo,\stm_i^+}$ defined
by $\gss_i$. Under our hypothesis, $\chi_{\gss_i}$ always extends 
to a quasi-character of $\G^i$ (see (\ref{lem: sstochar})). 
Let $\rho$ be an irreducible component of the $\G^0_{[\xo]}$-representation 
$\phi^{-1}\otimes(\pi|V_\pi^{\chi_\gss})$ 
where $\phi=\prod_i(\phi_i|\G^0_{[y]})$ and 
$V_\pi^{\chi_\gss}$ is the $\chi_\gss$ isotypic component in $V_\pi$.
To show that $(\vec\bG,\xo,\vec\rtm,\vec\phi,\rho)$ 
satisfies the desired properties $\mathbf D1$-- $\mathbf D5$, 
we need to verify that

\smallskip

(a) $\bZ_{\bG^0}/\bZ_{\bG}$ is anisotropic,

(b) $c\textrm{-Ind}_{\G^0_{[\xo]}}^{\G^0}(\rho)$ is irreducible 
and supercuspidal.

\smallskip

We prove (a) by analyzing appropriate Hecke algebras (\S13-\S14), 
and (b) by analyzing appropriate Jacquet modules (\S15-\S17). 

Let $\datum_Y$ be the generic $\G$-datum $(\vec\G,\xo,\vec\rtm,\vec\phi,\rho)$ 
associated to $(\pi,V_\pi)$ found as above.
Lastly, we show $\chi_\gss$  on $\vec\G_{\xo,\vec\stm^+}$
can be extended further to $\rho_{\datum_Y}$ on $K_{\datum_Y}$ (\S18), and
$(\pi,V_\pi)$ still contains $(K_{\datum_Y},\rho_{\datum_Y})$. Then,
by Frobenius reciprocity, we can conclude that $\pi$ is 
in fact the supercuspidal representation constructed from $\datum_Y$, that is,
$\pi\simeq\pi_{\datum_Y}$.

\smallskip

Generally speaking, different generic $\G$-datums can yield isomorphic 
supercuspidal representations (this is the case e.g.\ for $\G$-conjugate
generic $\G$-datums). The question of when exactly the resulting supercuspidal 
representations are isomorphic is settled by recent work of Hakim and 
Murnaghan \cite{HaMu}; we will not discuss it here.

\smallskip

If $\G$ is one of the classical groups considered in \cite{Kim},
the corresponding datum in \cite{Kim} is $(\gss,\G^0_{[\xo]},\rho)$.
The result of this paper also implies that the $\bK$-types
constructed in \cite{Kim} form a complete set for the classical groups
considered in that paper.

\smallskip

In the first three sections,
we review some facts about Moy-Prasad filtrations and relevant results
on buildings. Otherwise, reviews of many necessary results 
(in particular from \cite{Asym2}) are spread
throughout this paper before they are used.
The first step is carried out in \S4-\S11. We review
Yu's construction of supercuspidal representations in \S12.
The second step is done in \S13-\S17.
In \S18, we compare $(\vec\G_{\xo,\stm^+},\chi_\gss)$ and 
$(K_{\datum_Y},\rho_{\datum_Y})$. Finally in \S19, we conclude that 
all supercuspidal representations arise through Yu's construction
(Theorem \ref{thm: Main}).

Most notation is used throughout the paper once it is defined.
The table of some selected notation is available in the end of this paper.

\medskip

\noindent{\it Acknowledgments.} 
I thank Fiona Murnaghan for the joint work on character expansions,
which has been an important tool for this work.
Part of this paper is motivated from a conversation with Roger Howe
to whom I thank for the inspiration.
I would like to thank Joseph Bernstein, Allen Moy, Fiona Murnaghan 
and Jui-Kang Yu for many enlightening discussions, and Anne-Marie Aubert, 
Stephen DeBacker, Jonathan Korman and Loren Spice for their helpful comments 
on the earlier version of this paper. 
Finally, I would like to thank Tom Hales and Paul Sally for
their interest and warm support of this work.

\section*{\bf Notation and Conventions}\label{sec: notation}

Let $k$ be a $p$-adic field (a finite extension of $\bbQ_p$) 
with residue field $\bbF_{p^n}$.
Let $\val=\val_k$ be the valuation on $k$ such that $\val(k^\times)=\bbZ$.
Let $\overline k$ be an algebraic closure of $k$.
For an extension field $E$ of $k$, let $\val_E$ be the valuation
on $E$ extending $\val$. We will just write $\val$ for $\val_E$.
Let $\CaO_E$ be the ring of integers of $E$ with prime ideal $\pid_E$.

Let $\addch$ be a fixed additive character of $k$ such that
$\addch|{\CaO_k}\ne 1$ and $\addch|{\pid_k}=1$.

Let $\bG$ be a connected reductive group defined over $k$, and
$\blieG$ the Lie algebra of $\bG$.
Denote the group of $E$-rational points of $\bG$ by $\bG(E)$, and
the Lie algebra of $E$-rational points of $\blieG$ by $\blieG(E)$.
We denote $\bG(k)$ and $\blieG(k)$ by $\G$ and $\lieG$ respectively. 
Similarly, the linear duals of $\blieG$ and $\blieG(E)$ are denoted by
$\blieG^\ast$ and $\blieG^\ast(E)$ respectively.
We write $\lieG^\ast$ for $\blieG^\ast(k)$.
Let $\bZ_\bG$ denote the center of $\bG$, and $\blieZ_{\blieG}$ 
the Lie algebra of $\bZ_\bG$.
Let $\bG^\der$ denote the derived group of $\bG$, and $\blieG^\der$
the Lie algebra of $\bG^\der$.
In general, we use bold characters $\bH,\,\bM,\,\bN$, etc
to denote algebraic groups, and $\boldsymbol{\mathfrak h}$, 
$\boldsymbol{\mathfrak m}$, $\boldsymbol{\mathfrak n}$ 
to denote their Lie algebras. If they are defined over $k$, we will use
the corresponding roman characters $H,\,M$ and $N$ to denote the groups of
$k$-points, and $\mathfrak h,\,\lieM$ and $\lieN$ to denote
the Lie algebras of $H,\,M$ and $N$.

Let $\nil$ denote the set of nilpotent elements in $\lieG$.
There are different notions of nilpotency. However, since we assume
that $\textup{char}(k)=0$, those notions are all the same. 
We refer to \cite{dBP, MP} for more discussion of this.

If $X$ is a topological space with a Borel measure $dx$ and
if $Y$ is a Borel subset of $X$, $\vol_X(Y)$ denotes
the volume of $Y$ with respect to $dx$.

For any given set $W$, let $|W|$ denote the cardinality of $W$.

For any subset $S$ in $\lieG$ or in $\G$, 
we denote by $[S]$ the characteristic function on $S$,
and by $-S$ the set $\{-s\mid s\in S\}$.
For $g\in G$, ${ }^g\!Z$ denotes $gZg^{-1}$.

Let $\tbR:=\bbR\cup\{\rtm^+\mid\rtm\in\bbR\}$. We define an ordering
on $\tbR$ extending the one on $\bbR$: let $\rtm,\stm\in\bbR$.
Then, $\rtm<\rtm^+$. 
If $\rtm<\stm$, then $\rtm<\stm^+$, $\rtm^+<\stm^+$ and $\rtm^+<\stm$.
We define an addition on $\tbR$ extending the one on $\bbR$: 
for $\rtm,\stm\in\bbR$, $\rtm^++\stm=(\rtm+\stm)^+=\rtm^++\stm^+$.
Set $(\rtm^+)^+:=\rtm^+$.

Finally, we will not distinguish between representations and
their isomorphism classes.

\medskip

\section{\bf Moy-Prasad filtrations}

\begin{numbering} {\bf Apartments and buildings.}\rm \ 
For a finite extension $E$ of $k$, 
let $\Bd(\bG, E)$ denote the extended Bruhat-Tits building of $\bG$ over $E$. 
Recall that 
$\Bd(\bG, E)
=\Bd(\mathbf{DG},E)\times\left(\bX_\ast(\mathbf{Z_G},E)\otimes\bbR\right)$,
where $\mathbf{DG}$ is the derived group of $\bG$, and
$\bX_\ast(\mathbf{Z_G},E)$ is the abelian group 
of $E$-rational cocharacters of
the center $\mathbf{Z_G}$ of $\bG$.
For a maximal $E$-split torus $\bT$ in $\bG$, 
let $\Apt(\bT,E)$ be the corresponding apartment over $E$. 
It is known that for any tamely ramified
finite Galois extension $E'$ of $E$, 
$\Bd(\bG,E)$ can be embedded into $\Bd(\bG,E')$
and its image is equal to the set of the Galois fixed points in $\Bd(\bG,E')$
(see \cite[(5.11)]{Rous} or \cite{Pr}).

\end{numbering}

\begin{numbering} {\bf  Moy-Prasad filtrations.}  
\label{num: MP filtrations}\rm
Let $(x,r)\in\Bd(\bG,E)\times\bbR$. 
Regarding $\bG$ as a group defined over $E$, 
Moy and Prasad define $\blieG(E)_{x,r}$ and also $\bG(E)_{x,r}$ if 
$\rtm\ge0$ with respect to the valuation normalized as follows (\cite{MP2}):
let $E^u$ be the maximal unramified extension of $E$,
and $L$ the minimal extension of $E^u$ over which $\bG$ splits.
Then the valuation used by Moy and Prasad maps $L^\times$ onto $\bbZ$.

In a similar way, with respect to our normalized valuation $\nu$,
we can define filtrations in $\blieG(E)$ and $\bG(E)$.
Then our $\blieG(E)_{x,r}$ and $\bG(E)_{x,r}$ correspond to their 
$\blieG(E)_{x,elr}$ and $\bG(E)_{x,elr}$ respectively,
where $e=e(E/k)$ is the ramification index of $E$ over $k$ and
$\ell=[L:E^u]$. Hence, if $\varpi_{_E}$ is a uniformizing element of $E$,
our filtrations satisfy 
$\varpi_{_E}\blieG(E)_{x,\rtm}=\blieG(E)_{x,\rtm+\frac1e}$
while theirs satisfy
$\varpi_{_E}\blieG(E)_{x,\rtm}=\blieG(E)_{x,\rtm+\ell}$.

This normalization is chosen to have the following property 
(see also \cite[(1.4.1)]{Ad}): 
for a tamely ramified Galois extension $E'$ of $E$ 
and $x\in\Bd(\bG,E)\subset\Bd(\bG,E')$, we have
\[
\blieG(E)_{x,r}=\blieG(E')_{x,r}\cap\blieG(E).
\]
If $r>0$, we also have
\[
\bG(E)_{x,r}=\bG(E')_{x,r}\cap\bG(E).
\]

\smallskip

For simplicity, we put $\lieG_{x,\rtm}:=\blieG(k)_{x,\rtm}$, etc,
and $\Bd(\G):=\Bd(\bG,k)$. 
We will also use the following notation. Let $\rtm\in\bbR$.

\begin{enumerate}
\item 
$\lieG_{x,r^+}=\cup_{s>r} \lieG_{x,s}$ and
$\G_{x,|r|^+}=\cup_{s>|r|} \G_{x,s}$, \ $x\in\Bd(\G)$.
\item 
$\lieG^\ast_{x,r}=\left\{\chi\in\lieG^\ast
\mid\chi(\lieG_{x,(-r)^+})\subset\pid_k\right\}$, \ $x\in\Bd(\G)$.
\item
$\lieG_r=\cup_{x\in\Bd(\G)} \lieG_{x,r}$ and
$\lieG_{r^+}=\cup_{s>r} \lieG_s$
\item 
$\G_r=\cup_{x\in\Bd(\G)} \G_{x,r}$ and
$\G_{r^+}=\cup_{s>r} \G_s$ for $r\ge0$.
\end{enumerate}

The hypothesis $\HypB$ in \S\ref{num: hypotheses} 
is concerned with identifying $\lieG^\ast_{x,r}$
with $\lieG_{x,r}$ via an appropriate bilinear form $\bilinear$ on $\lieG$
(see \cite[(4.1)]{AR}). 

\end{numbering}

\begin{numbering} {\bf Root decomposition.}\label{num: root dec}
\rm \ 
Let $\bT$ be a maximal $k$-torus in $\bG$, and $E$ a finite extension
of $k$ over which $\bT$ splits.
Let $\Phi(\bG,\bT,E)$ be the set of $E$-roots of $\bT$ in $\bG$,
and let $\Psi(\bG,\bT,E)$ be the corresponding set of affine roots in $\bG$.
If $\psi\in\Psi(\bG,\bT,E)$, 
let $\dot\psi\in\Phi(\bG,\bT,E)$ be the gradient of $\psi$, 
and $\blieG(E)_{\dot\psi}\subset\blieG(E)$ the root space 
corresponding to $\dot\psi$.
We denote the open compact abelian group in $\blieG(E)_{\dot\psi}$ 
corresponding to $\psi$ by $\blieG(E)_\psi$ (\cite[(3.2)]{MP}).

Let $\bX_\ast(\bT,E)$ be the set of cocharacters of $\bT$, 
and let $\bX^\ast(\bT,E)$ be the set of characters of $\bT$. 
For $\rtm\in\tbR$, let 
\[
\blieT(E)_r
=\{\gss\in\blieT(E)\mid\nu(d\chi(\gss))\ge r\textrm{ for all }
\chi\in\bX^\ast(\bT,E)\}.
\]
Then, for $x\in\Apt(\bG,\bT,E)$, we have
\[
\blieG(E)_{x,r}
=\blieT(E)_r+\sum_{\psi\in\Psi(\bG,\bT,E),\ \psi(x)\ge r} \blieG(E)_\psi.
\]

Let $\bT$ be a maximal $k$-torus in $\bG$ which splits 
over a tamely ramified finite Galois extension $E$ of $k$.
Then, we write $\Apt(\bG,\bT,k)$ for $\Apt(\bG,\bT,E)\cap\Bd(\bG,k)$. This 
definition is independent of the choice of $E$ \cite{yu}. Moreover,
$\Apt(\bG,\bT,k)$ is the set of Galois fixed points in $\Apt(\bG,\bT,E)$.
\end{numbering}

\midvsp

\section{\bf Twisted Levi sequences}

\begin{definition} \label{def: levi seq} \cite{yu} \rm
Let $\bG$ be a connected reductive $k$-group.
Let $\vec\bG:=(\bG^0,\cdots,\bG^d=\bG)$ be a sequence of
connected reductive $k$-groups with
$\bG^0\subsetneq\bG^1\subsetneq\cdots\subsetneq\bG^d$.

\begin{enumerate}
\item
If each $\bG^i$ is a $k$-split Levi subgroup of $\bG$, 
$\vec\bG$ is called a {\it $k$-Levi sequence} in $\bG$.
\item
If there exists a (tamely ramified) finite extension $E/k$ 
such that $\bG^0\otimes E$ is split, 
and $\vec\bG\otimes E=(\bG^0\otimes E,\cdots,\bG^d\otimes E)$ 
is an $E$-Levi sequence in $\bG^d\otimes E$, then
$\vec\bG$ is called a {\it (tamely ramified) twisted Levi sequence in $\bG$}
and $E$ is called a {\it splitting field} of $\vec\bG$.
\end{enumerate}
\end{definition}

\smallskip

Note that any subsequence of a tamely ramified twisted Levi sequence
$\vec\bG$ is also a tamely ramified twisted Levi sequence.

\begin{numbering}\label{num: G'G}\rm
Let $(\bH,\bG)$ be a tamely ramified twisted Levi sequence, and
$E$ a tamely ramified Galois extension $E$ of $k$ over which
$(\bH,\bG)$ splits.
Since $\bH(E)$ is a Levi subgroup of $\bG(E)$, 
there is a Galois equivariant embedding of $\Bd(\bH,E)$ 
into $\Bd(\bG,E)$, which in turn induces an embedding of $\Bd(\bH,k)$ 
into $\Bd(\bG,k)$ (see \cite[\S1.9]{Ad} or \cite[(2.11)]{yu}). 
Such embeddings are unique modulo translation 
by $\bX_\ast(\bZ_{\bH},k)\otimes\bbR$. However, the images remain the same.

Fix an embedding $i:\Bd(\bH,k)\longrightarrow\Bd(\bG,k)$. Then
we will regard $\Bd(\bH,k)$ as a subset of $\Bd(\bG,k)$ and 
write simply $x$ for $i(x)$.
For any $x\in\Bd(\bH,k)$, the associated filtrations on $\rH:=\bH(k)$ 
and $\lieH:=\blieH(k)$ satisfy the following (\cite[(1.9.1)]{Ad}):
\begin{align*}
\rH_{x,r}&=\bG(E)_{x,r}\cap\rH=\G_{x,r}\cap\rH\quad\text{for }r>0,\\
\lieH_{x,r}&=\blieG(E)_{x,r}\cap\lieH=\lieG_{x,r}\cap\lieH
\quad\text{for any }r\in\tbR.
\end{align*}
For a $k$-torus $\bT\subset\bH$ which splits over $E$, 
we have $\Apt(\bH,\bT,k)=\Apt(\bG,\bT,E)\cap\Bd(\bH,k)$.
\end{numbering}

\begin{numbering}\label{num: embeddings}\rm
If $\vec\bG=(\bG^0,\cdots,\bG^d=\bG)$ is a tamely ramified twisted
Levi sequence, we can and will fix a sequence of embeddings
\[
\Bd(\bG^0,k)\hookrightarrow\Bd(\bG^1,k)\hookrightarrow\Bd(\bG^2,k)
\hookrightarrow\cdots\hookrightarrow\Bd(\bG^d,k).
\]
We will identify $\Bd(\bG^i,k)$ with a subset of $\Bd(\bG^j,k)$ for 
$i\le j$. Moreover, for $x\in\Bd(\bG^i,k)$ and $\rtm\in\tbR$, we have
\begin{align*}
\lieG^i_{x,r}&=\blieG^j(E)_{x,r}\cap\lieG^i=\lieG^j_{x,r}\cap\lieG^i,\\
\G^i_{x,r}&=\bG^j(E)_{x,r}\cap\G^i=\G^j_{x,r}\cap\G^i\quad\text{if }r>0.
\end{align*}
\end{numbering}

From now on, we say that a semisimple element $\gss\in\lieG$ 
\emph{splits} over a finite extension $E$ 
if $\gss$ lies on a $k$-torus which splits over $E$.
 
\begin{lemma}\label{lem: levi seq}
Let $\gss\in\lieG$ be a semisimple element
which splits over a tamely ramified finite extension $E$ of $k$. 
Set $\bH:=\bC_{\bG}(\gss)$, the centralizer of $\gss$ in $\bG$. 
Then $(\bH,\bG)$ is an $E$-split tamely ramified twisted Levi sequence.
\end{lemma}

\proof
Without loss of generality, we may assume that $\gss$ is
in a $k$-split torus $\blieT$. 
By (7.1) and (7.2) of \cite{yu}, $\bH$ is connected and reductive.
Since $\blieH$ is reductive, 
$\blieZ_{\blieH}\subset\blieT$, and thus
$\blieZ_{\blieH}$ is a $k$-split subtorus of $\blieT$.
Combining this with the fact that 
$\bH$ is the centralizer of $\blieZ_{\blieH}$ in $\bG$,
we conclude that $\bH$ is a $k$-Levi subgroup of $\bG$.
\qed

\begin{remarks}\rm \ 
\begin{enumerate}
\item
Note that if $\gss=0$, then $\bH=\bG$. 
\item
The above lemma is not valid if $\gss$ is replaced by a semisimple
element of $\G$. For example, in $\mathbf{Sp}_4$, the centralizer of
a semisimple group element can be $\mathbf{SL}_2\times\mathbf{SL}_2$.
\end{enumerate}
\end{remarks}

\begin{numbering}\label{num: G'G2}\rm
Let $(\bH,\bG)$ be a tamely ramified twisted Levi sequence. 
For $X\in\lieG$, denote the $\rH$-orbit ${}^{\rH}\!X$ of $X$ by 
$\orb^{\rH}_X$. For simplicity, we write $\orb_X$ for $\orb^{\G}_X$.
In general, we use the notation $\orb^\rH$ to denote $\rH$-orbits.
If $X\in\lieH$, $\orb({\rH},X)$ denotes
the set of all $\rH$-orbits whose closure in $\lieH$ contains $X$.
We write $\orb(X):=\orb(\G,X)$.

When $X\in\lieG$ is semisimple, $\orb(X)$ is described in \cite[\S 2]{HC}.
In the situation of Lemma \ref{lem: levi seq}, 
write $\gss=X$ and $\rH=C_\G(\gss)$, then
\[
\orb(\gss)=\{\orb_{\gss+n}\mid n\in \orb^{\rH}\in\orb(\rH,0)\}.
\]
Note that $\gss+n$ is already in the form of a Jordan decomposition.
The map $\orb_{\gss+n}\rightarrow\orb_n$, $n\in\orb(\rH,0)$
induces a bijection of $\orb(\gss)$ with $\orb({\rH},0)$,
the set of nilpotent $\rH$-orbits in $\lieH$.
\end{numbering}

\section{\bf Admissible sequences and lattices in $\lieG$ and $\lieG^\ast$}

We recall some definitions from \cite{BT1} and \cite{yu}.

\begin{definition} \rm
Let $E/k$ be a tamely ramified extension, and let $\bT$ be an
$E$-split maximal $k$-torus in $\bG$. Let $\Phi=\Phi(\bG,\bT,E)$ be
the corresponding root system.
Then, a function $f:\Phi\cup\{0\}\longrightarrow\tilde\bbR$ is
{\it concave} if for every non-empty finite subset 
$\{a_i\}\subset\Phi\cup\{0\}$
such that $\sum a_i\in\Phi\cup\{0\}$, we have
\[
f(\sum_i a_i)\,\le\,\sum_i f(a_i).
\]
\end{definition}

We keep the notation from the above definition.
If $\xx\in\Apt(\bG,\bT,E)$ and 
$f$ is a concave function on $\Phi(\bG,\bT,E)\cup\{0\}$,
there are a group $\bG(E)_{\xx,f}$ and a lattice $\blieG(E)_{\xx,f}$
associated to $\xx$ and $f$. (see \cite{AD2, yu} for details). 
If $f$ is $\Gal(E/k)$-invariant and 
$\xx\in\Apt(\bG,\bT,E)^{\Gal(E/k)}$,
we can define
\begin{align*}
\G_{\xx,f}&:=\bG(E)_{\xx,f}^{\Gal(E/k)}\cap\G_{\xx,0}, \\
\lieG_{\xx,f}&:=\blieG(E)_{\xx,f}^{\Gal(E/k)}.
\end{align*}
If $f$ is positive and $E/k$ is tame, then
$\G_{\xx,f}=\bG(E)_{\xx,f}^{\Gal(E/k)}.$

\begin{definition}\rm 
A sequence $\vec\utm:=(\utm_0,\cdots,\utm_d)$ 
in $\tilde\bbR$ is {\it admissible} if for some $0\le c\le d$, we have
\[
0\le\utm_0=\utm_1=\cdots=\utm_c\qquad\textrm{and}\qquad
\frac 12\utm_c\le\utm_{c+1}\le\cdots\le\utm_d.
\]
\end{definition}

\begin{numbering}\label{num: vecG}\rm
Let $E$ be a tamely ramified Galois extension of $k$.
Let $\vec\bG=(\bG^0,\cdots,\bG^d=\bG)$ be an $E$-split twisted Levi sequence, 
and let $\vec\utm=(\utm_0,\cdots,\utm_d)$ be
an admissible sequence. Then we can associate to $\vec\utm$ 
a $\Gal(E/k)$-invariant, concave function $f_{\vec\utm}$ on $\Phi(\bG,\bT,E)$
as follows: for an $E$-split maximal $k$-torus $\bT\subseteq\bG^0$, define
\[
f_{\vec\utm}(a)
=\left\{
\begin{array}{ll}
\utm_0\qquad&\textrm{if }\atm\in\Phi(\bG^0,\bT,E)\cup\{0\},\\
\utm_i\qquad&\textrm{if }\atm\in\Phi(\bG^i,\bT,E)
   \setminus\Phi(\bG^{i-1},\bT,E)\textrm{ for }i>0.
\end{array}\right.
\]
For $\xx\in\Apt(\bG,\bT,k)$, let
\begin{align*}
\vec\G_{\xx,\vec\utm}\ :=\ \G_{\xx,f_{\vec\utm}}\,,\qquad\qquad
\vec\lieG_{\xx,\vec\utm}\ :=\ \lieG_{\xx,f_{\vec\utm}}\,.
\end{align*}
These are well defined independent of the choice of $\bT$
with $\xx\in\Apt(\bG,\bT,k)$.
Moreover, if $\vec\utm$ is nondecreasing, we have
\[
\vec\G_{\xx,\vec\utm}
=\G^0_{\xx,\utm_0}\G^1_{\xx,\utm_1}\cdots\G^d_{\xx,\utm_d}\,,
\qquad\qquad
\vec\lieG_{\xx,\vec\utm}=\lieG^0_{\xx,\utm_0}+\lieG^1_{\xx,\utm_1}
+\cdots+\lieG^d_{\xx,\utm_d}\,.
\]
\end{numbering}

\begin{numbering} {\bf Hypotheses.}\label{num: hypotheses} \rm
We list the hypotheses used in this paper.
They are labeled by $\Hypk$, $\HypB$, $\HypGT$ and $\HypN$ respectively.
We will state explicitly whenever these hypotheses are necessary.

\medskip

\noindent
$\Hypk$
The residue characteristic $p$ is large enough (depending on
$\bG$ and $\val(p)$) such that the following hold.
\begin{enumerate}
\item 
The exponential map (resp. the logarithmic map) is defined 
on the subset $\lieG_{0^+}$ of $\lieG$ (resp. $\G_{0^+}$ of $\G$), 
and for a tamely ramified finite Galois extension $E$ over $k$, an $E$-split
maximal $k$-torus $\bT$ of $\bG$, $\xx\in\Apt(\bG,\bT,k)$,
and a $\Gal(E/k)$-invariant concave function $f$ 
on $\Phi(\bG,\bT,E)\cup\{0\}$ with $f(0)>0$, 
we have $\exp(\lieG_{x,f})=\G_{x,f}$ 
(resp. $\log(\G_{x,f})=\lieG_{x,f}$).
\item
For any $x\in\Bd(\bG,k)$, $X\in\lieG_{x,0^+}$ and $Y\in\lieG_{x,0}$, 
$\frac1p(\ad X)^{p-1}(Y)\in\lieG_{x,0^+}$.
\end{enumerate}

\bigskip

\noindent
$\HypB$ 
$\bG$ satisfies the condition in Proposition 4.1 in \cite{AR}.

\bigskip

\noindent
$\HypGT$
Every maximal $k$-torus $\bT$ in $\bG$ splits over a tamely ramified
Galois extension, and for any $r\in\bbR$, any nontrivial coset in $\lieT_r$
modulo $\lieT_{r^+}$ contains a good element (as defined in (\ref{def: gss}))
of depth $r$.

\bigskip

\noindent
$\HypN$
For any tamely ramified twisted Levi subgroup $\bH$ of $\bG$,
the hypotheses in \cite[\S4.2]{dBP} are valid.

\end{numbering}

\begin{remark}\rm
The condition in Proposition 4.1 of \cite{AR} requires that either
$\bG$ is a form of $\textrm{GL}_n$, or 
the absolute Dynkin diagram of $\bG$ has no bonds of order $p$ 
and $p$ does not divide $2\mathrm{k}(\bG)|\pi_1(\bG')|$
(see \cite{AR}, Proposition 4.1 for notation).
One sees easily that this is satisfied if $p$ is large enough.

One can use the Campbell-Hausdorff formula to determine
a sufficient condition on $k$ for $\Hypk$ to hold
(see \cite{Kim}, Proposition 3.1.1). 

In \cite{dBP}, under some hypotheses on $\bG$ and $k$ (see \cite{dBP}, \S4.2),
DeBacker gives a parameterization of nilpotent orbits in $\lieG$
via Bruhat-Tits theory (see Theorem 5.6.1 of \cite{dBP}). 
He uses this parameterization to get a homogeneity result in \cite{dBH}.
We need the hypothesis $\HypN$ to use the results of 
\cite{dBH} and \cite{dBP}. We refer the reader to \cite{dBP}
for precise statements. Again, if $p$ is large enough, $\HypN$ is valid.
\end{remark}

\begin{remark}\rm
If $\HypB$ is satisfied, there is a $\overline k$-valued, nondegenerate, 
$\bG(\overline k)$-invariant, symmetric, bilinear form $\bilinear$ 
on $\blieG$ satisfying the following (see the proof of \cite[(4.1)]{AR}): 
for any tamely ramified finite extension $E$ of $k$, 
$\bilinear$ induces an $E$-valued $E$-bilinear form on $\blieG(E)$ such that
\begin{enumerate}
\item
we can identify $\blieG^\ast(E)_{x,r}$ with $\blieG(E)_{x,r}$ via the map
$\Omega:\blieG(E)\to\blieG^\ast(E)$ defined by
$\Omega(X)(Y)=\bilinear(X,Y)$,
\item
if $\bT$ is a maximal $E$-split torus and
$(\bT,\{X_\alpha\})$ is a Chevalley splitting, then
we have 

({\it i}) $\bilinear(X_\alpha,X_\beta)=0$ unless $\alpha+\beta=0$,

({\it ii}) $\bilinear(X_\alpha,X_{-\alpha})\ne0$, and has valuation $0$,

({\it iii}) $\blieG=\blieZ_{\blieG}\oplus\blieG^\der$ 
is an orthogonal decomposition with respect to $\bilinear$. 
\item
For any Levi subgroup $\bM$ of $\bG$ which splits 
over a tamely ramified extension,
$\bilinear|\lieM\times\lieM$ satisfies (1) and (2).
\end{enumerate}

The above $\HypB$ implies the corresponding hypothesis 
$\HypB$ (labeled in the same way) in \cite{Asym2}.
Whenever we assume $\HypB$, we denote the associated bilinear form 
by $\bilinear$.
For more discussion on sufficient conditions for the above hypotheses, 
we refer to \cite{Asym2}.
\end{remark}

\begin{numbering}\label{num: orth dec}\rm
Let $\vec\utm=(\utm_0,\cdots,\utm_d)$ be a sequence of real numbers
which is not necessarily admissible. We can still define 
$f_{\vec\utm}$ and a subset 
$\vec\lieG_{\xx,\vec\utm}=\lieG_{\xx,f_{\vec\utm}}$ of $\lieG$
in a similar fashion as in (\ref{num: vecG}). 
Suppose $\HypB$ is valid.

For each $i=0,\cdots,d-1$, let $\lieG^{i\,\perp}$ 
be the orthogonal complement of $\lieG^i$ in $\lieG^{i+1}$ with respect
to $\bilinear$, and let $\lieG^{i}_{\perp}$ be the orthogonal complement of 
$\lieG^i$ in $\lieG$ with respect to $\bilinear$. Then, we have 
\[
\lieG^{i+1}=\lieG^i\oplus\lieG^{i\perp},\qquad
\lieG=\lieG^i\oplus\lieG^{i}_{\perp},
\]
and
\[
\lieG=
\lieG^0\oplus\lieG^{0\,\perp}\oplus\cdots\oplus\lieG^{d-1\,\perp},
\qquad
\lieG^i_\perp=
\lieG^{i\,\perp}\oplus\cdots\oplus\lieG^{d-1\,\perp}.
\]
For $x\in\Bd(\bG^i,k)$ and $\rtm\in\tbR$, 
let $\lieG^{i\,\perp}_{x,\rtm}=\lieG^{i+1}_{x,\rtm}\cap\lieG^{i\perp}$.
Then, for $\xx\in\Bd(\bG^0,k)$, we can write $\vec\lieG_{\xx,\vec\utm}$
more explicitly:
\[
\vec\lieG_{\xx,\vec\utm}=\lieG_{\xx,f_{\vec\utm}}
=\lieG^0_{\xx,\utm_0}\oplus\lieG^{0\,\perp}_{\xx,\utm_1}\oplus\cdots
\oplus\lieG^{d-1\,\perp}_{\xx,\utm_d}.
\]
\end{numbering}

\begin{definition}\label{def: f^ast}\rm\ 

\begin{enumerate}
\item 
Define an involution $\ast$ on $\tbR$ as follows:
for $r\in\bbR$, $\rtm^\ast:=(-\rtm)^+$ and $(\rtm^+)^\ast:= -\rtm$.
\item
For a sequence $\vec\rtm=(\rtm_1,\cdots,\rtm_d)$ of real numbers,
define $\vec\rtm^\ast$ as $(\rtm^\ast_1,\cdots,\rtm^\ast_d)$.
\item
Assume $\HypB$ is valid. When we identify $\blieG$ and $\blieG^\ast$ 
via $\bilinear$, the dual $(\lieG_{x,\rtm})^\ast\subset\lieG$ 
of $\lieG_{x,\rtm}$ with respect to $\bilinear$ is $\lieG_{x,\rtm^\ast}$. 
That is,
$(\lieG_{x,\rtm})^{\ast}
=\{Y\in\lieG\mid\bilinear(Y,\lieG_{x,\rtm})\subset\pid_k\}
=\lieG_{x,\rtm^\ast}$. Note that $(\lieG_{x,\rtm})^\ast\subset\lieG$ 
while  $\lieG_{x,\rtm}^\ast\subset\lieG^\ast$.
Generalizing this, for any subset $\CaL\subset\lieG$, 
we define the dual $\CaL^{\ast}$ of $\CaL$ in $\lieG$ as
\[
\CaL^{\ast}:=\{Y\in\lieG\mid \bilinear(Y,\CaL)\subset\pid_k\}.
\]
\end{enumerate}
\end{definition}

\begin{remark}\label{rmk: dualL}\rm
Let $\vec\bG:=(\bG^0,\cdots,\bG^d)$ be a tamely ramified twisted Levi sequence.
Let $\xx\in\Bd(\bG^0,k)$, and
$\vec\utm:=(\utm_0,\cdots,\utm_d)$ be an admissible sequence.
Assume the hypothesis $\HypB$ is valid. Then, we have (\cite{Asym2})
\[
(\vec\lieG_{\xx,\vec\utm})^{\ast}=\vec\lieG_{\xx,\vec\utm^\ast}
=\lieG^0_{\xx,\utm_0^\ast}\oplus\lieG^{0\,\perp}_{\xx,\utm_1^\ast}\oplus\cdots
\oplus\lieG^{d-1\,\perp}_{\xx,\utm_d^\ast}.
\]
\end{remark}

\begin{definition}\rm
Suppose $\HypB$ is valid.
Let $L$ be an open compact subgroup of $\G$ with $L=\exp(\CaL)$
for some lattice $\CaL$ in $\lieG$.
Let $\chi$ be a character of $L$. 
If there is a $\gamma\in\lieG$ such that 
$\chi(g)=\addch(\bilinear(\gamma,\log(g)))$, we say that 
$\chi$ is \emph{represented} by $\gamma$, and we write 
$\chi_\gamma$ for $\chi$.
In this case, any element $\gamma'$ in the coset 
$\gsc=\gamma+\CaL^{\ast}\in \lieG\left/\CaL^{\ast}\right.$ represents $\chi$.
We call $\gsc$ the \emph{dual blob} of $(L,\chi)$. 
\end{definition}

\begin{remark}\rm
It is possible that $\gamma\in\lieG$ represents characters on
different groups, say $L_1$ and $L_2$. However, since those characters
coincide on $L_1\cap L_2$, we will use $\chi_\gamma$ to denote 
both characters when there is no confusion. 

For a sufficient condition for $\gamma\in\lieG$
to represent a character of $\vec\G_{\xx,\vec\utm}$,
we refer to \cite[(3.3.4)]{Asym2}.
\end{remark}

The proof of the following lemma is similar to that of Lemma 3.1 in
\cite[p17]{HM1}.

\begin{lemma}\label{lem: dual blob}
Suppose $\HypB$ and $\Hypk$ are valid.
For $i=1,2$, let $L_i$ be an open compact subgroup with
$L_i=\exp(\CaL_i)$ for some $\CaL_i\subset\lieG$. Let 
$\chi_i$ be a character of $L_i$ with dual blob $\gamma_i+\CaL_i^\ast$.
Suppose $\chi_1=\chi_2$ on $L_1\cap L_2$.
Then, $(\gamma_1+\CaL_1^\ast)\cap(\gamma_2+\CaL_2^\ast)\neq\emptyset$.
\end{lemma}

\proof
Since both $\gamma_1$ and $\gamma_2$ represent 
$\chi_1|(L_1\cap L_2)=\chi_2|(L_1\cap L_2)$, we have 
$\gamma_1+(\CaL_1\cap\CaL_2)^\ast=\gamma_2+(\CaL_1\cap\CaL_2)^\ast$.
Now, the lemma follows from $(\CaL_1\cap\CaL_2)^\ast=\CaL_1^\ast+\CaL_2^\ast$.
\qed

\section{\bf Semisimple elements}

\begin{numbering} {\bf Depth functions and good elements.}
\label{num: depthgss} \rm
Recall that the depth function 
$\dth:\Bd(\G)\times\lieG\longrightarrow\bbR$ is defined as follows:
for $X\in\lieG$ and $x\in\Bd(\G)$, let
$\dth(x,X)=\rtm$ be the depth of $X$ in the $x$-filtration, that is,
$r$ is the unique real number such that 
$X\in\lieG_{x,\rtm}\setminus\lieG_{x,\rtm^+}$. 
We also define 
\begin{equation*}
\dth(X)=\sup_{x\in\Bd(\G)}\dth(x,X).
\end{equation*}
Note that if $\dth(X)<\infty$, then
the depth $\dth(X)$ of $X$ is the unique $\rtm$ 
in $\bbR$ such that $X\in\lieG_{\rtm}\setminus\lieG_{\rtm^+}$.
Moreover, $\dth$ is locally constant on $\lieG\setminus\nil$,
and it is $\infty$ on $\nil$ (see \cite[(3.3.7)]{AD}).
If $E$ is a finite extension of $k$, we can 
also define a depth function $\dth^E$ on $\Bd(\bG,E)\times\lieG(E)$.
If $E$ is tamely ramified over $k$, thanks to our normalization of valuation, 
we observe that for any $x\in\Bd(\G)$ and $X\in\lieG$, $\dth(x,X)=\dth^E(x,X)$
and $\dth(X)=\dth^E(X)$ (see \cite[(2.2.5)]{AD2}) 
Hence we may omit the superscript $^E$ in that case.
We remark that if $X$ has Jordan decomposition $X_s+X_n$, 
with semisimple part $X_s$, then $\dth(X)=\dth(X_s)$ (\cite[(3.3.8)]{AD}).

Let $\bT$ be a maximal $k$-torus in $\bG$, 
and let $\blieT$ be its Lie algebra. 
Then $T$ and $\lieT$ have the following filtrations: 
for $r\in\bbR$ 
\[
\lieT_r
=\{\gss\in\lieT\mid\nu(d\chi(\gss))\ge r\textrm{ for all }
\chi\in\bX^\ast(\bT)\}
\]
and for $r>0$,
\[
\T_{r}
=\{t\in\T\mid\nu(\chi(t)-1)\ge r\textrm{ for all }\chi\in\bX^\ast(\bT)\}.
\]
Note that 
if $\bT$ is $k$-split, the lattice $\lieT_r$ coincides with the lattice
$\blieT(k)_\rtm$ of (\ref{num: root dec}). 
The following definition is from \cite[\S5]{AR}.

\begin{definition}\label{def: gss}
\rm 
Let $\bT$ be a maximal $k$-torus in $\bG$ which splits over a tamely ramified
Galois extension of $k$, and let $\blieT$ be its Lie algebra.
\begin{enumerate}
\item
If $\gss\in\lieT_r\setminus\lieT_{r^+}$, we say that $\gss$ is 
{\it of depth $r$ with respect to $\bT$}, and we write $\dth_{\bT}(\gss)=r$. 
\item 
Let $\gss\in\lieT$ be of depth $r$.
Then $\gss$ is called {\it good with respect to $\bT$} 
if for every root $\alpha$ of $\bG$ with respect to $\bT$, 
$d\alpha(\gss)$ is either zero or has valuation $r$.
\end{enumerate}
\end{definition}

Note that $0\in\lieG$ is a good element of depth $\infty$. 
We remark that the depth and the goodness of a semisimple element
do not depend on the choice of $\bT$ (\cite[(5.1)]{AR}).

\end{numbering}

\begin{numbering}\label{num: recall facts}\rm
We recall some useful facts about degenerate cosets 
(see \cite[(3.2.6)]{AD} and \cite[(6.3)]{MP}).
\begin{enumerate}
\item 
Assume $X\in\lieG_{x,\rtm}\cap\lieG_{\rtm^+}$. Then,

(i) $X+\lieG_{x,\rtm^+}$ contains a nilpotent element, and

(ii) there is a $y\in\Bd(\bG,k)$ such that 
$X+\lieG_{x,\rtm^+}\subset\lieG_{y,\rtm^+}$.
\item 
Let $\stm\in\bbR$. Then,
$\lieG_s=\cap_{x\in \Bd(\bG,k)}\, (\nil+\lieG_{x,s})$.
\item  
If $x\in \Bd(\bG,k)$ and $\stm\in\bbR$, then
$\lieG_{x,s}\cap \lieG_{s^+}=(\nil\cap \lieG_{x,s}) + \lieG_{x,s^+}$.
\end{enumerate}
\end{numbering}

\begin{numbering}\label{num: Cent gss}\rm 
Let $\bT$ be a maximal $k$-torus in $\bG$
which splits over a tamely ramified Galois extension $E$.
Then we observe the following:
\begin{enumerate}
\item
Let $\Phi:=\Phi(\bG,\bT,E)$. 
Then we have 
\[
\blieG(E)
=\blieT(E) \ 
\oplus\sum_{\alpha\in\Phi}\blieG_\alpha(E),
\]
Let $\gss\in\lieT$ be a semisimple element. 
Let $\bH:=\bC_\bG(\gss)$ be the centralizer of $\gss$ in $\bG$.
The Lie algebra $\blieH(E)$ can be expressed as follows \cite[(7.1)]{yu}:
\[
\blieH(E)
=\blieT(E) \ 
\oplus\sum_{\alpha\in\Phi,\ d\alpha(\gss)=0}\blieG_\alpha(E).
\]
\item
If $\gss\notin\lieZ_\lieG$ is a good element of depth $\rtm$, 
then for any $\gamma\in\lieZ_{\lieG}$ with $\dth(\gamma)\ge\rtm$,
$\gss+\gamma$ is also a good element of depth $\rtm$.
\end{enumerate}
\end{numbering}

\begin{lemma}\label{lem: appr}
Let $\bT$ be a maximal $k$-torus in $\bG$
which splits over a tamely ramified Galois extension $E$.
Let $\gamma=\gamma_1,\cdots,\gamma_n \in\lieT$ be good elements of
depth $\btm=\btm_1,\cdots,\btm_n$ respectively. Let $\bH^0:=\bG$,
and $\bH^i:=\bC_{\bH^{i-1}}(\gamma_i)$.
\begin{enumerate}
\item
Let $X_1$, $X_2\in\gamma+\lieH^1_{\btm^+}$. 
If $g\in\G$ is such that $^gX_1=X_2$, then $g\in H^1$.
\item
Let $X\in\gamma+\lieH^1_{\btm^+}$. 
Then $C_\G(X)\subset H^1$.
\item
Suppose each $\gamma_i$ is a good element in $\bH^{i-1}$,
and $\btm_1<\btm_2<\cdots<\btm_n$.
Fix $i\,\in\,\{0,1,\cdots,d\}$ and
let $\gamma^i=\gamma_1+\gamma_2+\cdots+\gamma_i$.
Let $Y_1, Y_2\in\gamma^i+\lieH^i_{\btm_i^+}$.
If $^gY_1=Y_2$ for some $g\in\G$, then $g\in H^i$.
Moreover, $\bH^i=\bC_\bG(\gamma^i)$
\end{enumerate}
\end{lemma}

\proof 
(1) is \cite[(2.3.6)]{Asym}. 
(2) follows from (1).
For the first statement of (3), if $i=1$, it is (1). 
Assume the statement is true for $i$. 
Let $Y_1, Y_2\in\gamma^{i+1}+\lieH^{i+1}_{\btm_{i+1}^+}$.
Suppose $^gY_1=Y_2$ for some $g\in\G$. Since 
$Y_1,Y_2\in\gamma^i+\lieH^i_{\btm_i^+}$, we have $g\in H^i$
by induction hypothesis and $^g(Y_1-\gamma^i)=Y_2-\gamma^i$.
Since 
$Y_1-\gamma^i,\ Y_2-\gamma^i\in \gamma_{i+1}+\lieH^{i+1}_{\btm_{i+1}^+}$
and $\gamma_{i+1}$ is $\bH^i$-good, we have $g\in H^{i+1}$.
For the second statement, $\bH^i\subset\bC_\bG(\gamma^i)$ is obvious, and
$\bC_\bG(\gamma^i)\subset\bH^i$ follows from the first with $Y_1=Y_2=\gamma$.
\qed

\begin{lemma}\label{lem: sum gss}
Let $\bT$ be a maximal $k$-torus in $\bG$
which splits over a tamely ramified Galois extension $E$.
Let $\gamma_1,\gamma_2 \in\lieT$ be good elements of
depth $\btm_1,\btm_2$ respectively. Let $\bH:=\bC_\bG(\gamma_1)$. 
\begin{enumerate}
\item
Suppose $\btm=\btm_1=\btm_2$ and 
$\gamma_1\equiv\gamma_2\pmod{\lieT_{\btm^+}}$.
Then $\C_\G(\gamma_1)=\C_\G(\gamma_2)$.
\item
Suppose $\btm_1<\btm_2$ and $\gamma_1, \gamma_2\in\lieZ_{\lieH}$.
Then $\gamma_1+\gamma_2$ is also a $\G$-good element of depth $\btm_1$.
\end{enumerate}
\end{lemma}

\proof
(1) Note that $\gamma_1,\ \gamma_2\in\gamma_1+\lieH_{\btm^+}$.
Applying Lemma \ref{lem: appr}-(2), 
we have $\C_\G(\gamma_1)\subset\C_\G(\gamma_2)$.
Similarly, $\C_\G(\gamma_2)\subset\C_\G(\gamma_1)$. 
Hence $\C_\G(\gamma_1)=\C_\G(\gamma_2)$.

(2) Write $\gss=\gamma_1+\gamma_2$.
Let $\Phi=\Phi(\bG,\bT,E)$ be the set of $E$-rational $\bT$-roots in $\bG$.
Let $\alpha\in\Phi$. Since $\bH\subset\bC_{\bG}(\gamma_2)$,
by (\ref{num: Cent gss}), we see that
if $d\alpha(\gamma_1)=0$, then $d\alpha(\gamma_2)=0$.
Combining this with $\val(d\alpha(\gamma_2))\ge\btm_2>\btm_1$,
we see that $d\alpha(\gamma_1+\gamma_2)=0$ or 
$\val(d\alpha(\gamma_1+\gamma_2))
=\min(\val(d\alpha(\gamma_1)),\ \val(d\alpha(\gamma_2)))
=\val(d\alpha(\gamma_1))=\btm_1$.
Hence $\gamma_1+\gamma_2$ is a good element of depth $\btm_1$.
\qed

\begin{proposition}\label{prop: goodexp}
Suppose $\HypGT$ is valid. 
Let $\gamma\in\lieG$ be a semisimple element 
which splits over a tamely ramified Galois extension $E$. 
Then, $\gamma$ can be written as 
\[
\gamma=\gamma_{\btm_1}+\gamma_{\btm_2}+\cdots+\gamma_{\btm_n}+\gamma_\circ
\]
such that 
\begin{enumerate}
\item
each $\gamma_{\btm_i}$, $i=1,\cdots,n$ ,
is a $\G$-good element of depth $\btm_i$, and
$\gamma_\circ$ is a semisimple element with $\dth(\gamma_\circ)\ge0$, 

\item
$\btm_1<\btm_2<\cdots<\btm_n<0$, and

\item 
$\bH^{n}\subsetneq\bH^{n-1}\subsetneq\cdots\subsetneq\bH^1\subseteq\bG$
where $\bH^1=\bC_\bG(\gamma_{\btm_1})$, and
$\bH^i
=\bC_{\bH^{i-1}}(\gamma_{\btm_i})$.
\end{enumerate}

\noindent
Moreover, if 
$\gamma
=\gamma'_{\btm'_1}+\gamma'_{\btm'_2}+\cdots
+\gamma'_{\btm'_{n'}}+\gamma^{\prime}_\circ$
is another expression satisfying (1)-(3) with
$\bH^{\prime 1}=\bC_\bG(\gamma'_{\btm'_1})$ and
$\bH^{\prime i}=\bC_{\bH^{\prime {i-1}}}(\gamma'_{b'_i})$,
then we have
$n=n'$, $\btm_i=\btm_i'$ and $\bH^i=\bH^{\prime i}$.
\end{proposition}

\proof
Let $\lieT$ be a maximal $k$-torus with $\gamma\in\lieT$ which splits 
over $E$. If $\dth(\gamma)\ge0$, $\gamma=\gamma_\circ$ already satisfies
(1)-(3). Suppose $\atm_1=\dth(\gamma)<0$.
By $\HypGT$, $\gamma+\lieT_{\atm_1^+}$ 
contains a good element, say, $\tilde\gamma_{\atm_1}$ of 
depth $\atm_1$. Then 
$\gamma=\tilde\gamma_{\atm_1}+(\gamma-\tilde\gamma_{\atm_1})$ with
$\atm_2=\dth(\gamma-\tilde\gamma_{\atm_1})>\dth(\gamma)$. 
Applying the above process
for $\gamma-\tilde\gamma_{\atm_1}$, we find a $\G$-good element 
$\tilde\gamma_{\atm_2}\in \gamma-\tilde\gamma_{\atm_1}\,+\,\lieT_{\atm_2^+}$
such that 
$\gamma
=\tilde\gamma_{\atm_1}+\tilde\gamma_{\atm_2}
+(\gamma-\tilde\gamma_{\atm_1}-\tilde\gamma_{\atm_2})$
and $\atm_3=\dth(\gamma-\tilde\gamma_{\atm_1}-\tilde\gamma_{\atm_2})
<\dth(\gamma-\tilde\gamma_{\atm_1})$.
Repeatedly, we have
\[
\gamma=\tilde\gamma_{\atm_1}+\tilde\gamma_{\atm_2}
+\cdots+\tilde\gamma_{\atm_m}+\gamma_\circ
\]
where $\tilde\gamma_{\atm_i}$ is a $\G$-good element 
of depth $\atm_i$ with
$\atm_1<\atm_2<\cdots<\atm_m<0$ and $\dth(\gamma_\circ)\ge0$.
This procedure is finite because $\dth(\lieT)\subset\frac1{e(E/k)}\bbZ$.
Put $\atm_{m+1}=0$ and $\tilde\gamma_{\atm_{m+1}}=\gamma_\circ$.

Set $\Set:=\{\atm_1,\atm_2,\cdots,\atm_{m+1}\}$, and
for $\atm,\btm\in\tbR$, set 
$\tilde\gamma_{\atm,\btm}:=\sum_{\atm\le\atm_j<\btm}\tilde\gamma_{\atm_j}$.
We find a subsequence 
$\btm_1=\atm_1<\btm_2<\cdots<\btm_n<\btm_{n+1}=\atm_{m+1}=0$ 
of $\Set$ as follows:
let $b_1:=a_1$ and $\bH^1:=\bC_\bG(\tilde\gamma_{b_1})$.
Let $\btm_2$ be the maximal element in $\{\atm_2,\cdots,\atm_{m+1}\}$ 
with the property that 
if $a_j<\btm_2$, $\tilde\gamma_{a_j}\in\blieZ_{\blieH_{1}}$.
Note that $\bH^{1}=\bC_\bG(\tilde\gamma_{\btm_1,\btm_2})$.
Let $\bH^{2}:=\bC_{\bH^1}(\tilde\gamma_{\btm_2})$.
Then $\bH^{1}\supsetneq\bH^{2}$. 
Let $\gamma_{b_1}:=\tilde\gamma_{\btm_1,\btm_2}$.
Inductively, suppose $b_i$, $\bH^i$ and $\gamma_{b_{i-1}}$
are defined for $i\ge2$.
Let $\btm_{i+1}$ be the maximal element in 
$\{\atm_j\in\Set\mid\atm_j>\btm_i\}$ 
with the property that for any $\atm_j<\btm_{i+1}$, 
$\tilde\gamma_{a_j}\in\lieZ_{\lieH^{i}}$. 
Let $\bH^{i+1}:=\bC_{\bH^{i}}(\tilde\gamma_{\btm_{i+1}})$
and $\gamma_{b_{i}}=\tilde\gamma_{b_{i},b_{i+1}}$.
We repeat the process until $\btm_{n+1}=0$.
Then each $\gamma_{\btm_i}$ is also a $\G$-good element of depth $\btm_i$
by Lemma \ref{lem: sum gss}-(2), and 
$\bH^{i}=\bC_{\bH^{i-1}}(\tilde\gamma_{\btm_i})
=\bC_{\bH^{i-1}}(\gamma_{\btm_i})$.
Now, one can easily check
\[
\gamma =\gamma_{\btm_1}+\gamma_{\btm_2}+\cdots+\gamma_{\btm_n}+\gamma_\circ
\]
satisfies the required properties. 

\smallskip

For the second statement, 
let $\gamma
=\gamma'_{\btm'_1}+\gamma'_{\btm'_2}+\cdots
+\gamma'_{\btm'_{n'}}+\gamma^{\prime}_\circ$
be another expression satisfying (1)-(3).
Then $\btm_1=\dth(\gamma)=\btm'_1$.
Since $\gamma_{\btm_1}\equiv\gamma'_{\btm'_1}\pmod{\lieT_{\btm_1^+}}$
and $\gamma_{\btm_1}$, $\gamma'_{\btm'_1}$ are good, 
$\bH^{1}=\bH^{\prime 1}$ by  Lemma \ref{lem: sum gss}.
By induction, we assume that $\btm_j=\btm'_j$ and 
$\bH^{j}=\bH^{\prime j}$ for $1\le j\le i-1$.
Write $\gamma_{\btm_1,\btm_i}:=\gamma_{\btm_1}+\cdots+\gamma_{\btm_{i-1}}$,
and $\gamma'_{\btm_1,\btm_i'}:=\gamma'_{\btm_1}+\cdots+\gamma'_{\btm_{i-1}}$.
Suppose $\btm_i<\btm'_i$. Then, 
$\gamma'_{\btm_1,\btm_i'}-\gamma_{\btm_1,\btm_i}
\equiv\gamma_{\btm_i}\pmod{\lieT_{\btm_i^+}}$.
Since $\gamma_{\btm_i}$ is also $\bH^{i-1}$-good, 
by Lemma \ref{lem: appr}-(2), we have 
$\bH^{i-1}\subseteq\bH^{i}$, which is a contradiction. 
Hence $\btm_i=\btm'_i$.
Now we have
$(i)$ $\gamma_{\btm_i}\equiv\gamma'_{\btm_1,\btm_i}-\gamma_{\btm_1,\btm_i}
+\gamma'_{\btm_i} \pmod{\lieT_{\btm_i^+}}$, 
$(ii)$ $\gamma'_{\btm_1,\btm_i}-\gamma_{\btm_1,\btm_i}
\in\lieZ_{\lieH^{i-1}}$, and 
$(iii)$ $\gamma_{\btm_i}$, 
$\gamma'_{\btm_1,\btm_i}-\gamma_{\btm_1,\btm_i}+\gamma'_{\btm_i}$
are good in $\bH^{i-1}$ by (\ref{num: Cent gss}).
From these, it follows that 
$\bC_{\bH^{i-1}}(\gamma_{\btm_i})
=\bC_{\bH^{i-1}}(\gamma'_{\btm_i})$. 
Hence, $\bH^{i}=\bH^{\prime i}$.
\qed

\begin{remark}\rm 
By Lemma \ref{lem: appr}, we have
$\bH^{j}=\bC_\bG(\gamma_{\btm_1}+\cdots+\gamma_{\btm_j})$ and
$\bC_\bG(\gamma)\subset\bH^{n}$.
\end{remark}

\medskip

\section{\bf $\bK$-types: basic data and construction}
\label{sec: Ktypes}

In this section, we review the $\bK$-types constructed in \cite{Asym2}.
We will prove later that under some hypotheses, almost all irreducible
tempered representation of positive depth 
contains one of these types (see \S\ref{sec: tem}).
The construction is based on the following data:

\begin{definition} \label{def: posGdatum} \rm
A \emph{$\G$-datum of positive depth} is a quadruple 
$\datum=(\vec{\bG},\xo,\vec r,\vec\phi)$  satisfying the following 
conditions (D1)--(D4):

\midvsp

\item{(D1)}
$\vec{\bG}=(\bG^0\subsetneq\bG^1\subsetneq\cdots\subsetneq\bG^d=\bG)$ 
is a tamely ramified twisted Levi sequence.

\midvsp

\item{(D2)}
$\xo\in\Bd(\bG^0,k)$.

\midvsp

\item{(D3)}
$\vec\rtm=(\rtm_0,\rtm_1,\cdots,\rtm_{d-1},\rtm_d)$ is
a sequence of positive real numbers 
with $0<\rtm_0<\cdots<\rtm_{d-2}< \rtm_{d-1}\le\rtm_d$.

\midvsp

\item{(D4)}
$\vec\phi=(\phi_0,\cdots,\phi_d)$ is a sequence of quasi-characters,
where $\phi_i$ is a quasi-character of $\G^i$.
We assume that $\phi_i$ is trivial on $\G^i_{\xo,\rtm_i^+}$, but
non-trivial on $\G^i_{\xo,\rtm_i}$ for $0\le i\le d-1$.
If $\rtm_{d-1}<\rtm_d$, we assume 
$\phi_d$ is nontrivial on $\G^d_{\xo,\rtm_d}$ and trivial on 
$\G^d_{\xo,\rtm{}^+_d}$. Otherwise, we assume that $\phi_d=1$.

\medskip

We define the length $\ell(\datum)$ of the above $\G$-datum of positive depth
$\datum$ to be $d$.

\end{definition}

\begin{remark}\rm 
As mentioned in the introduction, Yu defined a notion of $\G$-datum for 
supercuspidal types (on which the above definition of  
$G$-datum of positive depth is based). However, the representations
containing $\bK$-types constructed out of the above data
are not necessarily supercuspidal.
\end{remark}

\medskip

\noindent
{\sc Notation and Conventions.}
\begin{enumerate}
\item
Let $\bG^{-1}:=\bG^0$,\ \  $\blieG^{-1}:=\blieG^{0}$,
and $\bG^{d+1}:=\bG^{d}$,\ \  $\blieG^{d+1}:=\blieG^{d}$.
\item
Let $\bZ_{\bG^i}=\bZ^i$ and $\blieZ_{\blieG^i}=\blieZ^i$ denote 
the centers of $\bG^i$ and $\blieG^i$ respectively. 
\end{enumerate}

\medskip

\begin{definition}\label{def: goodatum}\rm
Let $\datum=(\vec\bG,\xo,\vec\rtm,\vec\phi)$ be a $\G$-datum of positive
depth.

\begin{enumerate}
\item 
We say $\datum$ is \emph{\good} if each $\phi_i$, $0\le i\le d-1$, is good
and $\phi_d$ is either good or trivial. That is,
for $0\le i\le d-1$, \ $\phi_i|\G^i_{\xo,\rtm_i}$ is represented by 
a $\G^{i+1}$-good element $\gss_i\in\lieZ^i$ of depth 
$-\rtm_i$ such that $\bC_{\bG^{i+1}}(\gss_i)=\bG^i$, 
and for $i=d$, either $\phi_d$ is trivial or 
$\phi_d|\G^d_{\xo,\rtm_d}$ is represented 
by $\gss_d\in\lieZ^d$ of depth $-\rtm_d$.

\item
We say $\phi_i$ is \emph{strongly good} if 
$\phi_i|\G^i_{\xo,0^+}$ is represented by 
a $\G$-good element $\gss_i\in\lieZ^i$ of depth $-\rtm_i$ 
such that $\bC_{\bG^{i+1}}(\gss_i)=\bG^i$.

\item
We say $\datum$ is {\it strongly good} if each $\phi_i$, $0\le i\le d-1$, 
is strongly good and $\phi_d$ is either strongly good or trivial. 
Put $\gss_d=0$ if $\phi_d$ is trivial, and
define $\gss_\datum$ (or simply $\gss$)
and $\gss^i$ as follows:
\begin{align*}
&\gss^i:=\gss_d+\gss_{d-1}+\cdots+\gss_i,\\
&\gss_\datum=\gss:=\gss^0=\gss_d+\gss_{d-1}+\cdots+\gss_0.
\end{align*}
\end{enumerate}
\end{definition}

\begin{remarks}\rm 
Let $\datum$ be a $\G$-datum of positive depth.
\begin{enumerate}
\item
The expression of $\gss_\datum$ as $\gss_d+\gss_{d-1}+\gss_0$ satisfies 
the condition in Proposition \ref{prop: goodexp}. 
\item
In the above definition, $\gss_i$ being $\G$-good implies that
$\gss_i$ is $\G^{i+1}$-good. 
Hence, if $\datum$ is strongly good, it is also good.
Note also that our definition of strongly goodness
is stronger than the one in \cite{Asym2}.
\item
If  $\phi_i$ is good, it is also {\it generic} in the sense of \cite{yu}.
Hence we can apply most results in \cite{yu} to 
a (strongly) good $\G$-datum of positive depth.
\item
Observe that each $\gss_i\in\lieZ_{\lieG^0}$. Hence the $\gss_i$'s
commute with each other.
\end{enumerate}
\end{remarks}

\begin{lemma}\label{lem: sstochar}
Suppose $\HypB$ and $\Hypk$ are valid. 
Suppose $\gamma\in\lieZ_\lieG$ and $\dth(\gamma)<0$.
Then there is a character $\phi$ of $\G$ such that for any $x\in\Bd(\bG,k)$,
$\phi|\G_{x,0^+}$ is represented by $\gamma$.
\end{lemma}

\proof
Write $Z=Z_\G$ for simplicity. Set $Z_{0^+}=\exp(\lieZ_{0^+})$
where $\lieZ_{0^+}=\lieZ\cap\lieT_{0^+}$ with $\bT$
a maximal $k$-torus in $\bG$.
Since $Z$ is commutative, by $\Hypk$, $\chi_\gamma$ defines
a character of $Z_{0^+}$. 
Since the commutator $(\G,\G)$ is a subgroup of $\G^\der$, 
$\G/\G^\der$ is abelian.
Moreover, since we have an orthogonal decomposition 
$\lieG_{x,0^+}=\lieZ_{0^+}\oplus\lieG^\der_{x,0^+}$
with respect to the bilinear form $\bilinear$ by \cite[(3.2)]{AR},
we have $\G_{x,0^+}=\exp(\lieZ_{0^+}\oplus\lieG^\der_{x,0^+})
=Z_{0^+}\G^{\der}_{x,0^+}$ by $\Hypk$.
Hence $Z_{0^+}$ is embedded in $\G/G^\der$, and 
$\chi_\gamma$ defines a character $\overline\chi_\gamma$
of $Z_{0^+}\G^\der/\G^\der$. This character easily extends to
a character of $Z\G^\der/\G^\der$, which we again denote by 
$\overline\chi_\gamma$. 
Now, since $\G/(Z\G^\der)$ is finite, 
$\overline\chi_\gamma$ again extends to a character of $\G/G^\der$,
which induces a character $\phi$ of $\G$.

It remains to show that $\phi|\G_{x,0^+}$ is represented by $\gamma$.
For $g\in\G_{x,0^+}$, we can write $g=\exp(z)\exp(z')=\exp(z+z')$
for some $z\in\lieZ_{0^+}$ and $z'\in\lieG^\der_{x,0^+}$.
Then, since $\phi$ is trivial on $\G^\der$ and 
$\lieG=\lieZ\oplus\lieG^\der$ is an 
orthogonal decomposition with respect to $\bilinear$,
we have 
$\phi(g)=\phi(\exp(z))=\chi_\gamma(\exp(z))=\addch(\bilinear(\gamma,z))
=\addch(\bilinear(\gamma,z+z'))=\addch(\bilinear(\gamma,\log(g)))$.
Hence $\phi|\G_{x,0^+}$ is represented by $\gamma$.
\qed

\begin{remarks}\rm \ 
\begin{enumerate}
\item
Let $\phi$ be constructed as in the proof of Lemma \ref{lem: sstochar}.
Then $\phi$ is trivial on $\G^\der$ and the depth of $\phi$ is 
$\rtm=-\dth(\gamma)$ independent of $x\in\Bd(\bG,k)$. 
That is, for any $x\in\Bd(\bG,k)$, $\phi$ is trivial on $\G_{x,\rtm^+}$,
but, nontrivial on $\G_{x,\rtm}$.
\item
Combining Lemma \ref{lem: sstochar} and Proposition \ref{prop: goodexp}, 
for any semisimple element in $\lieG$ of negative depth, one can
associate a strongly good $\G$-datum of positive depth.
\end{enumerate}
\end{remarks}

\bigskip

Fix a good $\G$-datum of positive depth
$\datum=(\vec{\bG},\xo,\vec r,\vec\phi)$.
We review the construction of the $\bK$-type $(\rK^+_\datum,\phi_\datum)$ 
and some notation that we need in this paper. 

\begin{numbering}\rm
Let the embeddings
\[
\Bd(\bG^0,k)\hookrightarrow\Bd(\bG^1,k)\hookrightarrow\Bd(\bG^2,k)
\hookrightarrow\cdots\hookrightarrow\Bd(\bG^d,k)
\]
be fixed as in (\ref{num: embeddings}).
Let $\stm_i:=\frac{\rtm_i}2$ for $i=0,\cdots,d-1$.
For any $\epsilon\in\tbR$ with $0\le\epsilon<\stm_0$,
define the following sequences $\vec\stm$ and $\vec\stm^+$ of length $d$:
\[
\vec\stm(\epsilon)\ :=(\epsilon,\stm_0,\cdots,\stm_{d-1}), \qquad
\vec{\stm}^+(\epsilon)\ :=(\epsilon,\stm_0^+,\cdots,\stm_{d-1}^+).
\]
For simplicity, we write
\[
\vec\stm:=\vec\stm(0),\qquad \vec\stm^+:=\vec\stm(0^+).
\]
\end{numbering}

\begin{numbering}\label{num: construction Ktypes}\rm
Define some open compact subgroups associated to $\datum$ as follows:
\[
\begin{array}{lll}
\rK^{i+}_\datum & :=\ 
\G^0_{\xo,0^+}\G^1_{\xo,\stm_0^+}\cdots\G^i_{\xo,\stm_{i-1}^+}
&\subset\,\G^i_{\xo,0^+},\\
\rK^{i+1}_{\datum+} & :
\ =\ \G^{i+1}_{\xo,\stm_{i}^+}\cdots\G^d_{\xo,\stm_{d-1}^+}
&\subset\,\G^d_{\xo,\stm_{i}^+}, \\
\rK_\datum^+ & :=\ \rK^{i+}_\datum\rK^{i+1}_{\datum+}\ =\ 
\G^0_{\xo,0^+}\G^1_{\xo,\stm_0^+}\cdots\G^d_{\xo,\stm_{d-1}^+}
&=\ \vec\G_{\xo,\vec{\stm}^+}\ .
\end{array}
\]
Via the isomorphism 
$\G_{\xo,\stm_i^+}/\G_{\xo,\rtm_i^+}
\simeq\lieG_{\xo,\stm_i^+}/\lieG_{\xo,\rtm_i^+}$,
$\phi_i$ defines a character $\hat\phi_i$ of $\G_{\xo,\stm_i^+}$ 
such that $\hat\phi_i$ is trivial on 
$(\lieG_{\xo,\stm_i^+}\cap\lieG^i_\perp)
\left/(\lieG_{\xo,\rtm_i^+}\cap\lieG^i_\perp)\right.$, and 
$\hat\phi_i$ and $\phi_i$ coincide
on $\G^i_{\xo,\stm_i^+}$ (see \cite[\S4]{yu}). 
Since $\phi_i$ is already defined on $\rK^{i+}_\datum\subset\G^i_{\xo,0^+}$,
there is a unique character of $\rK^+_{\datum}\G_{\xo,\stm^+_i}$
extending $\phi_i$ and $\hat\phi_i$. We use the same notation $\hat\phi_i$ 
for this character of $\rK^+_{\datum}\G_{\xo,\stm^+_i}$
and its restriction to $\rK^+_\datum$.
Now, define the character $\phi_\datum$ of $\rK_\datum^+$ as 
$\prod_{i=0}^d\hat\phi_i$:
\[
\phi_\datum:=\prod_{i=0}^d\hat\phi_i\ .
\]
\end{numbering}

\begin{numbering}\rm 
Suppose $\Hypk$ is valid and $\datum$ is 
a strongly good $\G$-datum of positive depth. 
Let $\gss_i$ and $\gss$ be as in Definition \ref{def: goodatum}. 
We observe that $\hat\phi_i$ is represented by $\gss_i$
on $\G^i_{\xo,0^+}\G_{\xo,\stm_i^+}$ and hence on $\rK^+_\datum$, 
that is, $\hat\phi_i=\chi_{\gss_i}$ on $\G^i_{\xo,0^+}\G_{\xo,\stm_i^+}$
and on $\rK^+_\datum$. Moreover, we have
\[
\phi_\datum=\chi_{\gss}
\]
on $\rK^+_\datum$, and the dual blob of $\phi_\datum=\chi_\gss$ is
$\gss+(\vec\lieG_{\xo,\vec{\stm}^+})^{\ast}=\gss+\vec\lieG_{\xo,-\vec\stm}$,
where $-\vec\stm=(-\stm_{-1},-\stm_0,\cdots,-\stm_{d-1})$. 
We also observe 
\[
\gss^i:=\gss_d+\gss_{d-1}+\cdots+\gss_i
\]
defines a character $\chi_{\gss^i}$ of $\rK^{i+1}_{\datum+}$.
We also put $\gss^{-1}:=\gss^0=\gss$.
\end{numbering}

\begin{remarks}\label{rmk: xtoy} \rm 
Let $\datum=(\vec\bG,\xo,\vec\rtm,\vec\phi)$ be a strongly good
$\G$-datum of positive depth. 
\begin{enumerate}
\item
Since the construction of $(\rK_\datum^+,\phi_\datum)$ depends only 
on $\gss_i$ representing $\phi_i$ on $\G^i_{\xo,0^+}$, 
replacing $\phi_i$ with the characters constructed in Lemma \ref{lem: sstochar}
with $\gamma=\gss_i$ produces the same open compact subgroup 
and its representation.
Hence, without loss of generality, we may and {\bf will} assume that 
for a strongly good $\G$-datum of positive depth $\datum$,
each $\phi_i$ is represented by $\gss_i$ on $\G^i_{x,0^+}$
for any $x\in\Bd(\bG^i,k)$.
In this case, for any $x\in\Bd(\bG^0,k)$, 
$(\vec\bG,x,\vec\rtm,\vec\phi)$ is also a strongly good $\G$-datum 
of positive depth. 
We will often denote $(\vec\bG,x,\vec\rtm,\vec\phi)$ by $\datum_x$.
\item
Since $\gss$ determines $\vec\bG$, $\vec\rtm$ and $\phi_\datum$,
$\datum$ can be replaced by $(\gss,\xo)$ and they yield the same
$\bK$-type $(\rK^+_\datum=\vec\G_{\xo,\vec\stm^+},\chi_\gss)$.
\end{enumerate}
\end{remarks}

\begin{remarks} {\sc Some Properties of $(\rK_\datum^+,\phi_\datum)$.}\rm \ 
\begin{enumerate}
\item
If an irreducible admissible representation $(\pi,V_\pi)$ contains
$(\rK_\datum^+,\chi_\gss)$, the depth $\dpi(\pi)$ of $\pi$ is $\rtm_d$.
\item
Let $\datum=(\vec\bG,\xo,\vec\rtm,\vec\phi)$ be a good $\G$-datum 
of positive depth.
Suppose $\bZ_{\bG^0}/\bZ_{\bG}$ is anisotropic. Let
$\rho$ be an irreducible representation of $\G^0_{[\xo]}$, the stabilizer
in $\G^0$ of the image $[\xo]$ of $\xo$ in the reduced building of $\bG^0$,
such that $c\textrm{-Ind}_{\G^0_{[\xo]}}^{\G^0}\rho$ is
irreducible supercuspidal.
Then $(\vec\bG,\xo,\vec\rtm,\vec\phi,\rho)$ is a generic $\G$-datum 
(\cite{yu}, see also \S\ref{sec: sc}). 
Moreover, if $\pi$ is a supercuspidal representation constructed from 
the above generic datum, $\phi_\datum$ occurs in 
the restriction $\pi|\rK^+_\datum$ of $\pi$ to $\rK^+_\datum$.
\end{enumerate}
\end{remarks}

For more details and properties of the above $\bK$-types, 
we refer to \cite{Asym2}.

\begin{numbering}\label{num: aux lattices} \rm\ 
Let $\datum=(\vec\bG,\xo,\vec\rtm,\vec\phi)$ be as before.
We recall some lattices and open compact subgroups associated to $\datum$
from \cite{Asym2}.

\begin{enumerate}
\item
Let $\epsilon\in\bbR$ be such that $\lieG_{\epsilon}=\lieG_{0^+}$.
For $x\in\Bd(\bG^0,k)$, let
\[
\CaL_{x,\epsilon}\ =\CaL_x\ :=\vec\lieG_{x,\vec\stm^+(\epsilon)},\qquad
\CaL_{x,\epsilon}^\sharp\ =\CaL_x^\sharp\ 
:=\vec\lieG_{x,(\vec\stm^+(\epsilon))^\ast},
\]
and
\[
\rL_{x,\epsilon}\ =\rL_x\ :=\vec\G_{x,\vec\stm^+(\epsilon)}.
\]
Although these definitions $\CaL_{x,\epsilon},\ \rL_{x,\epsilon}$, etc, 
depend on $\vec\bG$, $\vec\rtm$ and the embedding of buildings, 
since we would not need their roles explicitly, 
we omit them from the subscripts  for simplicity.
If $\HypB$ is valid, $\CaL_{x,\epsilon}^\sharp$ is
just $\CaL_{x,\epsilon}^\ast$.
If there is no confusion, we will drop $\epsilon$ from the notation,
that is, we will write $\CaL_x$, $\CaL_x^\sharp$, 
$\rL_x$ for $\CaL_{x,\eps}$, $\CaL_{x,\eps}^\sharp$, $\rL_{x,\eps}$.
\item
Assume $\HypB$ is valid. 
For $i=1,\cdots,d+1$ and $x\in\Bd(\bG^i,k)$, define $\CaL^i_x$
and $\CaL^{i\,\sharp}_x$ as follows:
\begin{align*}
&\CaL^{d+1}_x\ 
=\ \CaL^d_x\ :=\ \lieG^d_{x,\stm_{d-1}^+},\\
&\CaL^i_x\ \quad
:=\ \lieG^i_{x,\stm_{i-1}^+}+\CaL^{i+1}_x 
\qquad\mathrm{for }\quad i=1,\cdots,d-1, \\
&\CaL^{i\,\sharp}_x\ \quad
:=\CaL^{i\,\ast}_x,
\qquad\qquad\qquad\mathrm{for }\quad i=1,\cdots,d+1.
\end{align*}
Set $\CaL_x^0:=\CaL_{x,\eps}$, $\CaL_x^{0\sharp}:=\CaL_{x,\eps}^\sharp$ 
and $\rL^0_x:=\rL_{x,\eps}$.
\end{enumerate}
\end{numbering}

For later use, we record the following which is 
a corollary of \cite[(5.3.2)]{Asym2}:

\begin{lemma} \label{lem: asym2appr} 
Assume the hypothesis $\Hypk$ is valid.
Let $x\in\Bd(\bG^0,k)$, and $X\in\lieG^0_{x,(-\rtm_0)^+}$. Then 
\[
\gss+X+\CaL_x^\sharp\subset\,^\G(\gss+X+\lieG^0_0).
\]
\end{lemma}

For our purpose, we also define a \emph{zero good $\G$-datum}.

\begin{definition}\rm
A \emph{zero $G$-datum} (or \emph{$0$-datum}) 
$\datum$ is a $\G$-datum of the form
$\datum=(\vec\bG,\xo,\vec\rtm,\vec\phi)$ 
where $d=0$, $\vec\bG=(\bG^0=\bG)$,
$\vec\rtm=(0)$, $\vec\phi=(1)$ and $\xo\in\Bd(\bG,k)$.
Then we associate the corresponding $\bK$-types and lattices as follows:
\[
\rK_\datum^+=\G_{\xo,0^+},\quad\phi_\datum=1,\quad\gss_\datum=\gss=0,
\]
and for $x\in\Bd(\bG,k)$, independent of the choice of $\epsilon$,
\[
\quad \CaL_x=\lieG_{x,0^+},\quad\CaL_x^\sharp=\lieG_{x,0}.
\]
We associate real numbers $\rtm_0$, $\rtm_{-1}$, 
$\stm_{-1}$ and $\stm_0$ as a zero datum as follows:
\[
\rtm_0=\rtm_{-1}=\stm_{-1}=\stm_0=0.
\]
\end{definition}

\begin{defrmk}\rm
Let $\datum$ be a zero datum.
For an irreducible admissible representation $(\pi,V_\pi)$,
if $(\rK_\datum^+,\phi_\datum)<\pi$,  $\pi$ is a depth zero representation.
By convention, we will call a zero datum both \emph{good} and 
\emph{strongly good}. Hence,
by a (strongly) good $\G$-datum, we mean either a zero or 
a (strongly) good $\G$-datum of positive depth.
\end{defrmk}

\begin{defrmk}\rm \ 
\begin{enumerate}
\item
Let $\admG$ denote the set of all irreducible admissible representations.
Let $\dualG$ (resp. $\temG$) denote 
the subset of $\admG$ which consists of unitarizable 
(resp. tempered) representations of $\G$.

We remark that $\admG$ carries a natural topology defined via approximation of
matrix cofficients.
Moreover, $\dualG$ is a closed subset of $\admG$, and the subspace topology
on $\dualG$ coincides with
the usual topology on the unitary dual.
See \cite{Tadic} for more details.

\item
Let $J$ be an open compact subgroup of $\G$, and $\sigma$ 
one of its irreducible representations. For $(\pi,V_\pi)\in\admG$, 
let $V_\pi^{(J,\sigma)}$ (or simply $V_\pi^\sigma$) denote the 
$\sigma$-isotypic component in $V_\pi$. If $V_\pi^\sigma\ne0$, that is,
if $\sigma$ occurs as a subrepresentation of the restriction of $\pi$ to $J$, 
we write $(J,\sigma)<\pi$. 
\end{enumerate}
\end{defrmk}

\section{\bf Review on Plancherel formulas on $\G$ and $\lieG$}

For $(\pi,V_\pi)\in\dualG$ and $\slf\in\Ccs(\G)$, 
we define a function $\widehat{\slf}$ on $\dualG$ as
\[
\widehat{\slf}(\pi) := \Tr(\pi(\slf)) =\charpi(\slf).
\]
Then the Plancherel formula on locally compact groups (see \cite{Dx}) 
states that there is a Borel measure $d\pi$ called \emph{Plancherel measure} 
on $\dualG$ such that the first equality in the following holds:
\begin{equation}\label{eq: PF on G}
\slf(1)=\int_{\dualG}\widehat{\slf}(\pi)\,d\pi
=\int_{\temG}\widehat{\slf}(\pi)\,d\pi.
\end{equation}
From Harish-Chandra's explicit Plancherel formula 
(\cite{HC:Pl, Wald:HC}), we have the second equality.

On the other hand, regarding $\lieG$ as a topological group
with respect to addition, we can formulate the Plancherel formula on $\lieG$
as follows: there is a Borel measure on $\dualg$ 
such that for $f\in\Ccs(\dualg)$, 
\[
f(0)=\int_{\dualg}\,\Ff(\chi)\,d \chi,
\]
where $\Ff\in\Ccs(\dualg)$ is the Fourier transformation of $f$
given by $\Ff(\chi)=\int_{\dualg}\,f(Y)\,\chi(Y)\,dY.$

Recall we have the following isomorphisms:
\[
\dualg\simeq\adualg \simeq\lieG,
\]
where $\dualg$ denotes the unitary dual of $\lieG$.
The first isomorphism is from Pontrjagin duality.
We have the second isomorphism via 
an additive character and an appropriate bilinear form on $\lieG$.
When $\HypB$ holds, we have $\dualg\simeq\lieG$ 
via $\addch$ and $\bilinear$, and we can rewrite the above formula as 
\begin{equation}\label{eq: PF on lieG}
f(0)=\int_{\lieG}\Ff(X)\,dX,
\end{equation}
where $\Ff(X)=\int_\lieG f(Y)\addch(\bilinear(X,Y))\,dY$.
In the above equation (\ref{eq: PF on lieG}), 
the $\G$-invariant measure $dX$ on $\lieG$ 
should satisfy $\vol_{\lieG}(\lieG_{x,r})\vol_{\lieG}(\lieG_{x,\rtm^\ast})=1$
for all $x\in\Bd(\bG,k)$ and $\rtm\in\bbR$.

To relate Plancherel formulas on $\lieG$ and $\G$, 
let $f\in\Ccs(\lieG)$ be supported 
in a sufficiently small neighborhood of $0$. 
Then $f\circ\log$ defines a function in $\Ccs(\G)$.
Combining the Plancherel formulas (\ref{eq: PF on G}) and 
(\ref{eq: PF on lieG}), we have
\begin{equation}\label{eq: Pl}
\int_{\lieG}\widehat f(X)\,dX=f(0)=\int_{\temG}\charpi(f\circ\log)\,d\pi
\end{equation}

\noindent
{\bf 6.4. Haar measures.} 
From now on, when $\HypB$ is valid, 
we fix a Haar measure on $\lieG$ so that (\ref{eq: PF on lieG}) is valid. 
When $\Hypk$ is valid, we fix a Haar measure on $\G$ 
so that $\vol_{\G}(G_{x,f})=\vol_{\lieG}(\lieG_{x,f})$ 
for any $x$ and $f$ as in $\Hypk$.
Then, the Plancherel measure $d\pi$ in (\ref{eq: PF on G}) 
is uniquely determined with respect to this Haar measure on $\G$.

\section{\bf Decomposition of $\lieG$}

Throughout this section, we assume that $\HypGT$ is valid.
The main result in this section is Proposition \ref{prop: dec dualg}, 
where we find a spectral decomposition of (\ref{eq: PF on lieG}).
Each spectral decomposition factor is parameterized by 
an equivalence class of semisimple elements (see (\ref{def: ss equivalence})).

\begin{numbering}\label{num: goodexp2}\rm 
We restate Proposition \ref{prop: goodexp} to fit better for our purpose: 
let $\gss\in\lieG$ be a semisimple element. 
Then $\gss$ can be written as 
\[
\gss=\gss_d+\gss_{d-1}+\cdots+\gss_0+\gamma
\]
satisfying the following:
\begin{enumerate}
\item 
$\gss_i$, $i=0,1,\cdots,d-1$, is a good element of depth $-\rtm_i$
with $-\rtm_{d-1}<-\rtm_{d-2}<\cdots<-\rtm_0<0$, and
$\dth(\gamma)\ge0$. 

\item
$\gss_d=0$ or
$\gss_d\in\lieZ^d=\lieZ_\lieG$ is nonzero with $\dth(\gss_d)<-\rtm_{d-1}$.

\item
$\bG^0\subsetneq\bG^1\subsetneq\cdots\subsetneq\bG^d=\bG$
where $\bG^i=\bC_{\bG^{i+1}}(\gss_i)=\bC_\bG(\gss^i)$
with $\gss^i=\gss_d+\cdots+\gss_{i-1}+\gss_i$.
\end{enumerate}
\noindent
Let $\rtm_d:=\rtm_{d-1}$ if $\gss_d=0$, and $\rtm_d:=-\dth(\gss_d)$ 
if $\gss_d\ne0$. 
Then, $\vec\rtm:=(\rtm_0,\cdots,\rtm_d)$ and
$\vec\bG:=(\bG^0,\bG^1,\cdots,\bG^d)$ are
determined uniquely independent of the choice of good elements $\gss_i$'s. 
Note that $\dth(\gss)=-\rtm_d$. 
\end{numbering}

\begin{defrmk} \rm 
We keep the notation from (\ref{num: goodexp2}).
Write $\Gamma=\Gamma^0+\gamma$ 
where $\Gamma^0=\gss_d+\gss_{d-1}+\cdots+\gss_0$. 
Define $\lieG^0_\Gamma$ as follows:
\[
\lieG^0_\Gamma:=\lieG^0_0=\cup_{x\in\Bd(\bG^0,k)}\lieG^0_{x,0}.
\]
Since $\bG^0=\bC_\bG(\gss^0)$ is well defined 
independent of the choice of $\Gamma^0$, so is $\lieG^0_\Gamma$.
\end{defrmk}

\begin{lemma}\label{lem: equiv ss}
Let $\Gamma$ and $\Gamma'$ be two semisimple elements in $\lieG$.
\begin{enumerate}
\item
$^\G(\Gamma^0+\lieG^0_\Gamma)$ is a $\G$-domain, that is, 
an open and closed $\G$-invariant subset of $\lieG$.
\item
$^\G(\Gamma^0+\lieG^0_\Gamma)$ and $^\G(\Gamma^{\prime0}+\lieG^0_{\Gamma'})$
are either disjoint or identical.
\end{enumerate}
\end{lemma}

\proof
(1) follows from \cite[(5.1.4)]{Asym2}.
To prove (2), suppose 
$^\G(\Gamma^0+\lieG^0_\Gamma)\cap\,^\G(\Gamma^{\prime0}+\lieG^0_{\Gamma'})
\neq\emptyset$. Without loss of generality, we may assume
$(\Gamma^0+\lieG^0_\Gamma)\cap(\Gamma^{\prime0}+\lieG^0_{\Gamma'})
\neq\emptyset$. Then there are $X\in\lieG^0_\Gamma$ and 
$X'\in\lieG^0_{\Gamma'}$ such that $Y:=\Gamma^0+X=\Gamma^{\prime0}+X'$. 
Let $X=X_s+X_n$ 
be the Jordan decomposition of $X$ 
in $\mathfrak c_{\lieG}(\Gamma^0)$.
Similary, write $X'=X'_s+X'_n$.
Then, $(\gss^0+X_s)+X_n$ and $(\gss^{\prime0}+X'_s)+X'_n$ 
are two expressions of the Jordan decomposition of $Y$. 
By the uniqueness of Jordan decomposition, we have 
$\gss'':=\gss^0+X_s=\gss^{\prime0}+X'_s$. 
Note that $\dth(X_s)=\dth(X)\ge0$ and $\dth(X'_s)=\dth(X')\ge0$.
Applying Proposition \ref{prop: goodexp} to $\gss''$, we conclude
$\C_\G(\Gamma^{\prime\prime0})=\C_\G(\Gamma^0)=\C_\G(\Gamma^{\prime0})$.
Hence, $\lieG^0_{\Gamma''}=\lieG^0_\Gamma=\lieG^0_{\Gamma'}$,
$(\Gamma^{\prime\prime}+\lieG^0_{\Gamma''})=(\Gamma^0+\lieG^0_\Gamma)
=(\Gamma^{\prime0}+\lieG^0_{\Gamma'})$, and thus 
$^\G(\Gamma^0+\lieG^0_\Gamma)=\,^\G(\Gamma^{\prime0}+\lieG^0_{\Gamma'})$.
\qed

\begin{definition}\label{def: ss equivalence}\rm\ 
Let $\Gamma$ and $\Gamma'$ be semisimple elements in $\lieG$.
\begin{enumerate}
\item
Define $\lieG_\Gamma:=\,^\G(\Gamma^0+\lieG^0_\Gamma)$.
\item
We say that $\Gamma$ and $\Gamma'$ are \emph{equivalent}
if $\lieG_\Gamma=\lieG_{\Gamma'}$.
In that case, we write $\Gamma\sim\Gamma'$.
Let $\sgsc$ be the set of equivalence classes of semisimple elements.
\item
If $\gss\sim0$, we set $d=0$, and the associated real numbers 
$\rtm_0=\rtm_{-1}=\stm_0=\stm_{-1}=0$.
\end{enumerate}
\end{definition}

By Lemma \ref{lem: equiv ss}, the above equivalence relation is well defined.

\medskip

\noindent
{\bf Examples.} 
\begin{enumerate}
\item
Any two semisimple elements $\Gamma$ and $\Gamma'$ with $\dth(\Gamma)\ge0$
and $\dth(\Gamma')\ge0$ are equivalent. In this case, we have
$\lieG_\Gamma=\lieG_{\Gamma'}=\lieG_0$.
\item
Let $\lieT$ be a $k$-torus in $\lieG$ which splits over a tamely ramified 
extension. Any $\gss,\gss'\in\lieT$ satisfying
$\gss\equiv\gss'\pmod{\lieT_0}$ are equivalent.
\end{enumerate}

\begin{lemma}\label{lem: dec lieG}
The Lie algebra $\lieG$ is a disjoint union of 
$\lieG_\Gamma=\,^\G(\gss^0+\lieG^0_\gss)$, $\gss\in\sgsc$:
\[
\lieG=\overset{\circ}\bigcup_{\gss\in\sgsc}\,^\G(\gss^0+\lieG^0_\gss)
=\overset{\circ}\bigcup_{\gss\in\sgsc}\lieG_\Gamma.
\]
\end{lemma}

\proof
`$\supset$' is obvious. For `$\subset$', let $X\in\lieG$. 
Let $X=X_s+X_n$ be the Jordan decomposition of $X$ and $\G':=C_\G(X_s)$.
Since $X_n$ is nilpotent and $X_n\in \lieG'$,
we have $X_n\in\lieG'_{x,0}$ for some $x\in\Bd(\G')$.
Hence $X\in\lieG_{X_s}$, and `$\subset$' follows. 
\qed

\begin{proposition}\label{prop: dec dualg}
For any $f\in\Ccs(\lieG)$, the integral $\int_\lieG f(X)\,dX$ 
is decomposed as follows:
\[
\int_\lieG f(X)\,dX
=\sum_{\Gamma\in\sgsc}\int_{\lieG_\Gamma} f(X)\,dX.
\]
\end{proposition}

\proof This follows from Lemma \ref{lem: equiv ss}-(1)
and Lemma \ref{lem: dec lieG}. \qed

\begin{corollary}\label{cor: dec dualg}
Suppose $\HypB$ is valid. Then, for any $f\in\Ccs(\lieG)$,
\[
f(0)=\int_\lieG \Ff(X)\,dX
=\sum_{\Gamma\in\sgsc}\int_{\lieG_\Gamma} \Ff(X)\,dX.
\]
\end{corollary}

\medskip

\section{\bf Decomposition of $\admG$}

\begin{numbering}\rm From now on, we
fix $\epsilon$ such that $\lieG_\epsilon=\lieG_{0^+}$.
Since $\{r\in\bbR\mid\lieG_r\supsetneq\lieG_{r^+}\}$ is a discrete
subset of $\bbR$, such an $\epsilon$ exists.

The choice of such an $\epsilon$ first appears in \cite{dBH}
to treat homogeneity in the depth zero case. This setting is carried to 
\cite{Asym2}. Via this choice, we can find nice test functions
which are in a certain sense dual to the orbital integrals in which 
we are interested
(see \cite[(9.1.6)]{Asym2} or (\ref{num: test ftns})).
\end{numbering}

\begin{numbering}\rm
Recall that for each twisted Levi sequence $\vec\bG$, we fixed
embeddings of buildings in (\ref{num: embeddings}).
Except when $\bZ^0/\bZ^d$ is $k$-anisotropic, 
most objects defined here might in principle depend on the 
$\bK$-types constructed in \S5, hence on the choice of embeddings.
In the course of the proof of our main result, we show
that if $\vec\bG$ is associated to a supercuspidal representation,
then $\bZ^0/\bZ^d$ is $k$-anisotropic (Proposition \ref{prop: cusp I}). 
Hence, the ambiguity in the choice of embedding is irrelevant 
for the particular argument pursued in this paper. 
\end{numbering}

Throughout this section, we assume $\Hypk$, $\HypB$ and $\HypGT$ are valid.

The main result in this section is Proposition \ref{prop: dec dualG}, 
where we find a spectral decomposition of (\ref{eq: PF on G}).
In \S11, we will see that each decomposition factor is 
parameterized by an element of $\sgsc$.

\begin{definition}\rm\ 
Let $\datum=(\vec\bG,\xo,\vec\rtm,\vec\phi)$ be a strongly good  
$\G$-datum of positive depth, and let $\gss_\datum$ be a semisimple element
associated to $\datum$ as in (\ref{def: goodatum}). 
Let $\Gamma$ be a semisimple element in $\lieG$.
\begin{enumerate}
\item
We say that $\datum$ and $\Gamma$ are {\it associated} if 
$\gss_\datum\sim\gss$. 
\item 
$\admG_\datum(\eps):=\{(\pi,V_\pi)\in\admG\mid 
\pi>(\rL_{x,\epsilon},\chi_{\gss_\datum+n}),\ \textrm{for some }
x\in\Bd(\bG^0,k), \ n\in\lieG^0_{x,-\epsilon}\cap\nil\}$.
\end{enumerate}
\end{definition}

\begin{remarks}\rm \ 
\begin{enumerate}
\item
Note that any two semisimple elements associated to $\datum$ 
as in (\ref{def: goodatum}) are equivalent. Hence, (1) is well defined.
\item
The definition in (2) is set up to achieve the equality in 
Theorem \ref{thm: main plancherel}.
Note that $n\in\lieG^0_{x,-\epsilon}$
defines a character on $\G_{x,\epsilon}$ and $\rL_{x,\epsilon}$
because $\G_{x,\epsilon}/\G_{x,\epsilon^+}$ is abelian.
Hence, $\gss_\datum+n$ also defines a character on $\rL_{x,\epsilon}$.
\end{enumerate}
\end{remarks}

\begin{definition}\label{def: admGamma} \rm 
Let $\Gamma$ be a semisimple element in $\lieG$.
\begin{enumerate}
\item 
If $\dth(\gss)<0$, 
\[
\admG_\Gamma(\eps):=\bigcup_{\gss_\datum\sim\Gamma}\admG_\datum(\eps)
\]
where the union runs over all strongly good $\G$-datum $\datum$ 
of positive depth such that $\gss_\datum\sim\Gamma$.
\item
If $\gss\sim 0$,  define $\admG_0=\admG_\gss(\eps)$ to be
the set of all depth zero representations (this is 
independent of the choice of $\epsilon$).
\end{enumerate}
\end{definition}

\begin{lemma}\label{lem: dec admG}\
Let $\Gamma$ and $\Gamma'$ be semisimple elements in $\lieG$.
\begin{enumerate}
\item
$\admG_\Gamma(\eps)=\admG_{\Gamma'}(\eps)$ 
if and only if $\Gamma\sim\Gamma'$.
\item
$\admG_\Gamma(\eps)$ and $\admG_{\Gamma'}(\eps)$ are either disjoint or
identical.
\end{enumerate}
\end{lemma}

\proof
(1) `$\Leftarrow$' is obvious. For `$\Rightarrow$',
if $\admG_\Gamma(\eps)\cap\admG_{\Gamma'}(\eps)\ne\emptyset$, then
there is a strongly good $\G$-datum of positive depth $\datum$ such that 
$\Gamma\sim\gss_\datum\sim\Gamma'$. Hence $\Gamma\sim\Gamma'$.

(2) Note that the depth of a representation in $\admG_\Gamma(\eps)$ is 
either 0 or $-\dth(\gss)$.
Suppose $\gss\sim0$. Then, $\admG_\Gamma(\eps)=\admG_0$, and
$\admG_{\Gamma'}(\eps)$ intersects $\admG_0$ 
if and only if $\Gamma'\sim0$. Hence, 
$\admG_{\Gamma'}(\eps)=\admG_\Gamma(\eps)$. 

Now, suppose $\dth(\gss),\ \dth(\gss')<0$.
Suppose $(\pi,V_\pi)\in\admG_\Gamma(\eps)\cap\admG_{\Gamma'}(\eps)$.
Then there are two strongly good $\G$-datums of positive depth
$\datum=(\vec\bG,\xo,\vec\rtm,\vec\phi)$
and $\datum'=(\vec\bG',\xo',\vec\rtm',\vec\phi')$ such that 

$(i)$ 
$\gss_\datum\sim\Gamma$ and 
$(\pi,V_\pi)>(\rL_{x,\epsilon},\chi_{\gss_\datum+n})$
for some $x\in\Bd(\bG^0,k)$ and $n\in\lieG^0_{x,-\epsilon}$, 

$(ii)$
$\gss_{\datum'}\sim\Gamma'$ and
$(\pi,V_\pi)>(\rL_{x',\epsilon},\chi_{\gss_{\datum'}+n'})$ for some
$x'\in\Bd(\bG^{\prime 0},k)$ and $n'\in\lieG^{\prime 0}_{x',-\epsilon}$.

\noindent
Since $(\pi,V_\pi)$ is irreducible, an argument similar to the one 
in \cite[(7.2)]{MP} shows that there is $g\in\G$
such that $^g\chi_{\gss_{\datum}+n}=\chi_{\gss_{\datum'}+n'}$ on 
$^g\rL_{x,\epsilon}\cap\rL'_{x',\epsilon}$. Then, 
by Lemma \ref{lem: dual blob},
$^g(\gss_\datum+n+\CaL_{x,\epsilon}^\sharp)
\cap(\gss_{\datum'}+n'+\CaL_{x',\epsilon}^{\prime\sharp})\ne\emptyset$.
Since $n\in\lieG^0_{\gss_\datum}$ and $n'\in\lieG^0_{\gss_{\datum'}}$,
by Lemma \ref{lem: asym2appr}, 
$\gss_\datum+n+\CaL_{x,\epsilon}^\sharp\subset\,
^\G\!(\gss_\datum+\lieG^0_{\gss_\datum})$
and 
$\gss_{\datum'}+n'+\CaL_{x,\epsilon}^\sharp\subset\,
^\G\!(\gss_{\datum'}+\lieG^0_{\gss_{\datum'}})$.
Hence, 
$^\G\!(\gss_\datum+\lieG^0_{\gss_\datum})\cap\,
^\G\!(\gss_{\datum'}+\lieG^0_{\gss_{\datum'}})\ne\emptyset$,
and $\gss_\datum\sim\gss_{\datum'}$ by Lemma \ref{lem: equiv ss}-(2).
Now, we have $\admG_{\Gamma}(\eps)=\admG_{\Gamma'}(\eps)$ by (1).
\qed

\begin{definition}\rm\ 
\begin{enumerate}
\item
Let $\temG_\Gamma(\eps):=\temG\cap\admG_\Gamma(\eps)$.
\item
For any pair $\type=(J,\sigma)$ consisting 
of an open compact subgroup $J$ and its
irreducible representation $\sigma$, we define subsets
$\admG_\type$ and $\temG_\type$ of $\admG$ as follows:
\[
\admG_\type:=\{\pi\in\admG\mid\type<\pi\},\qquad\qquad
\temG_{\type}:=\temG\cap\admG_{\type}.
\]
\end{enumerate}
\end{definition}

\begin{numbering}\label{num: tem structure} 
{\bf The structure of $\temG$.}\ \rm 
The following are some definitions and remarks
regarding the structure of $\temG$. The definitions here are made
in a parallel fashion to those in \cite{BK2} regarding Bernstein center, 
and most facts here can be deduced from \cite{HC:Pl} or \cite{Wald:HC}.

\begin{enumerate}
\item
Let $\frt=(\rM,\tau)_d$ be a pair of a Levi subgroup and a discrete series
on $\rM$. We will call such a pair a \emph{discrete pair}.
We say that two discrete pairs $(\rM,\tau)_d$ and $(\rM',\tau')_d$ are
\emph{equivalent} if there is a $g\in\G$ and an unramified 
unitary character on $M$
such that $\rM=\, ^g\rM'$ and $\tau\otimes\chi\simeq\,^g\tau'$.
\item
Let $\frB_d$ be the set of equivalence classes of discrete pairs.
Then, for any tempered representation $\pi$, there is a unique class 
$(\rM,\tau)_d$ in $\frB_d$ such that $\pi$ is a subquotient of a parabolically
induced representation from $(\rM,\tau\otimes\chi)_d$ 
for some unramified unitary character of $\rM$. In this case, we will say that
$(\rM,\tau)_d$ is the \emph{discrete support of $\pi$}.
\item
For $\frt\in\frB_d$, let $\temG(\frt)$ be the set of 
irreducible tempered representations whose discrete supports are $\frt$.
Then the result of \cite{HC:Pl} (see also \cite[III.4.1]{Wald:HC}) says that
$\temG$ is the disjoint union of $\temG({\frt})$, $\frt\in\frB_d$, and
the Plancherel formula can be written as follows: for $\slf\in\Ccs(\G)$,
\[
\slf(1)=\int_{\temG}\widehat{\slf}(\pi)\,d\pi.
=\sum_{\frt\in\frB_d}\int_{\temG(\frt)}\widehat{\slf}(\pi)\,d\pi.
\]
\item
Denote $\textrm{Irr}_{\frt}$ the set of all 
irreducibly induced representations 
$\Ind_{\rP}^\G\tau$ from some $(\rM,\tau)_d$ in $\frt$. 
Then 
\[
\int_{\temG(\frt)}\widehat{\slf}(\pi)\,d\pi
=\int_{\textrm{Irr}_{\frt}}\widehat{\slf}(\pi)\,d\pi
\]
and the Plancherel measure of 
$\temG(\frt)\setminus\textrm{Irr}_{\frt}$ is zero
(see \cite[IV.2.2, IV.3]{Wald:HC}).
\item
The representations in $\temG(\frt)$ have the same cuspidal support, that is,
there is a pair $(\rM',\sigma)_c$ of a Levi subgroup $\rM'$ 
and its supercuspidal representation $\sigma$ such 
that each representation in $\temG(\frt)$ is
a subquotient of the parabolically induced representation from 
$(\rM',\sigma\otimes\chi)_d$ for some unramified character $\chi$ of $\rM'$.
Moreover, representations in $\temG(\frt)$ have the same depth (\cite{MP2}).
\end{enumerate}
\end{numbering}

We also need the following lemma
for the proof of Lemma \ref{lem: closure}.
Here, by parabolic induction, we mean unitary parabolic induction.

\begin{lemma}\label{lem: tem structure}
Let $\frt\in\frB_d$. 
Let $(\rM_1,\tau_1)_d$ and $(\rM_2,\tau_2)_d$ be in the class of $\frt$.
Then, any parabolically induced representations 
$\Ind_{\rP_1}^\G\tau_1$ and $\Ind_{\rP_2}^\G\tau_2$,
where $P_i$, $i=1,2$ is a parabolic subgroup with Levi factor $M_i$,
are isomorphic when restricted to any 
special maximal compact subgroup.
\end{lemma}

\proof
Let $K$ be a special maximal compact subgroup, and 
$P_0$ a minimal parabolic subgroup of $\G$ such that $\G=P_0K=KP_0$.
Without loss of generality, we may assume that $P_1$ and $P_2$ contain $P_0$.

Denote the set of unitary unramified characters of $\rM_i$ by $X^u(\rM_i)$,
$i=1,2$.
Let $g\in G$ be 
such that $^g\!\rM_1=\rM_2$ and $^g\!\tau_1\simeq\tau_2\otimes\chi_2$
for some $\chi_2\in X^u(\rM_2)$.

Suppose $\tau_1$ is $\G$-regular (that is, $\tau_1\not\simeq\,^h\!\tau_1$ 
for any $h$ normalizing $\rM_1$). 
Then, $\tau_2\otimes\chi_2$ is also regular and 
$\Ind_{\rP_1}^\G\tau_1\simeq 
\Ind_{^g\!\rP_1}^\G\tau_2\otimes\chi_2
\simeq \Ind_{\rP_2}^\G\tau_2\otimes\chi_2$.
The second isomorphism follows from Frobenius reciprocity 
and \cite[III.7.3]{Wald:HC}.
Since $\left(\,\Ind_{\rP_2}^\G\tau_2\right)\mid K\,
\simeq\,\left(\Ind_{\rP_2}^\G\tau_2\otimes\chi_2\right)\mid K$
by \cite[(3.1.1)]{Cas}, the lemma follows.

In general, there is an unramified unitary character $\chi$ of $\rM_1$
such that $\tau_1\otimes\chi$ is $\G$-regular
This follows from the fact 
that $\{\chi\in X^u(\rM_1)\mid\tau_1\otimes\chi\simeq\tau_1\}$
is a finite set while $X^u(\rM_1)$ forms a complex manifold.
Using \cite[(3.1.1)]{Cas}, the general case follows from 
the $\G$-regular case.
\qed

\begin{lemma}\label{lem: closure}
Let $\gss\in\sgsc$ and $\frt\in\frB_d$.
Then, $\temG_\Gamma(\eps)\cap\temG(\frt)$ either contains 
$\mathrm{Irr}_{\frt}$, or has Plancherel measure zero.
\end{lemma}

\proof
Suppose $\gss\sim0$ and thus $\temG_\Gamma(\eps)$ is 
the set of depth zero representations.
By (\ref{num: tem structure})-(5), 
$\temG_\Gamma(\eps)\cap\temG(\frt)$ is either $\temG(\frt)$ or empty.

Suppose $\gss\not\sim0$.
Note that for any $\datum$, $\rL_{x,\eps}$ is contained in a special 
maximal compact subgroup. Then, by Lemma \ref{lem: tem structure},
if there is an irreducibly induced representation
$\Ind_{\rP}^\G\tau$ in $\temG_\Gamma(\eps)\cap\temG(\frt)$ for some
$(\rM,\tau)_d$ in $\frt$, 
$\textrm{Irr}_{\frt}$ is contained in $\temG_\Gamma(\eps)\cap\temG(\frt)$.
Otherwise, $\temG_\Gamma(\eps)\cap\temG(\frt)$ is contained in 
$\temG(\frt)\setminus\textrm{Irr}_{\frt}$ 
whose Plancherel measure is zero.
\qed

\medskip

In the following proposition, we denote the union of 
$\temG(\frt)$, $\frt\in\frB_d$ with
$\textrm{Irr}_{\frt}\subset\temG_\Gamma(\eps)\cap\temG(\frt)$
by $\overline{\temG_\gss}(\eps)$.

\begin{proposition} \label{prop: dec dualG} \ 

\begin{enumerate}
\item
For $\Gamma\in\sgsc$ and $\slf\in\Ccs(\G)$, we have
\[
\int_{\overline{\temG_\Gamma}(\eps)}\widehat\slf(\pi)\,d\pi
=\int_{\temG_\Gamma(\eps)}\widehat\slf(\pi)\,d\pi.
\]
\item
For $\slf\in\Ccs(\G)$, we have
\[
\slf(1)=\int_{\temG}\widehat\slf(\pi)\,d\pi=
\int_{\temG\setminus\lp\overset\circ\bigcup_{\Gamma\in\sgsc}
\temG_\Gamma(\eps)\rp}
\widehat\slf(\pi)\,d\pi \,+\,
\sum_{\Gamma\in\sgsc}\int_{\temG_\Gamma(\eps)}\widehat\slf(\pi)\,d\pi.
\]
\item
Each $\temG_\Gamma(\eps)$ has a finite Plancherel volume, 
that is, $\vol_{\dualG}(\temG_\Gamma(\eps))<\infty$.
\end{enumerate}
\end{proposition}

\begin{remark}\rm
Note that $\sgsc$ is a countable set. Moreover, for any $\slf\in\Ccs(\G)$,
$\int_{\temG_\Gamma(\eps)}\widehat\slf(\pi)\,d\pi\neq0$ for only finitely
many $\Gamma\in\sgsc$. Hence, the $\sum$ in (2) is well defined.
\end{remark}

\proof
(1) This follows from Lemma \ref{lem: closure} and 
(\ref{num: tem structure})-(4).

(2) We have 
\[
\begin{array}{ll}
\int_{\temG}\widehat\slf(\pi)\,d\pi
&\overset{(i)}=
\int_{\temG\setminus\lp\overset\circ\bigcup_{\Gamma\in\sgsc}
\overline{\temG_\Gamma}(\eps)\rp}
\widehat\slf(\pi)\,d\pi \,+\,
\sum_{\Gamma\in\sgsc}
\int_{\overline{\temG_\Gamma}(\eps)}\widehat\slf(\pi)\,d\pi
\\
&\overset{(ii)}=\int_{\temG\setminus\lp\overset\circ\bigcup_{\Gamma\in\sgsc}
\temG_\Gamma(\eps)\rp}
\widehat\slf(\pi)\,d\pi \,+\,
\sum_{\Gamma\in\sgsc}\int_{\temG_\Gamma(\eps)}\widehat\slf(\pi)\,d\pi.
\end{array}
\]
The equality $(i)$ follows from (1) and the decomposition 
in (\ref{num: tem structure})-(3).
Since $(\overline{\temG_\gss}(\eps)\cup\temG_\gss(\eps))\setminus
(\overline{\temG_\gss}(\eps)\cap\temG_\gss(\eps))$ has Plancherel measure zero
by Lemma \ref{lem: closure}, combining this with (1), the equality
$(ii)$ follows.

(3) Let $x\in\Bd(\G)$.
Since any $(\pi,V_\pi)\in\temG_\Gamma(\eps)$ is of depth $\rtm_d$,
$\temG_\Gamma(\eps)\subset\temG_{(\G_{x,\rtm_d+r'},1)}$ for a sufficiently
large $\rtm'\in\bbR$. 
Then, (3) is a result of the following lemma.
\qed

\begin{lemma}
For any open compact subgroup $J\subset\G$, 
let $\temG_{(J,1)}=\{(\pi,V_\pi)\in\temG\mid V_\pi^{(J,1)}\ne0\}$.
Then, $\temG_{(J,1)}$ has a finite Plancherel volume.
\end{lemma}

\proof 
Let $\slf\in\Ccs(\G)$ be given by 
the characteristic function of $J$ divided by $\vol_{\G}(J)$. Then,
from the Plancherel formula on $\G$, we have
\begin{center}
\quad$\vol_{\dualG}(\temG_{(J,1)})
\le\int_{\temG_{(J,1)}}\dim_{\bbC}(V_\pi^J)\,d\pi
=\int_{\temG}\dim_{\bbC}(V_\pi^J)\,d\pi  
=\int_{\temG}\widehat{\slf}(\pi)\,d\pi
=\slf(1)=\frac{1}{\vol_{\G}(J)}.$\hfill
\qed
\end{center}

\section
{\bf Review on $\gss$-asymptotic expansions}
\label{sec: gss exp}

As a preparation for the proof of Theorem \ref{thm: main plancherel}, 
we review $\gss$-asymptotic expansions and related materials
from \cite{Asym2}.

Let $\epsilon\in\bbR$ be as before.
Fix a strongly good $\G$-datum of positive depth
$\datum=(\vec\bG,\xo,\vec\rtm,\vec\phi)$, and let $\gss=\gss_\datum$.
In this case, set $\rtm_{-1}:=-\epsilon$.
Associated to $\gss$, we have a subspace $\distr^\gss$ 
of invariant distributions on $\lieG$.

\begin{numbering}\label{num: Index}\rm
Recall that $\gss^{-1}=\gss^0$. 
Let $\nil^i:=\nil\cap \lieG^i$, $0\le i\le d$, and let 
$\nil^{-1}=\nil^0$. Then $\nil^i$ is the set of nilpotent elements 
in $\lieG^i$.
Let $\Index$ be the set of all triples $(i,x,\stm)$
with $i\in \{\, 0,\dots,d\,\}$, $x\in \Bd(\bG^i,k)$, and
$\rtm_{i-1}\le \stm< \rtm_i$ if $i<d$, 
$\rtm_{d-1}\le \stm < \infty$ if $i=d$:
\[
\begin{array}{ll}
\Index&:=\{\, (i,x,\stm)\in {\bbZ} \times \Bd(\bG,k)
\times\tbR\ | \ 
i\in \{\,  0,\dots,d\,\},\ x\in  \Bd(\bG^i,k),  \\
&\qquad\qquad\qquad\qquad\qquad
\rtm_{i-1}\le \stm < \rtm_{i}\textrm{ if }i\not=d, \ 
\rtm_{d-1}\le \stm< \infty\textrm{ if }i=d\,\}.
\end{array}
\]
Let $\distr(\lieG)$ denote the set of $\G$-invariant distributions on $\lieG$.
In the following, we define a subspace $\distr^\gss$ of $\distr(\lieG)$.
\end{numbering}

\begin{definition}\label{def: distr}\rm \ 

\begin{enumerate}
\item
Let $(i,x,\stm)\in \Index$ with $s\le \rtm_{d-1}$.
Recall that $\CaL^i_x$ is defined in (4.3.1), and
define $\distr^{i,\gss}_{x,-\stm}$ as follows: 
If $\rtm_{i-1}< s<\rtm_i$, let
\[
\qquad
\begin{array}{ll}
\distr_{x,-\stm}^{i,\gss}:=&
\left\{\rT\in \distr(\lieG)\ | 
\ \textup{if}\ f\in C\lp\lp\gss^i + \lieG_{x,-\stm}^i +\CaL_{x}^{i\sharp}\rp
\left/\CaL_x^{i\sharp}\rp \ \textup{and}\right.\right.\\
&\qquad\qquad\qquad\qquad
\left.\supp(f)\cap\lp\gss^i + \lieG_{(-\stm)^+}^i\rp=\emptyset,\ 
\textup{then}\ \rT(f)=0\,
\right\}.
\end{array}
\]
Let
\[
\qquad\quad
\begin{array}{ll}
\distr_{x,-\rtm_{i-1}}^{i,\gss}:=&
\left\{ \rT\in \distr(\lieG)\ \left| \ 
\textrm{if}\ 
f\in 
\C\lp\lp\gss^i+\lieG_{x,-\rtm_{i-1}}^i+\CaL_{x}^{i\sharp}\rp
\left/\CaL_x^{i\sharp}\right.\rp\ 
\textrm{and}\ \right.\right. \\
&\qquad\qquad\quad
\left.\supp(f)\cap
{}^{\G^i}\!\lp\gss^{i-1}+\lieG^{i-1}_{(-\rtm_{i-1})^+}\rp=\emptyset,\ 
\textrm{then}\ \rT(f)=0\, \right\}.
\end{array}
\]
\item
If $(i,x,\stm)\in \Index$ and
$\stm>\rtm_{d-1}$, then $i=d$. In this case, let 
\[
\begin{array}{ll}
\distr^{d,\gss}_{x,-\stm}:=&
\left\{\, \rT\in \distr(\lieG)\ \left| \ 
\textrm{If}\ f\in C((\gss^d + \lieG_{x,-\stm})/\lieG_{x,-\stm_{d-1}})
\ \textrm{and}\ \right.\right.\\
&\qquad\qquad\qquad
\left.\supp(f)\cap (\gss^d + \lieG_{(-\stm)^+})=\emptyset,\ 
\textrm{then}\ \rT(f)=0\,
\right\}.
\end{array}
\]
\item
Set
\[
\distr^\gss := \bigcap_{(i,x,\stm)\in \Index}\, \distr_{x,-\stm}^{i,\gss}\ .
\]
\end{enumerate}
\end{definition}

For any $\tm\in\tbR$, define a subspace $D_\tm$ of $\Ccs(\lieG)$ as follows:
\[
D_{\tm}:=\sum_{x\in\Bd(\bG,k)} C_c(\lieG/\lieG_{x,\tm}).
\]
Note that $D_{\tm}$ is the space consisting of Fourier transforms of 
$f\in\Ccs(\lieG_{\tm^\ast})$ (\cite[(4.2.3)]{AD}).

\begin{theorem}\label{thm: final dec}\cite[(8.1.1)]{Asym2}
Suppose $\Hypk$ and $\HypN$ are valid. Let $\gss$ be as above.
Let $\distr_{\orb(\gss)}$ be the linear span of 
orbital integrals $\mu_\orb$ associated to $\orb$ with $\orb\in\orb(\gss)$.
Then, we have $\distr^\gss\equiv \distr_{\orb(\gss)}$ on $D_{-\stm_{d-1}}$.
\end{theorem}

For $(\pi,V_\pi)\in\admG$, we define the Fourier transform 
$\Fcharpi\in\distr(\lieG)$ of $\charpi$ as follows: for $f\in\Ccs(\lieG)$,
let $\resf= f\cdot[\lieG_{0^+}]$, and
\[
\Fcharpi(f)=\charpi(\widehat{\resf}\circ\log).
\]

\begin{theorem} \label{thm: stronger} \cite[(8.2.3)]{Asym2} 
Let $(\pi,V_\pi)\in\admG$.
Let $\datum=(\vec\bG,\xo,\vec\rtm,\vec\phi)$ be a strongly 
good $\G$-datum of positive depth.
Assume the hypotheses $\HypB$, $\Hypk$ and $\HypN$ are valid.
Let $x\in\Bd(\bG^0,k)$ and $\gss:=\gss_\datum$. 
Let $\CaL_{x}$ and $\rL_{x}$ be as in (\ref{num: aux lattices}).
Let $X\in\lieG^0_{x,-\epsilon}\cap\lieG^0_{(-\epsilon)^+}$, and
let $\chi_{\gss+X}$ be the character of $\rL_{x}$ with 
the dual blob $\gss+X+\CaL^{\sharp}_{x}$.
Suppose $(\pi,V_\pi)$ contains $(\rL_{x},\chi_{\gss+X})$. Then,
\begin{enumerate}
\item
$\Fcharpi\in\distr^{-\gss}$,
\item
$\Fcharpi$ is $\gss$-asymptotic on $\lieG_{\stm_{d-1}^+}$, that is,
there are $c_{\orb}(\pi)\in\bbC$ indexed by $\orb(\gss)$ such that
$\charpi(f\circ\log)
=\underset{\orb\in\orb(\gss)}\sum c_{\orb}(\pi)\cdot\widehat{\mu_{\orb}}(f)$ 
for all $f\in\Ccs(\lieG_{\stm_{d-1}^+})$, 
where $\mu_\orb$ is the orbital integral associated to $\orb\in\orb(\gss)$.
\end{enumerate}
\end{theorem}

\begin{numbering}\label{num: test ftns}\rm 
Recall $\bG^0=\bC_{\bG}(\gss)$. 
Label the elements of $\orb(\G^0,0)$, the set of nilpotent orbits in 
$\lieG^0$, as $\orb_1^0,\dots,\orb_m^0$.
Because $\gss$ is semisimple, the Jordan decomposition gives a bijection
between $\orb(\gss)$ and $\orb(\G^0,0)$ (see (\ref{num: G'G2})). 
If $\orb^0\in \orb(\G^0,0)$, the corresponding element of
$\orb(\gss)$ is ${}^\G(\gss + \orb^0)$.
Let $\orb_i={}^\G(\gss + \orb_i^0)$, $1\le i\le m$. 

In \cite{Asym2}, we also found some test functions 
$f^\gss_i$, $\orb_i\in\orb(\gss)$, such that $(i)$ the obvious pairing on 
$\langle \mu_{\orb_i}\mid\orb_i\in\orb(\gss)\,\rangle\times
\langle\,f^\gss_i\mid\orb_i\in\orb(\gss)\,\rangle$
is nondegenerate, and 
$(ii)$ for each $\orb_i\in\orb(\gss)$, there are $x_i\in\Bd(\bG^0,x)$
and $X_i\in\orb_i\cap\lieG^0_{x_i,(-\epsilon)}$ with 
$f^\gss_i=[\gss+X_i+\CaL^\sharp_{x_i}]$. 

This provides a way to compute the coefficients $c_\orb(\pi)$ 
in Theorem \ref{thm: stronger} by linear algebra. 
Moreover, if $\rT\,\in\distr^\gss$, then  the restriction of 
$\rT$ to $\D_{-\stm_{d-1}}$ is completely determined by
$\rT(f^\gss_i)$, $i=1,\cdots,m$, in the following sense.
\end{numbering}

\begin{theorem} \label{thm: test ftns} \cite[(9.1.6)]{Asym2}
Suppose $\Hypk$ and $\HypN$ are valid.
Let $\rT\in \distr^\gss$. Let $f^\gss_i$ be as in (\ref{num: test ftns}).
If $\rT(f^\gss_i)=0$ for all $i=1,\cdots,m$, 
then $\rT\equiv0$ on $\D_{-\stm_{d-1}}$.
\end{theorem}

\begin{corollary}\label{cor: test ftns}
Suppose $\Hypk$ and $\HypN$ are valid.
Let $\rT_1,\ \rT_2\in\distr^{\sgss}$. If \ 
$\rT_1([\sgss+\Xn+\CaL_{x}^\sharp])
=\rT_2([\sgss+\Xn+\CaL_{x}^\sharp])$ \ 
for each $x\in\Bd(\bG^0,k)$ and $\Xn\in\nil^0\cap\lieG^0_{x,(-\epsilon)}$,
then $\rT_1\equiv\rT_2$ on $\D_{-\stm_{d-1}}$.
\end{corollary}

\proof
Since $\left\{f^\gss_i\right\}_{i=1}^m\subset 
\left\{[\sgss+\Xn+\CaL_{x}^\sharp]\mid x\in\Bd(\bG^0,k),\ 
\Xn\in\nil^0\cap\lieG^0_{x,(-\epsilon)}\right\}$, the corollary follows from 
the above theorem.
\qed

\section{\bf Matching}\label{sec: match}

In the following theorem, we match the spectral decomposition factors
in (\ref{cor: dec dualg}) and (\ref{prop: dec dualG}).
This is a crucial step in proving Theorem \ref{thm: tem exhaustion}. 
Recall that if $\gss\sim0$, $\stm_{d-1}=\stm_{-1}=0$.

\begin{theorem}\label{thm: main plancherel}
Suppose $\Hypk$, $\HypB$, $\HypGT$ and $\HypN$ are valid.
Let $\sgss\in\sgsc$. 
Then for any $f\in\Ccs(\lieG_{\stm_{d-1}^+})$, we have
\[
\int_{\lieG_{\sgss}}\Ff(X)\,dX
=\int_{\temG_\sgss(\eps)} \widehat{f\circ\log}(\pi)\,d\pi.
\]
\end{theorem}

\proof
Suppose $\sgss\sim0$. Then $\lieG_{\sgss}=\lieG_0$ and 
$\temG_\sgss(\eps)$ is the set of all tempered depth zero representations.
Hence the above equality follows from \cite[(3.3.1)]{Dpl1}.

Now assume that $\dth({\sgss})<0$. Without loss of generality,
we may assume $\sgss=\gss^0$. Write $\sgss$ 
as in (\ref{num: goodexp2}):
\[
\sgss=\gss^0=\gss_d+\gss_{d-1}+\cdots+\gss_0.
\]
Let $\rT_{\sgss}^{\lieG}$ and $\rT_{\sgss}^{\G}$ be invariant distributions
defined as follows: for $f\in\Ccs(\lieG)$, let $\resf:=f\cdot[\lieG_{0^+}]$. 
Then,
\begin{align*}
\rT_{\sgss}^{\lieG}(f)
&:=\int_{\lieG_{\sgss}}\widehat{\resf}(X)\,dX , \\
\rT_{\sgss}^{\G}(f)
&:=\int_{\temG_{\sgss}(\eps)} \widehat{\resf\circ\log}(\pi)\,d\pi
=\int_{\temG_{\sgss}(\eps)} \charpi(\resf\circ\log)\,d\pi.
\end{align*}
Hence, to prove the theorem, it is enough to prove that
\[
\rT_{\sgss}^{\lieG}(\Ff)=\rT_{\sgss}^{\G}(\Ff)
\]
for any $f\in\D_{-\stm_{d-1}}$, 
since $\widehat\D_{-\stm_{d-1}}=\Ccs(\lieG_{\stm_{d-1}^+})$.
We first need the following lemma.

\begin{lemma} \label{lem: here}
$\widehat \rT_{\sgss}^{\lieG}$ and $\widehat \rT_{\sgss}^{\G}$ are 
elements of $\distr^{-\gss}$. 
\end{lemma}

\proof
By Theorem \ref{thm: stronger}-(1),
we have $\widehat \rT_{\sgss}^{\G}\in\distr^{-{\sgss}}$.
To show $\widehat \rT_{\sgss}^{\lieG}\in\distr^{-{\sgss}}$, note that
for any $f\in\Ccs(\lieG)$, we have
\[
\widehat \rT_{\sgss}^{\lieG}(f)=\int_{\lieG_{\sgss}} f(-X)\,dX
=\int_{-\lieG_{\sgss}} f(X)\,dX,
\]
and 
\[
\supp(\widehat \rT_{\sgss}^{\lieG})=-\lieG_{\sgss}=\lieG_{-\sgss}.
\]
If $\stm>\rtm_{d-1}$, then
\[
\supp(\widehat \rT_{\sgss}^{\lieG})=-\lieG_{\sgss}=\,^\G(-\sgss+\lieG^0_0)
\,\subset\,
^\G(-\sgss+\lieG^d_{(-\stm)^+})=\,^\G(-\gss_d+\lieG^d_{(-\stm)^+})
=-\gss_d+\lieG^d_{(-\stm)^+}.
\]
Hence $\widehat \rT_{\sgss}^{\lieG}\in\distr^{d,-\sgss}_{x,-\stm}$.
If $\rtm_{i-1}<\stm\le\rtm_i$, $i=0,\cdots,d$, then
\[
\supp(\widehat \rT_{\sgss}^{\lieG})=\,^\G(-\sgss+\lieG^0_0)\,\subset\,
^\G(-\sgss+\lieG^i_{(-\stm)^+})=\,^\G(-\gss^i+\lieG^i_{(-\stm)^+}),
\]
which implies that $\widehat \rT_{\sgss}^{\lieG}
\in\distr^{i,-\sgss}_{x,-\stm}$.
Therefore, 
$\widehat \rT_{\sgss}^{\lieG}
\in\bigcap_{(i,x,s)\in\Upsilon_\sgss}\distr^{i,-\sgss}_{x,-\stm}
=\distr^{-\gss}$.
\qed

\medskip

Continuing with the proof of Theorem \ref{thm: main plancherel}, 
by Lemma \ref{lem: here} and Corollary \ref{cor: test ftns}, 
it is enough to check
that for any $x\in\Bd(\bG^0,k)$ and $\Xn\in\nil^0\cap\lieG^0_{x,-\epsilon}$,

\begin{center}
(E)\hfill 
$\widehat \rT_{\sgss}^{\lieG}\lp[-{\sgss}+\Xn+\CaL_x^\sharp]\rp
=\widehat \rT_{\sgss}^{\G}\lp[-{\sgss}+\Xn+\CaL_x^\sharp]\rp.$
\qquad\qquad\qquad\qquad\qquad\qquad\qquad\qquad
\end{center}

\noindent
We first note that $\supp\lp[-{\sgss}+\Xn+\CaL_x^\sharp]\rp
\subset\,^\G(-\sgss+\Xn+\lieG^0_{0})
\,\subset\, ^\G(-\sgss+\lieG^0_{0}) \,=\,\lieG_{-\sgss}$.
The first inclusion follows from Lemma \ref{lem: asym2appr}, and
the second inclusion follows from $X$ being nilpotent and thus $X\in\lieG^0_0$.
Then we have 

\begin{center}
(a) \hfill
$\widehat \rT_{\sgss}^{\lieG}\lp[-{\sgss}+\Xn+\CaL_x^\sharp]\rp
=\vol_{\lieG}(\CaL_x^\sharp)=\frac{1}{\vol_{\lieG}(\CaL_x)}.$
\qquad\qquad\qquad\qquad\qquad\qquad\qquad\qquad
\end{center}

\noindent
Computing $\widehat \rT_{\sgss}^{\G}\lp[-{\sgss}+\Xn+\CaL_x^\sharp]\rp$,
write $f$ for $[-{\sgss}+\Xn+\CaL_x^\sharp]$. Then,

\smallskip

\qquad\qquad\qquad\qquad\qquad
$\Ff(Y)=\frac1{\vol_{\G}(\rL_x)}\cdot
  \addch(\bilinear(-{\sgss}+\Xn, Y))\cdot [\CaL_x](Y).$

\smallskip

\noindent
Note that $\lp\vol_{\G}(\rL_x)\cdot\Ff\circ\log\rp$ 
is a character of $\rL_x$ with dual blob $-{\sgss}+\Xn+\CaL_x^\sharp$,
that is, 
$\lp\vol_{\G}(\rL_x)\cdot\Ff\circ\log\rp\,
=\,\chi_{_{-{\sgss}+\Xn}}$ on $\rL_x$.
Then $\charpi(\Ff\circ\log)=\charpi(\chi_{_{-{\sgss}+\Xn}})$ 
is the multiplicity 
$\mult(\chi_{_{{\sgss}-\Xn}},V_\pi)$ of $\chi_{_{{\sgss}-\Xn}}$ in $V_\pi$. 
Since $\mult(\chi_{_{{\sgss}-\Xn}},V_\pi)=0$ 
unless $\pi\in\temG_{\sgss}(\eps)$, 
we have

\medskip

\noindent
(b)  \qquad
$\begin{array}{ll}
\frac{1}{\vol_{\lieG}(\CaL_x)}
&=\Ff(\log(1))=\int_{\temG}\charpi(\Ff\circ\log)\,d\pi
=\int_{\temG}\mult(\chi_{_{{\sgss}-\Xn}},V_\pi)\,d\pi
=\int_{\temG_\sgss(\eps)}\mult(\chi_{_{{\sgss}-\Xn}},V_\pi)\,d\pi \\
&=\int_{\temG_{\sgss}(\eps)}\charpi(\Ff\circ\log)\,d\pi
=\widehat \rT_{\sgss}^{\G}\lp[-{\sgss}+\Xn+\CaL_x^\sharp]\rp.
\end{array}$

\medskip
\noindent
Now the equality (E) follows from (a) and (b), and Theorem 
\ref{thm: main plancherel} is proved.
\qed

\smallskip

\section{\bf Tempered representations}\label{sec: tem}

Throughout this section, we assume $\Hypk$, $\HypB$, $\HypGT$ 
and $\HypN$ are valid.

\smallskip

We prove that for almost every tempered representation $(\pi,V_\pi)$, 
there is a strongly good $\G$-datum $\datum$ such that
$(\pi,V_\pi)$ contains the $\bK$-type $(\rK^+_\datum,\phi_\datum)$.
As a corollary of the proof, we get a spectral decomposition of 
the delta distribution on $\G$, where each decomposition factor
is parameterized by an element of $\sgsc$.

\begin{theorem} \label{thm: tem exhaustion}
Almost every tempered representation is an element of $\temG_\sgss(\eps)$ 
for some $\sgss\in\sgsc$. 
That is, 
$\temG\setminus\overset\circ\bigcup_{\sgss\in\sgsc}\,\temG_\sgss(\eps)$
has Plancherel measure zero.
\end{theorem}

In the following, for $\CaA$, $\CaA'\subset\temG$, if
$(\CaA\cup\CaA')\setminus(\CaA'\cap\CaA)$ has measure zero,
we will write $\CaA\overset\mu=\CaA'$.

\proof
We fix $\xx\in\Bd(\bG,k)$. For $\itt\in\bbZ_{>0}$,
let $\slft:=[\G_{\xx,\itt^+}]$. Then, for any $(\pi,V_\pi)\in\admG$,
we have $\widehat\slft(\pi)=\Fcharpi([\G_{\xx,\itt^+}])
=\vol_\G(\G_{\xx,\itt^+})\dim V_\pi^{(\G_{\xx,\itt^+},1)}\ge0$.
Let $\temG{}^\itt:=\temG_{(\G_{\xx,\itt^+},1)}$, 
and let $\temG_\sgss(\eps){}^\itt:=\admG_\sgss(\eps)\cap\temG{}^\itt$.

Computing
\begin{align*}
&\int_\lieG\widehat{\slft\circ\exp}(X)\,dX
=\int_{\temG}  \widehat\slft(\pi)\,d\pi 
=\int_{\temG{}^\itt}  \widehat\slft(\pi)\,d\pi 
\overset{\mathit{(i)}}
\ge\sum_{\sgss\in\sgsc}\int_{\temG_\sgss(\eps)^\itt} \widehat\slft(\pi)\,d\pi\\
&=\sum_{\sgss\in\sgsc,\,\temG_\sgss(\eps)^\itt\ne\emptyset}
\int_{\temG_\sgss(\eps)^\itt} \widehat\slft(\pi)\,d\pi
\overset{\mathit{(ii)}}
=\sum_{\sgss\in\sgsc,\,\temG_\sgss(\eps)^\itt\ne\emptyset}
\int_{\lieG_{\sgss}} \widehat{\slft\circ\exp}(X)\,dX\\
&\overset{\mathit{(iii)}}
=\sum_{\sgss\in\sgsc}\int_{\lieG_{\sgss}}\widehat{\slft\circ\exp}(X)\,dX
\overset{\mathit{(iv)}}
=\int_\lieG\widehat{\slft\circ\exp}(X)\,dX.
\end{align*}
The relations $(i),\ (iii)$ and $(iv)$ are straightforward.
We use Theorem \ref{thm: main plancherel} to verify $(ii)$. 

The first inequality $\mathit{(i)}$ follows from 
$\overset\circ\bigcup_{\sgss\in\sgsc}\temG_\sgss(\eps)^\itt\subset\temG{}^\itt$
and $\widehat\slft(\pi)\ge0$.

To prove $\mathit{(ii)}$, it is enough to verify that for $\sgss\in\sgsc$,
if $\temG_\sgss(\eps)^\itt\ne\emptyset$, 
\[
\int_{\temG_\sgss(\eps)^\itt} \widehat\slft(\pi)\,d\pi
=\int_{\temG_\sgss(\eps)} \widehat\slft(\pi)\,d\pi
=\int_{\lieG_{\sgss}} \widehat{\slft\circ\exp}(X)\,dX.
\]
Since for $(\pi,V_\pi)\in\temG_\sgss(\eps)$, $\widehat\slft(\pi)\ne0$ 
only if $(\pi,V_\pi)\in\temG_\sgss(\eps)^\itt$, the first equality
follows. To prove the second one,
we note that $\widehat\slft(\pi)\ne0$ only if $\itt\ge\dpi(\pi)=\rtm_{d}$
for $\pi\in\temG_\sgss(\eps)$. Hence
$\temG_\sgss(\eps)^\itt\ne\emptyset$ only if $\itt\ge\rtm_d$.
In that case, since 
$\slft\circ\exp\in\Ccs(\lieG_{\rtm_d^+})
\subset\Ccs(\lieG_{\stm_{d-1}^+})$, 
the equality in Theorem \ref{thm: main plancherel} is valid for $\slft$.
Hence the above equalities are valid, and so is $\mathit{(ii)}$. 

To prove $\mathit{(iii)}$, we note that $\widehat{\slft\circ\exp}$ is
a scalar multiple of the characteristic function of $\lieG_{\xx,-\itt}$.
Then since $\supp(\widehat{\slft\circ\exp})\cap\lieG_{\sgss}
=\lieG_{\xx,-\itt}\cap\lieG_{_\sgss}\ne\emptyset$
only if $\temG_\sgss(\eps)^\itt\ne\emptyset$, 
the equality $\mathit{(iii)}$ holds.
The equality $\mathit{(iv)}$ follows from Lemma \ref{lem: dec lieG}.

Hence, $\ge$ in $(i)$ is in fact an equality for all $\itt\in\bbZ_{>0}$.
Therefore, $\temG{}^\itt\overset\mu=
\bigcup_{\sgss\in\sgsc}\temG_\sgss(\eps)^\itt$.
Since $\temG=\bigcup_{\itt\in\bbZ_{>0}}\temG{}^\itt$ 
and $\temG_\sgss(\eps)=\bigcup_{\itt\in\bbZ_{>0}}\temG_\sgss(\eps)^\itt$,
we can now conclude

\smallskip\

\noindent
$(\ast)$ \hspace{5cm}
$\temG 
\overset\mu=\bigcup_{\sgss\in\sgsc}\,\temG_\sgss(\eps)\ .$

\smallskip
\noindent
Combining this with Lemma \ref{lem: dec admG}, the above is a disjoint union.
\qed

\begin{corollary}
For any $f\in\Ccs(\G)$, we have
\[
f(1)=\sum_{\sgss\in\sgsc}\int_{\temG_\sgss}\charpi(f)\,d\pi.
\]
\end{corollary}

\proof 
This follows from Proposition \ref{prop: dec dualG}-(1) 
and $\temG\,\overset\mu=\,\overset\circ\bigcup_{\sgss\in\sgsc}\temG_\sgss$.
\qed

\begin{remark}\label{rmk: indep epsilon}\rm 
Note that $\temG\overset\mu
=\overset\circ\bigcup_{\gss\in\sgsc}\temG_\sgss(\epsilon)$ 
independent of the choice of $\epsilon>0$ such that 
$\lieG_\epsilon=\lieG_{0^+}$. 
Since $\temG_\sgss(\epsilon')\subset\temG_\sgss(\epsilon)$ 
for $0<\epsilon'<\epsilon$ and the above $(\ast)$ is a disjoint union,
$\temG_\sgss(\epsilon)\overset\mu=\temG_\sgss(\epsilon')$.
By Lemma \ref{lem: closure}, we also have 
\[
\bigcup{}_\epsilon\,\temG_\sgss(\epsilon)\,\overset\mu=\,
\bigcap{}_\epsilon\,\temG_\sgss(\epsilon)
\]
where $\epsilon$ runs over real numbers such that 
$\lieG_\epsilon=\lieG_{0^+}$, and thus
\[
\temG\,\overset\mu=\,\overset\circ\bigcup_{\sgss\in\sgsc}
\left(\bigcap{}_\epsilon\,\temG_\sgss(\epsilon)\right).
\]
\end{remark}

\begin{theorem}\label{thm: tempered}
For $(\pi,V_\pi)\in\overset\circ\bigcup_{\sgss\in\sgsc}
\left(\bigcap{}_\epsilon\,\temG_\sgss(\epsilon)\right)$, 
there is a strongly good $\G$-datum $\datum$ such that
$(\pi,V_\pi)$ contains $(\rK_\datum^+,\phi_\datum)$.
\end{theorem}

\proof
This is obvious if $\dpi(\pi)=0$. 
Now, suppose $\dpi(\pi)>0$. Let $\sgss\not\sim0$ be such that 
$(\pi,V_\pi)\in\bigcap{}_\epsilon\,\temG_\sgss(\epsilon)$, and
it splits over a tamely ramified extension $E$.
Let $\vec\bG$, $\vec\rtm$ and $\vec\phi$ be the sequences
associated to $\sgss$.
Fix a subset $\chamberD$ in $\Bd(\bG^0,k)$ such that 
$\chamberD$ is compact and $\G^0\cdot\chamberD=\Bd(\bG^0,k)$.
If $\bG^0$ is semisimple, the closure of a chamber in $\Bd(\bG^0,k)$ 
is such a subset. Otherwise, since $\Bd(\bG^0,k)$ is the product of 
$\Bd(\bG^{0\der},k)$ and $\bX_\ast(\bZ_{\bG^0},k)\otimes\bbR$,
and since $(\bX_\ast(\bZ_{\bG^0},k)\otimes\bbR)/Z_{\G^0}$ is compact,
such a subset $\chamberD$ exists.
In particular, we can choose $\chamberD$ such that 
$\chamberD\subset\Apt(\bG^0,\bT,k)$ 
for some $E$-split maximal $k$-torus $\bT$ in $\bG^0$.

Fix $\epsilon>0$ such that $\lieG_\epsilon=\lieG_{0^+}$.
Choose a decreasing sequence $\epsilon_j$, $j=1,2,\cdots,$ such that 
$\epsilon>\epsilon_j>0$ and  $\epsilon_j\rightarrow0$.

For each $\epsilon_j$, there is a strongly good $\G$-datum of positive depth
$\datum_{x_j}=(\vec\bG,x_j,\vec\rtm,\vec\phi)$ with $x_j\in\chamberD$ 
and $n_j\in\lieG^0_{x_j,-\epsilon_j}\cap\nil^0$ such that 
$(\rL_{x_j,\epsilon_j},\chi_{\sgss+n_j})<(\pi,V_\pi)$
(such an $x_j$ exists by conjugation if necessary).
Since $\chamberD$ is compact, there is
a cluster point $y\in\chamberD$ of $x_j$'s.

Let $\epsilon'>0$ be such that 
$\lieG_{y,0^+}=\lieG_{y,\alpha}=\lieG_{y,\alpha^+}$,
and $\lieG_{y,\stm_i^+}=\lieG_{y,\stm_i+\alpha}$ 
for any $\alpha<\epsilon'$.
Let $\delta$ be such that for any affine root $\psi$ in $\Psi(\bG,\bT,E)$
and $x,x'\in\Apt(\bG,\bT,E)$,
if $dist(x,x')<\delta$, $|\psi(x-x')|<\epsilon'/3$. 
Then, in that case, for any $\rtm\in\bbR$,
\[
\blieG(E)_{x,\rtm+\epsilon'/3}\overset{(\dag)}\subset\blieG(E)_{x',\rtm}
\subset\blieG(E)_{x,\rtm-\epsilon'/3}.
\]
Since $y$ is a cluster point, there are $\epsilon_j$ and $x_j$ such that 
$|y-x_j|<\delta$ and $\epsilon_j<\epsilon'/3$. 
Fix such $\epsilon_j$ and $x_j$.
Then we have the following:

({\it i}) With $\rtm=\epsilon_j$, $x=y$ and $x'=x_j$ in $(\dag)$, \ 
$\lieG_{y,0^+}=\lieG_{y,\epsilon_j+\epsilon'/3}\subset\lieG_{x_j,\epsilon_j}$,

({\it ii}) with $\rtm=\epsilon_j+\epsilon'/3$, \ 
$\lieG_{y,\epsilon_j+2/3\epsilon'}\subset\lieG_{x_j,\epsilon_j+\epsilon'/3}$.
Since $\lieG_{x_j,\epsilon_j+\epsilon'/3}\subset\lieG_{x_j,\epsilon_j^+}$
and $\lieG_{y,\epsilon_j+2/3\epsilon'}=\lieG_{y,0^+}$, we have
$\lieG_{y,0^+}\subset\lieG_{x_j,\epsilon_j^+}$ and
$\lieG_{x_j,-\epsilon_j}\subset\lieG_{y,0}$,

({\it iii}) with $\rtm=\stm_i+\epsilon'/3$, \  
$\lieG_{y,\stm_i^+}=\lieG_{y,\stm_i+2\epsilon'/3}
\subset\lieG_{x_j,\stm_i+\epsilon'/3}\subset\lieG_{x_j,\stm_i^+}$.

\noindent
From ({\it i})--({\it iii}), we see that
$\rK^+_{\datum_y}=\vec\G_{y,\vec\stm^+}\subset 
L_{x_j,\eps_j}=\vec\G_{x_j,\vec\stm^+(\eps_j)}$
and $n_j\in\lieG^0_{x_j,-\epsilon_j}\subset\lieG^0_{y,0}$. Hence,
\[
\chi_{\sgss+n_j}|\vec\G_{y,\vec\stm^{+}}=\chi_\sgss,
\]
and $(\pi,V_\pi)$ contains $(\rK^+_{\datum_y},\chi_\sgss)$.
\qed

\begin{corollary}\label{cor: tempered}
Let $(\pi,V_\pi)$ be an irreducible supercuspidal representation. 
Then, there is a strongly good $\G$-datum $\datum$ such that
$(\pi,V_\pi)$ contains $(\rK_\datum^+,\phi_\datum)$.
\end{corollary}

\proof
Suppose $(\pi,V_\pi)$ is tempered.
Then $\frt:=(\G,\pi)\in\frB_t$. For any $\gss\in\sgsc$ and
$\epsilon>0$ such that $\lieG_{0^+}=\lieG_\epsilon$, we have either
$\temG(\frt)\subset\temG_\Gamma(\eps)$ or 
$\temG(\frt)\cap\temG_\Gamma(\eps)=\emptyset$. Now, since
$\temG(\frt)$ has strictly positive Plancherel measure, 
the corollary follows from Theorem \ref{thm: tempered}.
Now, any supercuspidal representation is tempered 
up to twisting by an unramified character, and 
the general case follows from the tempered case.
\qed

\section{\bf Review of Yu's construction}\label{sec: sc}

The construction of supercuspidal representations by Yu
is based on a {\it generic $\G$-datum} (see \cite{yu} for details
and \cite{AD2} for a summary).

\begin{definition} \rm
A {\it generic $\G$-datum} is a quintuple 
$\datum_Y=(\vec\bG,\xo,\vec\rtm,\vec\phi,\rho)$ satisfying the following:

\item{${\mathbf{D}}1.$}
$\vec{\bG}=(\bG^0,\bG^1,\cdots,\bG^d=\bG)$ is a tamely ramified twisted 
Levi sequence such that $\bZ_{\bG^0}/\bZ_{\bG}$ is anisotropic.

\midvsp

\item{${\mathbf{D}}2.$}
$\xo\in\Bd(\bG^0,k)$.

\midvsp

\item{${\mathbf{D}}3.$}
$\vec\rtm=(\rtm_0,\rtm_1,\cdots,\rtm_{d-1},\rtm_d)$ is
a sequence of positive real numbers 
with $0<\rtm_0<\cdots<\rtm_{d-2}< \rtm_{d-1}\le\rtm_d$ if $d>0$,
$0\le\rtm_0$ if $d=0$.

\midvsp

\item{${\mathbf{D}}4.$}
$\vec\phi=(\phi_0,\cdots,\phi_d)$ is a sequence of quasi-characters,
where $\phi_i$ is a generic quasi-character of $\G^i$; 
$\phi_i$ is trivial on $\G^i_{\xo,\rtm_i^+}$, but
non-trivial on $\G^i_{\xo,\rtm_i}$ for $0\le i\le d-1$.
If $\rtm_{d-1}<\rtm_d$, 
$\phi_d$ is nontrivial on $\G^d_{\xo,\rtm_d}$ and trivial on 
$\G^d_{\xo,\rtm{}^+_d}$. Otherwise, $\phi_d=1$.

\midvsp

\item{${\mathbf{D}}5.$}
$\rho$ is an irreducible representation of $G^0_{[\xo]}$, 
the stabilizer in $\G^0$ of the image $[\xo]$ of $\xo$
in the reduced building of $\bG^0$, 
such that $\rho|G^0_{\xo,0^+}$ is isotrivial
and $c\textrm{-Ind}_{G^0_{[\xo]}}^{\G^0}\rho$ is irreducible and supercuspidal.
\end{definition}

\begin{remark}\rm\ 
By (6.6) and (6.8) of \cite{MP2}, 
${\mathbf{D}}5$ is equivalent to the condition that 
$\G^0_{\xo,0}$ is a maximal parahoric subgroup in $\G^0$ and
$\rho|\G^0_{\xo,0}$ induces a cuspidal representation of 
$\G^0_{\xo,0}/\G^0_{\xo,0^+}$.
\end{remark}

As remarked before, under our hypotheses, a generic quasi-character 
is also good. Instead of reviewing the definition of genericity,
we will briefly describe the construction when $\phi_i$ is good.

\begin{numbering}\label{num: notation yu}\rm 
From now on, we fix a generic $\G$-datum 
$\datum_Y=(\vec\bG,\xo,\vec\rtm,\vec\phi,\rho)$. 
Associated to  $\datum_Y$, 
we have the following open compact subgroups. Note that
their definitions depend only on $\vec\bG,\ \xo$ and $\vec\rtm$.

\begin{enumerate}
\item
$\rK^0:=\G^0_{[\xo]}$ ; $\rK^{0+}:=\G^0_{\xo,0^+}$.
\item
$\rK^i:= \G^0_{[\xo]}\G^1_{\xo,s_0}\cdots\G^i_{\xo,s_{i-1}}$ ; 
$\rK^{i+}:= \G^0_{\xo,0^+}\G^1_{\xo,s_0^+}\cdots\G^i_{\xo,s_{i-1}^+}$ \ 
for $1\le i\le d$.

\noindent
Then note that $\rK^{i+}=\rK^{i+}_\datum$ 
(see (\ref{num: construction Ktypes})).
\item 
$J^i:=(\bG^{i-1},\bG^i)(k)_{\xo,(\rtm_{i-1},\stm_{i-1})}$ ; 
$J^i_+:=(\bG^{i-1},\bG^i)(k)_{\xo,(\rtm_{i-1},\stm_{i-1}^+)}$.
\item
For $i>0$, we have 
$\rK^{i-1} J^i=\rK^i$, and $\rK^{i-1\,+} J^i_+=\rK^{i+}$.
\end{enumerate}
\end{numbering}

\begin{numbering} \label{num: yusc} {\bf Construction.} \rm 
Let $\hat\phi_i$ denote a character of $\rK^0\G^i_{\xo,0}\G_{\xo,\stm_i^+}$
defined similarly as in (\ref{num: construction Ktypes}) 
(replacing $\G^0_{\xo,0^+}$ with $\rK^0\G^i_{\xo,0}$).
For $0\le i<d$, there exists a representation $\tilde\phi_i$ of 
$\rK^i\ltimes J^{i+1}$ such that $\tilde\phi_i|J^{i+1}$ is
$\hat\phi_i|J^{i+1}_+$-isotypical and $\tilde\phi_i|\rK^{i+}$ is isotrivial.

Let $\textrm{inf}(\phi_i)$ denote the inflation of $\phi_i|\rK^i$
to $\rK^i\ltimes J^{i+1}$. Then $\mathrm{inf}(\phi_i)\otimes\tilde\phi_i$
factors through a map 
\[
\rK^i\ltimes J^{i+1}\longrightarrow \rK^i\rJ^{i+1}=\rK^{i+1}.
\]
Let $\kappa_i$ denote the corresponding representation of $\rK^{i+1}$.
Then it can be extended trivially to $\rK^d$, and we denote
the extended representation again by $\kappa_i$. We also extend $\rho$ 
trivially to a representation of $\rK^d$ and denote 
this extended representation again by $\rho$.
Define a representation $\kappa$ and $\rho_{\datum_Y}$
of $\rK^d$ as follows:
\[
\begin{array}{ll}
\kappa&:=\kappa_0\otimes\cdots\otimes\kappa_{d-1}\otimes(\phi_d|\rK^d),\\
\rho_{\datum_Y}&:=\rho\otimes\kappa.
\end{array}
\]
Note that $\kappa$ is defined only from $(\vec\bG,\xo,\vec\rtm,\vec\phi)$
independent of $\rho$.
\end{numbering}

\begin{theorem}\cite{yu}
$\pi_{\datum_Y}=c\textrm{-}\mathrm{Ind}_{\rK^d}^\G\rho_{\datum_Y}$ 
is an irreducible supercuspidal representation.
\end{theorem}

\begin{remark}\label{rmk: depth zero exhaustion}\rm
In $\datum_Y$, if $d=0$ and $\rtm_0=0$, 
$\pi_{\datum_Y}=c\textrm{-}\mathrm{Ind}_{\G_{[\xo]}}^\G\rho$ is a depth zero
supercuspial representation. 
Moreover, by \cite{Mr} or \cite[(6.6), (6.8)]{MP2},
any depth zero supercuspidal representation is of this form.
Hence, Yu's construction includes all depth zero supercuspidal 
representations.
\end{remark}

\begin{remark}\label{rmk: nec cond}\rm 
Suppose $\Hypk$, $\HypB$, $\HypGT$ and $\HypN$ are valid.
Let $(\pi,V_\pi)$ be a supercuspidal representation of positive depth. 
Then by Corollary \ref{cor: tempered}, 
there is a strongly good $\G$-datum $\datum$ of positive depth
so that $(\rK^+_\datum,\phi_\datum)<(\pi,V_\pi)$.
Comparing $\datum$ and $\datum_Y$, we note that Yu imposed certain 
sufficient conditions on $\datum_Y$ to get supercuspidal representations
such as $(i)$ $\bZ_{\bG^0}/\bZ_{\bG}$ is anisotropic, 
({\it ii}) $\G^0_{\xo,0}$ is a maximal parahoric subgroup in $\G^0$, and 
({\it iii}) $\rho$ is irreducible and induces a cuspidal representation of 
$\G^0_{\xo,0}/\G^0_{\xo,0^+}$.
Now, we will prove that those are also necessary conditions (\S\S13--17).
We start with some preparation. 
\end{remark}

\section{\bf Preparation}

To prove that the condition $(i)$ in Remark \ref{rmk: nec cond} is also 
a necessary condition for supercuspidality 
(see Proposition \ref{prop: cusp I}),
we will use some basic properties of Hecke algebras.
In this section, we review Hecke algebras
and deduce various results regarding them (in particular, 
Corollary \ref{cor: extend} and Proposition\ref{prop: support}).

\begin{numbering}\label{num: HA}
{\bf Hecke algebras.} \rm
Let $J$ be an open compact subgroup of $\G$ and let $(\sigma,V_\sigma)$ be an
irreducible representation of $J$. Denote the contragredient of $\sigma$
by $\tilde\sigma$. Then, we can associate
a Hecke algebra $\CaH(\G/\!/J,\sigma)$ to $(J,\sigma)$ as follows:
\[
\CaH(\G/\!/J,\sigma)=\{f\in\Ccs(\G,\textrm{End}_{\bbC}(V_{\tilde\sigma}))\mid 
f(jgj')=\tilde\sigma(j)f(g)\tilde\sigma(j'),
\textrm{ for } j,j'\in J,\ g\in\G\}.
\]
This is an algebra under convolution with the identity
$\frac1{\vol(J)}1_{\tilde\sigma}\cdot[J\,]$. Here, $1_{\tilde\sigma}$ is
the identity map on $V_{\tilde\sigma}$. It is well known 
there is a one-to-one correspondence between the set of simple nondegenerate
modules of $\CaH(\G/\!/J,\sigma)$ and $\admG_{(J,\sigma)}$.
For more details, see for example \cite{BK, BK2, HM}. 

We say $g\in\G$ is in the support of $\CaH(\G/\!/J,\sigma)$, if
there is an $f\in\CaH(\G/\!/J,\sigma)$ such that $g$ is in the support of
$f$. Denote the set of all such $g$ as $\supp(\CaH(\G/\!/J,\sigma))$.
\end{numbering}

\begin{center}
{\bf The case $\vec\bG=(\bG',\bG)$}
\end{center}

\begin{numbering}\label{num: prepare}\rm
Let $\vec\bG=(\bG',\bG)$ be a tamely ramified twisted Levi sequence, and
let  $\yo\in\Bd(\bG',k)$. 
Let $\phi$ be a $\G$-good character of $\G'$ such that on $\G'_{\yo,\rtm}$,
it is represented by a $\G$-good element $\gamma\in\lieZ_{\lieG'}$ of depth
$-\rtm<0$ with $\bC_\bG(\gamma)=\bG'$. Set $\stm:=\frac\rtm2$.

Let $\bT\subset\bG'$ be an $E$-split maximal $k$-torus such that 
$\yo\in\Apt(\bG,\bT,k)$. Let $\bS$ be a $k$-split subtorus in $\bT$.
Let $\bM$ be the $k$-Levi subgroup of $\bG$ defined as $\bC_\bG(\bS)$.
Since $\bT\subset\bM$, $\yo\in\Bd(\bM,k)$.
Let $\bP$ be a $k$-parabolic subgroup with Levi decomposition $\bP=\bM\bU$, 
and let $\overline\bU$ be the unipotent subgroup opposite to $\bU$.
Define $f_\tm:\Phi(\bG,\bT,E)\cup\{0\}\longrightarrow \tbR$ and
$f'_\tm:\Phi(\bG',\bT,E)\cup\{0\}\longrightarrow \tbR$ as 
\[
f_\tm(a)=\left\{
\begin{array}{ll}
\tm^+&\textrm{if }a\in\Phi(\bM,\bT,E)\cup\{0\},\\
\tm^+&\textrm{if }a\in\Phi(\overline\bU,\bT,E),\\ 
\tm&\textrm{if }a\in\Phi(\bU,\bT,E),
\end{array}\right.
\qquad
f'_\tm(a)=\left\{
\begin{array}{ll}
\tm^+&\textrm{if }a\in\Phi(\bM\cap\bG',\bT,E)\cup\{0\},\\
\tm^+&\textrm{if }a\in\Phi(\overline\bU\cap\bG',\bT,E),\\ 
\tm&\textrm{if }a\in\Phi(\bU\cap\bG',\bT,E),
\end{array}\right.
\]
and
\[
\tilde f_t(a)=\left\{
\begin{array}{ll}
\tm&\textrm{if }a\in\Phi(\bM,\bT,E)\cup\{0\},\\
\tm^+&\textrm{if }a\in\Phi(\overline\bU,\bT,E),\\
\tm&\textrm{if }a\in\Phi(\bU,\bT,E).
\end{array}\right.
\]

\noindent
Here and later, $\Phi(?,\bT,E)$ denotes
the set of $E$-roots of $\bT$ in $?$. To avoid any confusion, 
it should be pointed out that the above $\tilde f_t$ is totally unrelated to 
$\tilde\sigma$ in (\ref{num: HA}).

\medskip

Note that $f_\tm^\ast=\tilde f_{-\tm}$.
We write 
\[
\overline\rU_{\yo,t}=\overline\rU\cap\G_{\yo,t},\quad 
\overline\lieU_{\yo,t}=\overline\lieU\cap\lieG_{\yo,t},\quad
\overline\rU'_{\yo,t}=\overline\rU\cap\G'_{\yo,t},\quad
\overline\lieU'_{\yo,t}=\overline\lieU\cap\lieG'_{\yo,t}.
\]

\smallskip

Recall we have an orthogonal decomposition
$\lieG=\lieG'\oplus\lieG^{\prime\perp}$ with respect to $\bilinear$ when 
$\HypB$ is valid. Otherwise, we can take $\lieG^{\prime\perp}=[\gamma,\lieG]$
to have a linear decomposition of $\lieG$. In both case, 
$\lieG^{\prime\perp}$ is invariant under $\ad(\gamma)$.
Let $\lieG^{\prime\perp}_{\yo,\tilde f_t}
:=\lieG_{\yo,\tilde f_t}\cap\lieG^{\prime\perp}$.
\end{numbering}

\begin{lemma}\label{lem: approx} \ 
Let $Y'\in\lieG'_{\yo,(-\rtm)^+}$, and let $\tm\in\tbR$.
\begin{enumerate}
\item
Let $t\in\bbR$. Then,
$\ad(\gamma+Y'):
\lieG^{\prime\perp}_{\yo,\tilde f_\tm}/
\lieG^{\prime\perp}_{\yo,\tilde f_{\tm+1}}
\longrightarrow
\lieG^{\prime\perp}_{\yo,\tilde f_{\tm-\rtm}}/
\lieG^{\prime\perp}_{\yo,\tilde f_{\tm-\rtm+1}}$ \ 
is an isomorphism.
\item
Let $t\in\bbR$ with $\tm>-\rtm$. Then, 
$^{\G_{\yo,\tilde f_{\rtm+\tm}}}(\gamma+\lieG_{\yo,\tilde f_{\tm}}\cap\lieG')
=\gamma+\lieG_{\yo,\tilde f_{\tm}}$.
\item
If $Y'\in\lieG'_{\yo,(-\rtm)^+}\cap\overline\lieU$, then
$\ad(\gamma+Y'):
(\lieG^{\prime\perp}_{\yo,\tm}\cap\overline\lieU)/
(\lieG^{\prime\perp}_{\yo,\tm+1}\cap\overline\lieU)
\longrightarrow
(\lieG^{\prime\perp}_{\yo,\tm-\rtm}\cap\overline\lieU)/
(\lieG^{\prime\perp}_{\yo,\tm-\rtm+1}\cap\overline\lieU)$ \ 
is an isomorphism.
\item
Let $Y'\in\lieG'_{\yo,(-\rtm)^+}\cap\overline\lieU$.
Suppose $t>-\rtm$. 
Then 
$^{\overline U_{\yo,\rtm+\tm}}
(\gamma+Y'+\overline\lieU'_{\yo,{\tm}})
=\gamma+Y'+\overline\lieU_{\yo,{\tm}}.$
\item
Let $u\in\overline\rU_{\yo,\stm}$.
Then,
$^u(\gamma+\overline\lieU_{\yo,(-\stm)^+})
=(\gamma+\overline\lieU_{\yo,(-\stm)^+})$ if and only if  
$u\in \overline\rU_{\yo,\stm^+}\overline\rU'_{\yo,\stm}
=\overline\rU'_{\yo,\stm}\overline\rU_{\yo,\stm^+}$.
\end{enumerate}
\end{lemma}
\proof
(1) Suppose $Y'=0$.
Since $\tilde f_\tm$ is Galois invariant, we have
$\blieG(E)_{\yo,\tilde f_\tm}\cap\lieG =\lieG_{\yo,\tilde f_\tm}$ 
(\cite{AD2}). Moreover, from the definition of $\tilde f_\tm$
and the properties of $\bilinear$, we have 
$\blieG(E)_{\yo,\tilde f_\tm}\cap \blieG'(E)^{\perp} \cap\lieG
=\lieG_{\yo,\tilde f_\tm}^{\prime\perp}$
where $\blieG(E)=\blieG'(E)\oplus\blieG'(E)^{\perp}$. 
Hence we may assume that $\gamma$ is $k$-split.
Then (1) results from $\gamma$ being good.
If $Y'\ne0$, since $\ad(\gamma)$ is an isomorphism and
$\ad(Y')$ acts pronilpotently, $\ad(\gamma+Y')$ is also
an isomorphism (see also \cite[(2.3.4)]{Asym}).

(2) Having (1), (2) follows from the usual approximation argument 
(see \cite[(6.3)]{AR} or \cite[(8.5)-(8.6)]{yu}). 
We omit the details.

(3) Since $\ad(\gamma+Y')(\overline\lieU)\subset\overline\lieU$,
(3) follows from (1). 

(4) Using (3), one can use an approximation argument.

(5) Since $^u\overline\lieU_{\yo,(-\stm)^+}=\overline\lieU_{\yo,(-\stm)^+}$, 
$^u(\gamma+\overline\lieU_{\yo,(-\stm)^+})
=(\gamma+\overline\lieU_{\yo,(-\stm)^+})$ 
if and only if $^u\gamma-\gamma\in\overline\lieU_{\yo,(-\stm)^+}$.
We claim $^u\gamma-\gamma\in\overline\lieU_{\yo,(-\stm)^+}$
if and only if $u\in \overline U_{\yo,\stm^+}\overline U'_{\yo,\stm}$.
Suppose $u\in \overline\rU_{\yo,\stm^+}\overline\rU'_{\yo,\stm}$.
Then, $^u\gamma\in\,  
^{\overline U_{\yo,\stm^+}}(\gamma+\overline\lieU'_{\yo,(-\stm)^+})
=\gamma+\overline\lieU_{\yo,(-\stm)^+}$ by (4), and `$\Leftarrow$' follows.
Conversely, suppose $^u\gamma\in\gamma+\,{\overline\lieU_{\yo,(-\stm)^+}}$. 
Then by (4), there are $u_1\in\overline U_{\yo,\stm^+}$ 
and $X'\in\overline\lieU'_{\yo,(-\stm)^+}$
such that $^u\gamma=\,^{u_1}(\gamma+X')$. Since $\gamma$ is good and
$X'\in\lieG'_{\yo,(-\rtm)^+}$, $u^{-1}u_1\in(\G'\cap\overline U_{\yo,s})$
by Lemma \ref{lem: appr}-(1).
Hence, $u\in\overline U_{\yo,\stm^+}\overline U'_{\yo,\stm}$.
Now, (5) follows from the claim.
\qed

\begin{numbering}\label{num: more prepare}\rm 
Define
\[
\rK_\vdash^0:=\G'_{\yo,f'_0},\qquad
\rK_\vdash^1:=\G_{\yo,f_\stm},\qquad
\Jv:=\rK_\vdash^0\rK_\vdash^1
\qquad\textrm{and}\qquad 
\tilde\Jv:=\G_{\yo,\tilde f_s}.
\]
Suppose $\phi$ is represented on $\G_{\yo,\stm^+}$ again by $\gamma$.
Then, by \cite[(6.4.44)]{BT1}, we have
$[\rK_\vdash^1,\rK_\vdash^1]\subset\vec\G_{\yo,(\rtm^+,\stm^+)}
\subset\ker(\phi)$. 
Hence, there is a unique character of $\rK_\vdash^1$
with dual blob $\gamma+\lieG_{\yo,f_\stm^\ast}
=\gamma+\lieG_{\yo,\tilde f_{-\stm}}$. 
As before, we denote this character by $\chi_\gamma$.
\end{numbering}

\begin{lemma} \label{lem: d1} \ 
Assume $\HypB$ is valid. Suppose $Y'\in\lieG'_{\yo,(-\rtm)^+}$ and 
$\gamma+Y'+\lieG_{\yo,f_\stm^\ast}$ represent a character of $\rK_\vdash^1$.
Then, the support of $\CaH=\CaH(\G/\!/\rK_\vdash^1,\chi_{\gamma+Y'})$ is 
contained in $\tilde\Jv\G'\tilde\Jv$.
\end{lemma}

\proof
Note that $\tilde\chi_{\gamma+Y'}=\chi_{-\gamma-Y'}$.
Suppose $g\in\G$ is in the support of $\CaH$. Then, 
$^g\chi_{-\gamma-Y'}=\chi_{-\gamma-Y'}$ on 
$^g\!\rK_\vdash^1\cap\rK_\vdash^1$. By Lemma \ref{lem: dual blob}, we have
$^g(-\gamma-Y'+\lieG_{\yo,f_\stm^\ast})\cap(-\gamma-Y'+\lieG_{\yo,f_\stm^\ast})
\neq\emptyset$, and thus 
$^g(\gamma+Y'+\lieG_{\yo,f_\stm^\ast})\cap(\gamma+Y'+\lieG_{\yo,f_\stm^\ast})
\neq\emptyset$. 
By Lemma \ref{lem: approx}-(2), there are $j,\ j'\in \tilde\Jv$ 
and $X'_1,\ X'_2\in\lieG'\cap\lieG_{\yo,f_\stm^\ast}$ such that 
$^{gj'}(\gamma+Y'+X'_1)=\,^{j}(\gamma+Y'+X'_2)$. 
Since $\gamma+Y'+X'_1,\ \gamma+Y'+X'_2\in\gamma+\lieG'_{(-\rtm)^+}$, 
we have $j^{-1}gj'\in\G'$ by Lemma \ref{lem: appr}-(1).
Hence $g\in\tilde\Jv\G'\tilde\Jv$.  
\qed

\begin{remark}\label{rmk: approx}\rm
All the statements in Lemmas \ref{lem: approx} and \ref{lem: d1} 
are also valid when we replace $\gamma$ with $\gamma+Z$ 
for $Z\in \lieZ_{\lieG}$.
\end{remark}

\begin{numbering}\rm
Now, suppose that $\phi'$ is a character of $\G'$ so that
$\phi'|\G'_{\yo,0^+}$ is represented by a $\G$-good element
$\gamma'\in\lieZ_{\lieG'}$ of depth $-\rtm$. However, we do not assume that
$\bC_\bG(\gamma')=\bG'$.
Note that $\gamma'$ still defines a character $\chi_{\gamma'}$ 
on $\Jv$ such that 
$\chi_{\gamma'}|\rU'_{\yo,0}$ and $\chi_{\gamma'}|\overline\rU'_{\yo,0^+}$ 
are trivial.
The following lemma will be used for the proof of 
Proposition \ref{prop: extend}.
\end{numbering}

\begin{lemma}\label{lem: ext d1} 
Suppose $\HypB$ and $\Hypk$ are valid. Let $\vec\stm=(0,\stm)$ and
$\vec\stm^+=(0^+,\stm^+)$. Let
$u\in\overline\rU_{\yo,\stm}$. Let $\gamma'$ and $\chi_{\gamma'}$ be as above. 
Then $^u\chi_{\gamma'}\equiv\chi_{\gamma'}$ on $\vec\G_{\yo,\vec\stm^+}$ and 
$^u\chi_{\gamma'}|\rU'_{\yo,0}=1$.
\end{lemma}

\proof
Note that $\overline\rU_{\yo,\stm}$ normalizes $\Jv$.
Write $\log(u)=(X'+X^\perp)$ for $X'\in\overline\lieU'_{\yo,\stm}$
and $X^\perp\in\overline\lieU_{\yo,\stm}\cap\lieG^{\prime\perp}$.
Then, by $\Hypk$, we have
\[
\Ad(u)(\gamma')=\sum_{n=1}^\infty\frac1{n!}(\ad(X'+X^\perp))^n(\gamma').
\]
Since $[\gamma',\lieG']=0$, \  
$[\gamma', \lieG^{\prime\perp}\cap\overline\lieU]
\subset \lieG^{\prime\perp}\cap\overline\lieU$, \ 
$[\lieG^{\prime\perp}\cap\overline\lieU,\lieG^{\prime\perp}\cap\overline\lieU]
\subset \lieG^{\prime\perp}\cap\overline\lieU$ \ and \ 
$[\lieG'\cap\overline\lieU,\lieG^{\prime\perp}\cap\overline\lieU]
\subset \lieG^{\prime\perp}\cap\overline\lieU$, \ 
from the above formula and $\Hypk$, we conclude 
$^u\gamma'-\gamma'\in(\overline\lieU_{\yo,-\stm}\cap\lieG^{\prime\perp})
\subset(\vec\lieG_{\yo,-\vec\stm}\cap\lieG^{\prime\perp})
=(\vec\lieG_{\yo,-\vec\stm}\cap\lieG^{\prime\perp})$. 
Then $^u\gamma'-\gamma'\in\lieG^{\prime\perp}$ implies that 
$1=\chi_{\gamma'}|\rU'_{\yo,0}=\,^u\chi_{\gamma'}|\rU'_{\yo,0}$, and
$^u\gamma'-\gamma'\in\vec\lieG_{\yo,-\vec\stm}
=(\vec\lieG_{\yo,\vec\stm^+})^\ast$ implies that 
$^u\chi_{\gamma'}|\vec\G_{\yo,\vec\stm^+}
=\chi_{\gamma'}|\vec\G_{\yo,\vec\stm^+}$. 
Hence the lemma follows.
\qed

\begin{center}
{\bf General cases}
\end{center}

\begin{numbering} \label{num: Jv}\rm
Let $\datum_\yo=(\vec\bG,\yo,\vec\rtm,\vec\phi)$ be a good $\G$-datum 
of positive depth.
Let $\bT\subset\bG^0$ be an $E$-split maximal $k$-torus such that 
$\yo\in\Apt(\bG,\bT,k)$. Let $\bS$ be a $k$-split subtorus in $\bT$.
Let $\bM$ be the $k$-Levi subgroup of $\bG$ defined by $\bC_\bG(\bS)$.
Then, $\yo\in\Bd(\bM,k)$. 
Fix a $k$-parabolic subgroup $\bP$ with Levi decomposition $\bP=\bM\bU$, 
and let $\overline\bU$ be the unipotent subgroup opposite to $\bU$.
Let 
\[
\bM^i:=\bM\cap\bG^i=\bC_{\bG^i}(\bS),\qquad\bU^i:=\bU\cap\bG^i,\qquad
\overline\bU^i=\overline\bU\cap\bG^i
\]
for $i=0,1,\cdots,d$. Note that each $\bM^i$ is a $k$-Levi subgroup of $\bG^i$.

Define $f^i:\Phi(\bG^i,\bT,E)\cup\{0\}\longrightarrow\tbR$ as follows: 
if $i\ge1$, 
\[
f^i(a)=\left\{
\begin{array}{ll}
\stm_{i-1}^+&\textrm{if }a\in\Phi(\bM^i,\bT,E)\cup\{0\},\\
\stm_{i-1}^+&\textrm{if }a\in\Phi(\overline\bU^i,\bT,E),\\
\stm_{i-1}&\textrm{if }a\in\Phi(\bU^i,\bT,E),
\end{array}\right.
\quad\textrm{and}\quad
f^0(a)=\left\{
\begin{array}{ll}
0^+&\textrm{if }a\in\Phi(\bM^0,\bT,E)\cup\{0\},\\
0^+&\textrm{if }a\in\Phi(\overline\bU^0,\bT,E),\\
0&\textrm{if }a\in\Phi(\bU^0,\bT,E).
\end{array}\right.
\]
For $i=1,\cdots,d$, define 
$\tilde f^i:\Phi(\bG^i,\bT,E)\cup\{0\}\longrightarrow\tbR$ as 
\[
\tilde f^i(a)=\left\{
\begin{array}{ll}
\stm_{i-1}&\textrm{if }a\in\Phi(\bM^i,\bT,E)\cup\{0\},\\
\stm_{i-1}^+&\textrm{if }a\in\Phi(\overline\bU^i,\bT,E),\\
\stm_{i-1}&\textrm{if }a\in\Phi(\bU^i,\bT,E).
\end{array}\right.
\]
We also define corresponding open compact subgroups of $\G^i$
as follows: for $i=0,1,\cdots,d$, 
\[
\Jvy^i:=\G^i_{\yo,f^i}\subset\G^i,
\]
and 
\[
\tilde\rK_{\yo\vdash}^0:=1,\qquad
{\tilde K}_{\yo\vdash}^i:=\G^i_{\yo,\tilde f^i}\subset\G^i.
\]
If there is no confusion, we omit $\yo$ from the above notation and
simply write $\Jv^i$ and ${\tilde K}_{\vdash}^i$.
Note that we have $\Jv^0\subset\G^0$, and 
$\Jv^0=\G^0_{\yo,0^+}$ if $\bG^0\subset\bM$.
Define
\[
\Jvy=\Jv:=\Jv^0\Jv^1\cdots \Jv^d,
\qquad\textrm{and}\qquad
{\tilde K}_{\yo\vdash}=\tilde\Jv
:={\tilde K}_\vdash^0{\tilde K}_{\vdash}^1\cdots{\tilde K}_{\vdash}^d.
\]
Note that the above open compact subgroups depend also on the choice of 
$\bM$ and $\bU$. 

If $\datum_y$ is strongly good, each $\phi_i$ defines a character of $\Jv$ 
represented by $\gss_i$ and $\vec\phi$ defines a character of $\Jv$ 
represented by $\sgss$. We again denote these characters by 
$\chi_{\gss_i}$ and $\chi_\sgss$ respectively.

For $i=0,\cdots,d$ and $t\in\tbR_{\ge0}$, we have 
$\rM^i_{\yo,\tm}=\G^i_{\yo,\tm}\cap \rM$. Write
\[
\overline\rU^i_{\yo,\tm}:=\G^i_{\yo,\tm}\cap\overline\rU,\qquad\qquad
\rU^i_{\yo,\tm}:=\G^i_{\yo,\tm}\cap\rU.
\]
For any admissible sequence $\vec\tm$ of length $\ell(\datum)$, we also write
\[
\overline\rU_{\yo,\vec\tm}:=\vec\G_{\yo,\vec\tm}\cap\overline\rU, \qquad
\rM_{\yo,\vec\tm}:=\vec\G_{\yo,\vec\tm}\cap \rM, \qquad
\rU_{\yo,\vec\tm}:=\vec\G_{\yo,\vec\tm}\cap\rU.
\]
Then note that $\rU_{\yo,\vec\stm}=\Jv\cap\rU=\tilde\Jv\cap\rU$, 
$\overline\rU_{\yo,\vec\stm^+}=\Jv\cap\overline\rU=\tilde\Jv\cap\overline\rU$,
$\rM_{\yo,\vec\stm^+}=\Jv\cap\rM$, etc.
\end{numbering}

\begin{lemma}\label{lem: more approx}\ 
For (2) and (3), we assume $\Hypk$ and $\HypB$ are valid.
\begin{enumerate}
\item
For $i=0,\cdots,d-1$, the following is an isomorphism:
\[
\ad(\gss_i):(\overline\lieU_{\yo,\vec\stm}\cap\lieG^{i\perp})/
(\overline\lieU_{\yo,\vec\stm^+}\cap\lieG^{i\perp})
\longrightarrow
(\overline\lieU_{\yo,-\vec\stm}\cap\lieG^{i\perp})/
(\overline\lieU_{\yo,(-\vec\stm)^+}\cap\lieG^{i\perp}).
\]
\item
$|\overline\lieU_{\yo,\vec\stm}/\overline\lieU_{\yo,\vec\stm^+}|
=|\overline\lieU_{\yo,-\vec\stm}/\overline\lieU_{\yo,(-\vec\stm)^+}|
=\left|\rU_{\yo,\vec\stm}/\rU_{\yo,\vec\stm^+}\right|
=\left|\Jv/\rK^+_\datum\right|$.
\item 
$|\overline\lieU_{\yo,\vec\stm(0^+)}/\overline\lieU_{\yo,\vec\stm^+}|
=\left|\rU_{\yo,\vec\stm(0^+)}/\rU_{\yo,\vec\stm^+}\right|
=\left|\Jv/(\rU^0_{\yo,0}\rK^+_\datum)\right|$ 
where $\vec\stm(0^+)=(0^+,\stm_0,\cdots,\stm_{d-1})$.
\end{enumerate}
\end{lemma}

\proof
(1) Since $\gss_i$ is a $\G^{i+1}$-good element of depth $-\rtm_i$ and
$\ad(\gss_i)(\overline\lieU)\subset(\overline\lieU)$, 
the above $\ad(\gss_i)$ is an isomorphism
by \cite[(2.3.1)]{Ad}. 

(2) Since $\lieG=\lieG^0\oplus\lieG^{0\perp}\oplus\cdots\lieG^{d-1\perp}$, 
we have
\[
\begin{array}{ll}
|\overline\lieU_{\yo,\vec\stm}/\overline\lieU_{\yo,\vec\stm^+}|
&=|\overline\lieU^0_{\yo,0}/\overline\lieU^0_{\yo,0^+}| 
\cdot\prod_{i=0}^{d-1}
|(\overline\lieU_{\yo,\vec\stm}\cap\lieG^{i\perp})/
(\overline\lieU_{\yo,\vec\stm^+}\cap\lieG^{i\perp})| \\
&=|\overline\lieU^0_{\yo,0}/\overline\lieU^0_{\yo,0^+}| 
\cdot\prod_{i=0}^{d-1}
|(\overline\lieU_{\yo,-\vec\stm}\cap\lieG^{i\perp})/
(\overline\lieU_{\yo,(-\vec\stm)^+}\cap\lieG^{i\perp})|
=|\overline\lieU_{\yo,-\vec\stm}/\overline\lieU_{\yo,(-\vec\stm)^+}|.
\end{array}
\]
The second equality follows from (1).
By duality, 
\[
\left|\overline\lieU_{\yo,-\vec\stm}/\overline\lieU_{\yo,(-\vec\stm)^+}\right|
=\left|(\lieU^\ast_{\yo,\vec\stm^+}\cap\overline\lieU)
/(\lieU^\ast_{\yo,\vec\stm}\cap\overline\lieU)\right|
=\left|\lieU_{\yo,\vec\stm}/\lieU_{\yo,\vec\stm^+}\right|.
\]
Since $\log(\Jv)\cap\lieU=\lieU_{\yo,\vec\stm}$, 
$\log(\rK^+_\datum)\cap\lieU=\lieU_{\yo,\vec\stm^+}$, and
$\log(\Jv)\cap(\overline\lieU\oplus\lieM)
=\log(\rK^+_\datum)\cap(\overline\lieU\oplus\lieM)$,
we have 
$\left|\rU_{\yo,\vec\stm}/\rU_{\yo,\vec\stm^+}\right|
=\left|\lieU_{\yo,\vec\stm}/\lieU_{\yo,\vec\stm^+}\right|
=\left|\log(\Jv)/\log(\rK^+_\datum)\right|
=\left|\Jv/\rK^+_\datum\right|$.

(3) By duality, we have
$|\overline\lieU^0_{\yo,0}/\overline\lieU^0_{\yo,0^+}|
=|\lieU^0_{\yo,0}/\lieU^0_{\yo,0^+}|$. Denote this number by $m$.
Dividing each term in (2) by $m$, we get the equalities in (3).
\qed

\begin{proposition}\label{prop: extend}
Suppose $\datum=\datum_\yo$ is a strongly good $\G$-datum of positive depth.
Suppose $\HypB$ and $\Hypk$ are valid. 
Let $\chi$ be an irreducible representation of $\Jv$ such that
$\chi|\rK^+_\datum$ contains $\chi_\gss$ and 
$\chi|\rU^0_{\yo,0}\equiv 1$.
Then there is a $u\in\overline U_{\yo,\vec\stm(0^+)}$ 
such that $\chi\equiv\, ^u\chi_{\sgss}$ on $\Jv$.
\end{proposition}

\proof
Let $u\in\overline U_{\yo,\vec\stm(0^+)}$. Then, we can write 
$u=u_d\cdots u_1u_0$ for some  $u_0\in\overline U^0_{\yo,0^+}$ and  
$u_i\in\overline U^i_{\yo,\stm_{i-1}}$, $i=1,\cdots,d$.

For the character $\chi_{\gss_j}$ of $\Jv$, if  $j\ge i$, 
we have $^{u_i}\chi_{\gss_j}\equiv\chi_{\gss_j}$ on $\Jv$
since $u_i\in\G^j$.
Then $^{u_du_{d-1}\cdots u_0}\chi_{\gss_{i-1}}
=\,^{u_du_{d-2}\cdots u_i}\chi_{\gss_{i-1}}
=\,^{u_iu_i^{-1}(u_du_{d-2}\cdots u_{i+1})u_i}\chi_{\gss_{i-1}}
=\,^{u_i}\chi_{\gss_{i-1}}$.
The last equality follows from 
$u_i^{-1}(u_du_{d-2}\cdots u_{i+1})u_i\in\overline\rU_{\yo,\stm_i}
\subset\ker(\chi_{\gss_{i-1}}|\G_{\yo,\stm_{i-1}^+})$.
Hence $^{u_du_{d-1}\cdots u_0}\chi_\gss
=\,^{u_1}\chi_{\gss_0}{}^{u_2}\chi_{\gss_1}
\cdots{}^{u_d}\chi_{\gss_{d-1}}\chi_{\gss_d}$ on $\Jv$.
Since $\rK^+_\datum\subset (\bG^{i-1},\bG)(k)_{\yo,(0^+,\stm_{i-1}^+)}$
and $\rU^0_{\yo,0}\subset\rU^{i-1}_{\yo,0}$, 
from Lemma \ref{lem: ext d1}, we deduce that
$^{u_i}\chi_{\gss_{i-1}}\equiv\chi_{\gss_{i-1}}$
on $\rK^+_\datum$ and $^{u_i}\chi_{\gss_{i-1}}| U^0_{\yo,0}\equiv1$.
Hence, for any $u\in\overline U_{\yo,\vec\stm(0^+)}$, 
we have $^u\chi_\gss|\rK^+_\datum=\chi_\gss|\rK^+_\datum$ 
and $^u\chi_\gss| U^0_{\yo,0}=1$.

We claim that $^u\chi_{\gss}\equiv\chi_{\gss}$ on $\Jv$ 
if and only if $u\in\overline U_{\yo,\vec\stm^+}$. 
Since $u\in\overline U_{\yo,\vec\stm^+}\subset\Jv$, `$\Leftarrow$' is obvious.
To prove `$\Rightarrow$', we first observe that 
$^{u}\chi_\gss
=\chi_{\gss_0}\chi_{\gss_1}\cdots\,^{u_i}\!\chi_{\gss_{i-1}}
\cdots\chi_{\gss_d}$
on $\G^{i-1}_{\yo,\stm_{i-2}^+}\Jv^i$.
Hence, $^u\chi_{\gss}\equiv\chi_{\gss}$ on $\Jv$ implies
$^{u_i}\chi_{\gss_{i-1}}=\chi_{\gss_{i-1}}$ 
on $\G^{i-1}_{\yo,\stm_{i-2}^+}\Jv^i$. Then,
$u_i\in\overline U^i_{\yo,\stm_{i-1}^+}\overline U^{i-1}_{\yo,\stm_{i-1}}$
by Lemma \ref{lem: approx}-(5).
Hence, $u=u_d\cdots u_0\in\overline U^d_{\yo,\stm_{d-1}^+}\cdots
\overline U^1_{\yo,\stm_{0}^+}\overline U^{0}_{\yo,0^+}$. 
The claim is now proved.

From the claim and Lemma \ref{lem: more approx}, 
$|\{\,^u\chi_\gss\mid u\in\overline U_{\yo,\vec\stm(0^+)}\}|
=|\overline U_{\yo,\vec\stm(0^+)}/\overline U_{\yo,\vec\stm^+}|
=\left|\Jv/(\rU^0_{\yo,0}\rK^+_\datum)\right|$.
Hence, by counting, any irreducible extension of $\chi_\gss$ of $\rK^+_\datum$ 
to $\Jv$ which is trivial on $U^0_{\yo,0}$ 
is of the form $^u\chi_\gss$ for some $u\in\overline U_{\yo,\vec\stm(0^+)}$.
\qed

\medskip

For later use, we record the following corollary:

\begin{corollary}\label{cor: extend}
Suppose $\datum=\datum_\yo$ is a strongly good $\G$-datum of positive depth. 
Assume $\HypB$ and $\Hypk$ are valid. 
Let $(\pi,V_\pi)\in\admG$. Suppose $(\pi,V_\pi)$ 
contains $(\rK^+_\datum,\chi_\gss)$. Then
$(\pi,V_\pi)$ also contains $(\Jv,\chi_\gss)$ if
\begin{enumerate}
\item
$\rU^0_{\yo,0}=\rU^0_{\yo,0^+}$, or 
\item
there is a nonzero vector 
in the $(\rK^+_\datum,\chi_\gss)$-isotypic component of $V_\pi$ 
which is invariant under $\rU^0_{\yo,0}$.
\end{enumerate}
\end{corollary}

The following lemma is a preparation for Proposition \ref{prop: support}
and Proposition \ref{prop: cusp I}. 
For any open compact subgroup $J\subset\G$, we write 
$J_\ell:=J\cap\overline U$, $J_M:=J\cap M$ and $J_u:=J\cap U$.

\begin{lemma}\label{lem: Jdec}\ 
\begin{enumerate}
\item
$\tilde\rK_{\vdash}^d
=\tilde\rK_{\vdash\ell}^d\tilde\rK_{\vdash\rM}^d\tilde\rK_{\vdash u}^d
=\tilde\rK_{\vdash\rM}^d\rK_{\vdash\ell}^d\rK_{\vdash u}^d
=\rK_{\vdash\ell}^d\rK_{\vdash u}^d\tilde\rK_{\vdash\rM}^d$.
\item
$\tilde\Jv=\tilde\rK_{\vdash\ell}\tilde\rK_{\vdash\rM}\tilde\rK_{\vdash u}$.
\item
Suppose $\bG^0\subset\bM$. Then,
$\Jv=\rK_{\vdash\ell}\rK_{\vdash\rM}\rK_{\vdash u}$.
\end{enumerate}
\end{lemma}

\proof
Proving (3), we first note that each $f^i$, $i=1,\cdots,d$ is concave 
and positive.
From \cite[(6.4.9)]{BT1} and the concavity of $f^i$, we deduce that 
$\bG^i(E)_{\yo,f^i}$ is decomposable with respect to $\bP(E)=\bM(E)\bU(E)$.
That is, $\bG^i(E)_{\yo,f^i}
=(\bG^i(E)_{\yo,f^i}\cap\overline\bU(E))\cdot
(\bG^i(E)_{\yo,f^i}\cap\bM(E))\cdot
(\bG^i(E)_{\yo,f^i}\cap\bU(E)).$
From the positivity of $f^i$, we have 
$\Jv^i=\G^i_{\yo,f}=\bG^i(E)_{\yo,f^i}^{\Gal(E/k)}$.
Since $\bM$, $\bU$, and $\overline\bU$ are $\Gal(E/k)$-stable, 
$\Jv^i=\G^i_{\yo,f}$ is also decomposable with respect to $P=MU$, that is,
$\Jv^i=\rK^i_{\vdash\ell}\rK^i_{\vdash\rM}\rK^i_{\vdash u}$. 
If $i=0$, since $\bG^0\subset\bM$, $\Jv^0=\rK^0_{\vdash M}$.

Now, let $\Jv(i)=\Jv^i\cdots\Jv^d$. We claim that each $\Jv(i)$ is
decomposable with respect to $P=MU$. If $i=d$, it follows from the previous
paragraph. Suppose $\Jv(i+1)$ is decomposable. 
Note that $\Jv(i+1)$ is normalized by $\Jv^i$
and $\Jv(i)=\Jv^i\Jv(i+1)$. Let $g\in\Jv(i)$. 
Then $g=ab$ for some $a\in\Jv^i$ and $b\in\Jv(i+1)$. 
Write $a=a_\ell a_\rM a_u$ compatible with the decomposition of $\Jv^i$. 
Since $b':=a_u b a_u^{-1}\in\Jv(i+1)$, we can also write 
$b'=b'_\ell b'_\rM b'_u$. Then,
$ab=a_\ell a_\rM a_u b a_u^{-1}a_u
=(a_\ell a_\rM b'_\ell a_\rM^{-1})(a_\rM b'_\rM)(b'_u a_u)$, where
$a_\ell a_\rM b'_\ell a_\rM^{-1}\in\Jv(i)_\ell$, 
$a_\rM b'_\rM\in\Jv(i)_\rM$ and $b'_u a_u\in\Jv(i)_u$.
Hence, $\Jv(i)$ is decomposable. 
Since $\Jv=\Jv(0)$, $\Jv$ is also decomposable.

The second and the third equalities in (1) follow from 
$\tilde\rK_{\vdash\ell}^d=\rK_{\vdash\ell}^d$ and
$\tilde\rK_{\vdash u}^d=\rK_{\vdash u}^d$.
The other equalities are proved similarly as in (3).
\qed

\begin{proposition}\label{prop: support}
Let $\datum=\datum_\yo$ be a strongly good $\G$-datum of positive depth. 
Suppose $\HypB$ is valid. 
The support of $\CaH(\G/\!/\Jv,\chi_\gss)$ is contained in 
$\tilde\Jv\G^0\tilde\Jv$. 
\end{proposition} 

The following is a modification of the proof of \cite[(4.1)]{yu}.

\proof
We prove the proposition by induction. If $d=0$, it is obvious.
Now, assume $d\ge1$. Let $\datum'=(\vec\bG',\xo,\vec\rtm',\vec\phi')$ be
such that $\vec\bG'=(\bG^0,\cdots,\bG^{d-1})$, 
$\vec\phi'=(\phi_0,\cdots,\phi_{d-1})$, 
and $\vec\rtm'=(\rtm_0,\cdots,\rtm_{d-1})$. 
Then note that $\ell(\datum')=d-1$. Let 
\[
\Jv':=\Jv^0\Jv^1\cdots \Jv^{d-1},\qquad
\textrm{and}\qquad
{\tilde K}_\vdash'
:={\tilde K}_\vdash^0{\tilde K}_\vdash^1\cdots{\tilde K}_\vdash^{d-1}.
\]
Suppose $g\in\G$ is in the support of $\CaH(\G/\!/\Jv,\chi_\gss)$.
Since $\phi_d|\Jv$ is the restriction of $\phi_d$ which is defined on 
the whole of $\G$, $g$ also intertwines 
$\theta'=(\chi_\gss\phi_d^{-1})|\Jv=\prod_{j=0}^{d-1}\chi_{\gss_j}|\Jv$.
Therefore, $g$ also intertwines $\theta'|\Jv^d$. 
Note that if $\rtm_j\le\stm$, $\chi_{\gss_j}$ is trivial on $\Jv^d$, and 
$\theta'|\Jv^d=\chi_{\gss_{d-1}+Y}$ for $Y=\sum_{\rtm_j>\stm_{d-1}}\gss_j$.
By Lemma \ref{lem: d1} and Lemma \ref{lem: Jdec}, 
there are $j_1, j_2\in{\tilde K}_{\vdash M}^d=\rM^d_{\yo,\stm_{d-1}}$,
$h_1, h_2\in\rK_{\vdash\ell}^d\rK_{\vdash u}^d$ and
$g'\in\G^{d-1}$ such that $g=h_1j_1g'j_2h_2$. 
Since $h_1,h_2\in\Jv$, $j_1g'j_2$ intertwines $\theta'$. 
Now, from Lemma \ref{lem: sup}-(1), we may further assume 
$j_1,j_2\in\rJ_M^d$, where $\rJ_M^d$ is as in that lemma. Then, 
by Lemma \ref{lem: sup}-(2), $g'$ also intertwines 
$\theta'|\Jv'=\chi_{\gss'}$ where 
$\gss'=\gss_{d-1}+\gss_{d-1}+\cdots+\gss_0$.
By the induction hypothesis, 
$g'\in{\tilde\rK}_\vdash'\G^0{\tilde\rK}_\vdash'$. Hence 
$g\in {\tilde\rK}_\vdash^d{\tilde\rK}_\vdash'\G^0{\tilde\rK}_\vdash'
{\tilde\rK}_\vdash^d=\tilde\Jv\G^0\tilde\Jv$.
\qed

\begin{lemma}\label{lem: sup}
We keep the notation from the proof of Proposition \ref{prop: support}. Let 
\[
\rJ_M^d=(\bM^{d-1},\bM^d)(k)_{\yo,(\rtm_{d-1},\stm_{d-1})}.
\]
\begin{enumerate}
\item 
$\rM^d_{\yo,\stm_{d-1}}\G^{d-1}\rM^d_{\yo,\stm_{d-1}}
\subset\rJ_M^d\G^{d-1}\rJ^d_M$.
\item 
If $g\in\G^{d-1}$ intertwines $\theta'$, $j\in\rJ_M^d$, 
then $gj$ and $jg$ also 
intertwines $\theta'$.
\end{enumerate}
\end{lemma}

\proof
(1) follows from 
$\rM^d_{\yo,\stm_{d-1}}=\rM^{d-1}_{x,\stm_{d-1}}\rJ_M^d
=\rJ_M^d\rM^{d-1}_{x,\stm_{d-1}}\subset\rJ_M^d\G^{d-1}
\subset\rJ_M^d\G^{d-1}\rJ_M^d$.
(2) can be proved as in \cite[(4.3)]{yu}.
\qed

\section{\bf Supercuspidality I}

Let $\datum=\datum_\yo$ be a strongly good $\G$-datum of positive depth, 
and let $\bM$, $\bP=\bM\bU$ and 
$\Jv$ be as in (\ref{num: Jv}).
In this section, we prove that if there is 
a supercuspidal representation $(\pi,V_\pi)$ such that 
$(\rK^+_\datum,\chi_\gss)<(\pi,V_\pi)$,
$\bZ^0/\bZ^d$ is $k$-anisotropic  (see Proposition \ref{prop: cusp I}).

We first recall the following from \cite{BK2}. 
For any open compact subgroup $J\subset\G$, $J_\ell$, $J_M$ and
$J_u$ are defined as before (see Lemma \ref{lem: Jdec}).

\begin{definition}\label{def: decomposable} \cite[(6.16)]{BK2}\rm \ 
Let $P=MU$ be as above. Suppose $J$ is an open compact subgroup
decomposable with respect to $P=MU$, that is,
$J=J_\ell J_M J_u$. Then $\zeta\in Z_M$ is {\it strongly positive 
relative to $(P,J)$} if 
\begin{enumerate}
\item
$\zeta J_u\zeta^{-1}\subset J_u,\ \ \zeta^{-1}J_{\ell}\zeta\subset J_{\ell}$;
\item
for any compact open subgroups $H_1,\ H_2$ of $U$, there exists an integer
$m\ge0$ such that $\zeta^mH_1\zeta^{-m}\subset H_2$;
\item
for any compact open subgroups $K_1,\ K_2$ of $\overline U$, 
there exists an integer $m\ge0$ such that $\zeta^{-m}K_1\zeta^{m}\subset K_2$.
\end{enumerate}
\end{definition}

Combining \cite[(6.14)]{BK2} and the remark after \cite[(6.16)]{BK2},
a strongly positive element always exists 
relative to $(P,J)$ such that $J=J_\ell J_M J_u$.

The following can be regarded as a corollary of the proof of 
\cite[(6.10)]{BK2}.

\begin{lemma}\label{lem: BK2}
Let $J$, $\overline U$ and $P=MU$ be as in Definition \ref{def: decomposable}.
Then, $JMJ\cap \overline U = J_\ell$. 
\end{lemma}

\proof
$\supset$ is obvious. To prove $\subset$,
suppose $u\in JMJ\cap\overline U$. Let $t\in\rM$ be so that
$u\in JtJ$. Then, there are $j,h\in J$ so that $uj=ht$.
Write $j=j_\ell j_Mj_u$ and $h=h_\ell h_Mh_u$ compatible with
the decomposition $J=J_\ell J_M J_u$.
Then, $(uj_\ell) j_M j_u=h_\ell(h_M t)(t^{-1}h_ut)$.
Hence, by the uniqueness of Iwahori decomposition,
we have $u=h_\ell j_\ell^{-1}\in J_\ell$, and $\subset$ follows.
\qed

\begin{proposition}\label{prop: strongly positive}
Assume $\HypB$ is valid. Suppose $\bG^0\subset\bM\subsetneq\bG$. 
Then there is a strongly positive element $\zeta$ 
relative to $(P,\Jv)$ such that an element 
$f_\zeta\in\CaH(\G/\!/\Jv,\chi_\gss)$ supported on $\Jv\zeta\Jv$ is invertible.
\end{proposition}

\proof
We fix a Haar measure on $\G$ so that $\vol(\Jv)=1$. 
From $\bG^0\subset\bM$ and Proposition \ref{prop: support},  the support of 
$\CaH:=\CaH(\G/\!/\Jv,\chi_\gss)$ is contained in $\Jv M\Jv$. 
Let $\zeta\in Z_M\subset\G^0$ be a strongly positive element such that 
$\overline U=\cup_{n\ge0}\zeta^n K_{\vdash\ell} \zeta^{-n}$. 
For $j=j_\ell j_M j_u\in\Jv\cap(\zeta\Jv\zeta^{-1})$, since 
$^\zeta\chi_\gss(j)=\chi_\gss(\zeta^{-1} j \zeta)
=\chi_\gss(\zeta^{-1}j_\ell\zeta)\,\chi_\gss(\zeta^{-1}j_M\zeta)\,
\chi_\gss(\zeta^{-1}j_u\zeta)
=\chi_\gss(j_M)=\chi_\gss(j)$, we have $\zeta\in\supp(\CaH)$. 
Here, the third equality follows from $\rK_{\vdash\ell}$, 
$\rK_{\vdash u}\subset\ker(\chi_\gss)$ and $\zeta^{-1}j_M\zeta=j_M$.
Similarly, $\zeta^{-1}\in\supp(\CaH)$.
Let $f_\zeta$ (resp. $f_{\zeta^{-1}}$) be an element of $\CaH$ which is 
the function supported on $\Jv\zeta \Jv$ (resp. $\Jv\zeta^{-1}\Jv$) with
$f_\zeta(\zeta)=1$ (resp. $f_{\zeta^{-1}}(\zeta^{-1})=1$).
We claim 
\[
f_\zeta\conv f_{\zeta^{-1}}=f_{\zeta^{-1}}\conv f_\zeta=
\vol(\Jv\zeta\Jv)f_1
\]
where $f_1$ is the identity in $\CaH$ given by
the function supported on $\Jv$ with $f_1(1)=1$.
Note that $\supp(f_\zeta\conv f_{\zeta^{-1}})
\subset (\Jv\zeta\Jv\zeta^{-1}\Jv)\cap(\tilde\Jv\G^0\tilde\Jv)$. 
From Lemma \ref{lem: support}, 
$\supp(f_\zeta\conv f_{\zeta^{-1}})\subset \Jv$.
Hence $f_\zeta\conv f_{\zeta^{-1}}=\,{c}\cdot f_1$ for some constant $c$. 
Evaluating both sides at the identity, $c=\vol(\Jv\zeta\Jv)$. 
Similar computation shows that 
$f_{\zeta^{-1}}\conv f_\zeta=\vol(\Jv\zeta\Jv)f_1$.
Hence, the inverse of $f_\zeta$ is $\frac1{\vol(\Jv\zeta\Jv)}f_{\zeta^{-1}}$.
\qed

\begin{lemma}\label{lem: support} 
We keep the notation and the situation from the proof of
Proposition \ref{prop: strongly positive}.
\[
(\Jv\zeta \Jv \zeta^{-1}\Jv)\cap (\tilde\Jv\G^0\tilde\Jv)= \Jv
=(\Jv\zeta^{-1} \Jv \zeta\Jv)\cap (\tilde\Jv\G^0\tilde\Jv).
\]
\end{lemma}

\proof
It is obvious that 
$\Jv\subset(\Jv\zeta\Jv\zeta^{-1}\Jv)\cap(\tilde\Jv\G^0\tilde\Jv)$.
Suppose $g\in (\Jv\zeta \Jv \zeta^{-1}\Jv)\cap(\tilde\Jv\G^0\tilde\Jv)$.
Since $\zeta\rK_{\vdash\rM}\zeta^{-1}=\rK_{\vdash\rM}$
and $\zeta\rK_{\vdash u}\zeta^{-1}\subset\rK_{\vdash u}$,
We can write $g=j_1\zeta j_{\ell}\zeta^{-1} j_2$ 
for some $j_1,j_2\in\Jv$ and $j_{\ell}\in\rK_{\vdash\ell}$. 
Since $\zeta j_{\ell}\zeta^{-1}\in\tilde\Jv\G^0\tilde\Jv\cap\overline U$,
by Lemma \ref{lem: BK2}, 
$\zeta j_{\ell}\zeta^{-1}\in\tilde\rK_{\vdash\ell}=\rK_{\vdash\ell}$.
Hence, $g\in\Jv$. The other equality can be proved in a similar way.
\qed

\begin{proposition}\label{prop: cusp I}
Assume $\HypB$ and $\Hypk$ are valid.
Let $(\pi,V_\pi)$ be a supercuspidal representation.
Let $\datum=(\vec\bG,\yo,\vec\rtm,\vec\phi)$ be a strongly
good $\G$-datum of positive depth.
Suppose that $(\rK^+_\datum,\phi_\datum)$ is contained in $(\pi,V_\pi)$.
Then, $\bZ^0/\bZ^d$ is $k$-anisotropic.
\end{proposition}

\proof
Suppose $\bZ^0/\bZ^d$ is not anisotropic. Then, there is
a $k$-split subtorus $\bS$ of $\bZ^0$ such that 
$\bM=\bC_\bG(\bS)$ is a proper $k$-Levi subgroup of $\bG$. 
Note that $\bG^0\subset\bM$. Hence $\yo\in\Bd(\bM,k)$.
Fix a $k$-parabolic subgroup $\bP$ with Levi decomposition $\bP=\bM\bU$, 
and let $\overline\bU$ be the unipotent subgroup opposite to $\bU$.
Since $\bG^0\subset\bM$, $\rU^0_{\yo,0}=1$ and 
$(\pi,V_\pi)$ contains $(\Jv,\chi_\gss)$ by Corollary \ref{cor: extend}.
Now, consider $\CaH:=\CaH(\G/\!/\Jv,\chi_\sgss)$. 
By Proposition \ref{prop: strongly positive}, there is a strongly positive
element $\zeta$ relative to $(P,\Jv)$ such that
$f_\zeta $, an element of $\CaH(\G/\!/\Jv,\chi_\sgss)$ 
supported on $\Jv\zeta \Jv$, is invertible. 
This implies that the Jacquet map restricted to 
$V_\pi^{\chi_\sgss}$, the $(\Jv,\chi_\gss)$-isotypic component in $V_\pi$,
is injectively mapped into $V_{\pi U}^{\chi_\sgss}$, the 
$(\Jv\cap M,\chi_\gss)$-isotypic component in the Jacquet module $V_{\pi U}$ 
(\cite{BK2}).
Since $V_\pi^{\chi_\sgss}\ne0$, this implies that 
the Jacquet module $V_{\pi U}$ of $V_\pi$ with respect to $U$ is nontrivial,
which contradicts to the supercuspidality of $(\pi,V_\pi)$.
Hence, $\bZ^0/\bZ^d$ is $k$-anisotropic.
\qed

\section{\bf Depth zero types}
\label{sec: depth zero}

In this section and the following one, we deduce some facts that we need
for \S\ref{sec: cusp II} (Proposition \ref{prop: cusp II}).
Most results here are from \cite{MP2}.

\begin{numbering}\label{num: setup}\rm
Let $(\bG',\bG)$ be a tamely ramified twisted Levi sequence, which 
splits over a tamely ramfied extension $E$. For any $y\in\Bd(\bG',k)$,
let $\sfM'_y$ be the quotient of the reduction mod $\pid_k$ of the 
$\CaO_k$-group scheme $\sfP'_y$
associated to $\G'_{y,0}$ by its unipotent radical. 
Then, $\G'_{y,0}/\G'_{y,0^+}\simeq\sfM'_y(\bbF_q)$.

Let $\bT\subset\bG'$ be an $E$-split maximal $k$-torus
such that $\yo\in\Apt(\bG',\bT,k)$ and $\Apt(\bG',\bT,k)$ is 
of maximal dimension. Let $\bS$ be a maximal $k$-split subtorus in $\bT$.
We attach a $k$-Levi subgroup $\bM$ of 
$\bG$ to $\G'_{\yo,0}$ (see also \cite[(6.3)]{MP2}):
first note that $\bS$ gives rise to a maximal $\bbF_q$-torus 
$\sfS$ of $\sfM'_{\yo}$. 
Let $\sfC$ be the
maximal $\bbF_q$-split torus contained in the center of of $\sfM'_{\yo}$. 
Lift $\sfC$ to $\bS$ to get a subtorus $\bC$ of $\bS$.  
Let $\bM$ be the centralizer of $\bC$ in $\bG$.
Note that since $\bT\subset\bM$, $\yo\in\Bd(\bM,k)$. Moreover, if
$\G'_{y,0}$ is not a maximal parahoric subgroup of $\G'$, $\rM$ is a proper 
$k$-Levi subgroup.

For a given $\G'_{y,0}$, if $\bM$ is chosen in the above manner, we 
will say that $\bM$ is \emph{adapted} for $\G'_{\yo,0}$.
In \cite[(6.3)]{MP2}, the case $\bG'=\bG$ is considered.
\end{numbering}

\begin{remark} \rm
We keep the situation from (\ref{num: setup}). 
Let $\bM':=\bM\cap\bG'=\bC_{\bG'}(\bC)$. Then, $\bM'$ is 
a $k$-Levi subgroup of $\bG'$ associated to $\G'_{\yo,0}$ as in 
\cite[(6.3)]{MP2}. 
Set $\rM'_{\yo,\rtm}:=\rM\cap\G'_{\yo,\rtm}$ for $\rtm\in\tbR$ as before. 
Then, by \cite[(6.4)]{MP2}, we have the following:
\begin{enumerate}
\item
$\rM'_{\yo,0}$ is a maximal parahoric subgroup of $\rM'$.
\item
$\G'_{\yo,0}/\G'_{\yo,0^+}\simeq \rM'_{\yo,0}/\rM'_{\yo,0^+}\simeq
\sfM'_{\yo}(\bbF_q)$.
\end{enumerate}
\end{remark}

\begin{numbering}\label{num: fact2}\rm
Let $(\pi,V_\pi)$ be a smooth representation of $\G$.
Let $x,y,y'\in\Bd(\bG,k)$ be such that $\G_{y,0}$ and $\G_{y',0}$
are proper parahoric subgroups in $\G_{x,0}$. Let $\bM$ be 
the Levi subgroup adapted for $\G_{y,0}$ with $x,y,y'\in\Bd(\bM,k)$.
Suppose $\rM_{y,0}=\rM_{y',0}$, and
$\G_{y,0}$ and $\G_{y',0}$ project onto opposite parabolic subgroups
in $\G_{x,0}/\G_{x,0^+}$ with common Levi factor 
$(\G_{y,0}\cap\G_{y',0})/\G_{x,0^+}
\simeq\rM_{y,0}/\rM_{y,0^+}=\rM_{y',0}/\rM_{y',0^+}$. 
Suppose $\rP=\rM\rU$ and $\overline\rP=\rM\overline\rU$ are
two opposite parabolic subgroups in $\G$ so that 
$\G_{y,0}\cap\rP=\rM_{y,0}\rU_{y,0}$ projects onto $\G_{y,0}/\G_{x,0^+}$, 
and $\G_{y',0}\cap\overline\rP=\rM_{y',0}\overline\rU_{y',0}$ 
projects onto $\G_{y',0}/\G_{x,0^+}$.
\end{numbering}

The following is a corollary of the proof of \cite[(6.7)]{MP2}, which
in turn results from \cite[(6.1)]{MP2}. 
Note that $\G_{y,0}\cap\rU=\G_{y,0^+}\cap\rU$. 

\begin{lemma}\label{lem: fact2}
We keep the notation and situation from (\ref{num: fact2}). 
If a nonzero vector $v\in V_\pi$ is fixed under $\G_{y',0^+}$, 
then the following integral does not vanish:
\[
\int_{\rU\cap\G_{y,0}}\pi(n)v\,dn \ne0.
\]
Let $\chi$ be a quasi character of $\G$ and let $\pi_\chi:=\chi\otimes\pi$. 
Suppose $\chi$ is trivial on $\rU$ and $\overline\rU$.
Then, if $v$ is a nonzero vector in $V_{\pi_\chi}^{(\G_{y',0^+},\chi)}$,
\[
\int_{\rU\cap\G_{y,0}}\pi_\chi(n)v\,dn=\int_{\rU\cap\G_{y,0}}\pi(n)v\,dn
\ne0.
\]
\end{lemma}

\section{\bf Heisenberg representation}

\begin{numbering}{\bf The case $(\bG',\bG)$.} \rm  \ 

\noindent
Let $(\bG',\bG)$ be a tamely ramified twisted Levi sequence, and
let  $\yo\in\Bd(\bG',k)$. Let $\bT$, $\bM$, $\bU$, $\overline\bU$
and $\phi$ be as in (\ref{num: prepare}), and let $\bM'=\bG\cap\bM$.
Define concave functions
$\heish_+,\ \heish:\Phi(\bG,\bT,E)\cup\{0\}\longrightarrow \tbR$ as 
\[
\heish_+(a)=\left\{
\begin{array}{ll}
\rtm&\textrm{if }a\in\Phi(\bM,\bT,E)\cup\{0\},\\
\stm&\textrm{if }a\in\Phi(\bG',\bT,E)\setminus\Phi(\bM,\bT,E),\\ 
\stm^+&\textrm{otherwise},
\end{array}\right.
\qquad
\heish(a)=\left\{
\begin{array}{ll}
\rtm&\textrm{if }a\in\Phi(\bM,\bT,E)\cup\{0\},\\
\stm&\textrm{if }a\in\Phi(\bG',\bT,E)\setminus\Phi(\bM,\bT,E),\\ 
\stm&\textrm{otherwise}.
\end{array}\right.
\]
Let $\tilde J_+:=\G_{\yo,\heish_+}$ and $\tJ:=\G_{\yo,\heish}$.
Note that $\tJ\cap U=\G_{\yo,\stm}\cap U$.
Then since 
$\phi$ defines a character $\hat\phi$ on $\G'_{\yo,0^+}\G_{\yo,\stm^+}$
(see \ref{num: construction Ktypes}), $\hat\phi$ can be restricted to 
$\tJ_+\subset\G'_{\yo,0^+}\G_{\yo,\stm^+}$.
Let $\rN:=\ker(\hat\phi)$.

\end{numbering}

\begin{lemma}\ 
\begin{enumerate}
\item
The pairing $\langle a,b\rangle=\hat\phi(aba^{-1}b^{-1})$
defined on $\tJ/\tJ_+\times\tJ/\tJ_+$ is nondegenerate.
\item
$\tJ/\rN$ is a Heisenberg $p$-group with center $\tJ_+/\rN$.
\end{enumerate}
\end{lemma}

The proof is similar as that of \cite[(11.1)]{yu}. We would not repeat it here.

\begin{lemma}\label{lem: fact1}\ 
Let $(\phi^\rh,V_{\phi^\rh})$ be 
a representation of the Heisenberg group $\tJ/\rN$.
We use the same notation $\phi^\rh$ for the inflated representation of $\tJ$.
Let $v\in V_{\phi^\rh}$. If $v\ne0$ and it is fixed
by $\tJ\cap\overline\rU$, then the following integral does not vanish:
\[
\int_{\tJ\cap\rU}\phi^\rh(u)v\,du\ne0.
\]
\end{lemma}

\proof 
Note that $(\tJ\cap\overline\rU)\left/(\tJ_+\cap\overline\rU)\right.$ and 
$(\tJ\cap\rU)\left/(\tJ_+\cap\rU)\right.$ are isotropic subspaces
of $\tJ/\tJ_+$ such that their sum form a complete polarization of $\tJ/\tJ_+$.
Then the above lemma is a result of the following general fact on 
Heisenberg representations:

\begin{lemma}
Let $\bbH=\bbW\oplus\bbF_q$ be a Heisenberg group where $\bbW$ is
a symplectic space over $\bbF_q$. 
Let $\bbW=W_1\oplus W_2$ be 
a complete polarization of $\bbW$. Let $(\rho,V_\rho)$ be a finite dimensional
representation of $\bbH$.
If a nonzero vector $v\in V_{\rho}$ is fixed under $W_1$, 
then 
\[
\int_{W_2}\rho(\ty)v\,d\ty\ne0.
\]
\end{lemma}

\proof
Without loss of generality, we may assume that $\rho$ is irreducible.
If $\rho$ is the trivial representation, the lemma is trivial.
Otherwise, there is a nontrivial character $\chi$ of
the center $\bbF_q$ such that  $\rho$ is isomorphic to 
the Heisenberg representation $\chi^\rh$ 
of $\bbH$ with the central character $\chi$.
Then, we can realize $\chi^\rh$ on the function space $C(W_1)$ 
as follows: for $f\in C(W_1)$,
\[
\begin{array}{ll}
\chi^\rh(\tx)f(\txx)=f(\txx+\tx), &\textrm{for}\ \tx\in W_1,\\
\chi^\rh(\ty)f(\txx)=\chi(\frac12\langle \txx,\ty\rangle)f(\txx),\quad
&\textrm{for}\ \ty\in W_2,\\
\chi^\rh(z)f(\txx)=\chi(z)f(\txx),&\textrm{for}\ z\in\bbF_q.
\end{array}
\]
Denote the characteristic function supported at $\txx\in W_1$ by $f_{\txx}$.
Suppose $f\in C(W_1)$ is a nonzero function invariant by $W_1$. 
Then from the above formula,
$f$ is a nonzero constant function, that is, $f=c\,\sum_{\txx\in W_1} f_w$ for
some constant $c\ne0$.
Now, we have 

\medskip

$\begin{array}{ll}
\int_{W_2}\chi^\rh(\ty)\,f\,d\ty
&=\,\int_{W_2}\chi^\rh(\ty) (\sum_{\txx} c\,f_{\txx})\,d\ty
=\,\int_{W_2}\sum_{\txx}
c\,(\chi(\frac12\langle \txx,\ty\rangle)f_{\txx})\,d\ty\\
&=\,\int_{W_2} c\, f_0\,d\ty
=\,\textrm\vol(W_2)\cdot c\cdot f_0\ne0.
\end{array}$
\qed

\begin{corollary}\label{cor: fact1}\ 
We keep the notation from Lemma \ref{lem: fact1}.
Let $\chi$ be a character of $\tJ$ such that $\chi$ is trivial
on $\tJ\cap\rU$ and $\tJ\cap\overline\rU$.
Let $v\in V_{\chi\otimes\phi^\rh}$. If $v\ne0$ and it is fixed
by $\tJ\cap\overline\rU$, then 
\[
\int_{\tJ\cap\rU}(\chi\otimes\phi^\rh)(u)v\,du\ne0.
\]
\end{corollary}

\begin{numbering} \label{num: Jv2} {\bf The general case.} \rm
Let $\datum_\yo:=(\vec\bG,\yo,\vec\rtm,\vec\phi)$ be a good $\G$-datum
of positive depth. Let $\bT$, $\bM$, $\bU$ and $\overline\bU$ 
be as in (\ref{num: Jv}).
Recall $\bM^i=\bM\cap\bG^i$, $\bU^i=\bU\cap\bG^i$ and 
$\overline\bU^i=\overline\bU\cap\bG^i$.
For $i=1,\cdots,d$, define
concave functions $\heish^i_+,\ \heish^i:\Phi(\bG^i,\bT,E)\cup\{0\}\longrightarrow \tbR$ as 
\[
h^i_+(a)=\left\{
\begin{array}{ll}
\rtm_{i-1}&\textrm{if }a\in\Phi(\bM^i,\bT,E)\cup\{0\},\\
\stm_{i-1}&\textrm{if }
a\in\Phi(\bG^{i-1},\bT,E)\setminus\Phi(\bM^i,\bT,E),\\ 
\stm_{i-1}^+&\textrm{otherwise},
\end{array}\right.
\]
and
\[
h^i(a)=\left\{
\begin{array}{ll}
\rtm_{i-1}&\textrm{if }a\in\Phi(\bM^i,\bT,E)\cup\{0\},\\
\stm_{i-1}&\textrm{if }
a\in\Phi(\bG^{i-1},\bT,E)\setminus\Phi(\bM^i,\bT,E),\\ 
\stm_{i-1}&\textrm{otherwise}.
\end{array}\right.
\]
We also define corresponding open compact subgroups of $\G^i$
as follows: for $i=1,\cdots,d$,
\[
\tJ^i_{\yo+}=:\G^i_{\yo,h^i_+}\subset\G^i,\qquad
\tJ^i_{\yo}=:\G^i_{\yo,h^i}\subset\G^i.
\]
Note that 
$\tJ^i\cap U=\G^i_{\yo,\stm_{i-1}}\cap U$.
Then, $\phi_{i-1}$ defines a character $\hat\phi_{i-1}$ of $\tJ^i_{\yo+}$
as in (\ref{num: construction Ktypes}), 
and each $\tJ^i_\yo/\ker(\hat\phi_{i-1})$ is 
a Heisenberg $p$-group with center $\tJ^i_{\yo+}/\ker(\hat\phi_{i-1})$.
When $\datum$ is strongly good, $\hat\phi_{i-1}$ coincides with 
$\chi_{\gss_{i-1}}$.
\end{numbering}

\section{\bf Supercuspidality II}\label{sec: cusp II}

In this section,
let $(\pi,V_\pi)$ be a fixed supercuspidal representation of positive depth.
Suppose $\datum=(\vec\bG,\yo,\vec\rtm,\vec\phi)$ is
a strongly good $\G$-datum of positive depth such that 
$(\rK^+_\datum,\phi_\datum=\chi_\gss)$ is contained in $(\pi,V_\pi)$. 
Then, by Proposition \ref{prop: cusp I}, $\bZ^0/\bZ^d$ is $k$-anisotropic.

\begin{numbering}\label{num: depth zero data}\rm
Consider the $(\rK^+_\datum,\chi_\sgss)$ isotypic component 
$V_\pi^{\chi_\gss}$ of $V_\pi$.
It can be found as the image of the projection
\[
v\longrightarrow 
\frac1{\vol(\rK^+_{\datum})}\int_{\rK^+_{\datum}}
\chi_\gss(g^{-1})\pi(g)v\,dg.
\]
Recall $\G^0_{[\yo]}$ denotes the stabilizer in $\G^0$ 
of the image $[\xo]$ of $\xo$ in the reduced building of $\bG^0$.
Since $\G^0_{[\yo]}$ stabilizes $(\rK^+_{\datum},\chi_\gss)$,
$V_\pi^{\chi_\gss}$ is a finite dimensional representation 
of $\G^0_{[\yo]}$. 
Let $\overline\cuspd$ be an irreducible subrepresentation of 
$\G^0_{[\yo]}$ in $V_\pi^{\chi_\sgss}$. Denote the
character of $\G^0_{[y]}$ given by $\prod_i(\phi_i|\G^0_{[y]})$ again
by $\phi_\datum$. 
Then, $\cuspd=\phi_\datum^{-1}\otimes\overline\cuspd$ is an irreducible
representation of $\G^0_{[\yo]}$ which factors through $\G^0_{\yo,0^+}$.

Let $\cuspf$ be an irreducible subrepresentation of $\cuspd|\G^0_{\yo,0}$.
Then $\cuspf$ induces an irreducible representation 
of $\G^0_{\yo,0}/\G^0_{\yo,0^+}$, which we again denote by $\cuspf$.
\end{numbering}

\begin{proposition}\label{prop: cusp II}
Suppose $\Hypk$ and $\HypB$ are valid.
We keep the notation from above.
\begin{enumerate}
\item
$\G^0_{\yo,0}$ is a maximal parahoric subgroup of $\G^0$.
\item
$\cuspf$ is a cuspidal representation of 
$\G^0_{\yo,0}/\G^0_{\yo,0^+}$.
\end{enumerate}
\end{proposition}

\begin{numbering}\label{num: Uy}
{\it Proof of  Proposition \ref{prop: cusp II}-(1).} \rm
We generalize the proof of \cite[(6.7)]{MP2}.

Let $\bM$ be a $k$-Levi subgroup of $\bG$ 
adapted for $\G^0_{\yo,0}$. Suppose $\G^0_{y,0}$ is not maximal. Then,
$\bM$ is a proper Levi subgroup of $\bG$.
Let $\bC$ be the center of $\bM$.
Let $\bP$ be a $k$-parabolic subgroup with $\bP=\bM\bU$
and let $\overline\bU$ be the unipotent subgroup opposite to $\bU$.
Then, $\rU\cap\G^0_{\yo,0}=\rU\cap\G^0_{\yo,0^+}$ 
and $\overline\rU\cap\G^0_{\yo,0}=\overline\rU\cap\G^0_{\yo,0^+}.$ 
By Corollary \ref{cor: extend}-$(ii)$, 
$(\pi,V_\pi)$ contains $(\Jvy,\chi_\gss)$.

For any $\xx\in\Bd(\bG^0,k)\cap\Bd(\bM,k)$, let
$\datum_\xx=(\vec\bG,\xx,\vec\rtm,\vec\phi)$.
As in (\ref{num: Jv}), for $\datum_\xx$ and $\bM\bU$,
we can define $\rK_{\xx\vdash}$, and  
$\vec\phi$ defines a character $\chi_\sgss$ 
of $\rK_{\xx\vdash}$ realized by $\sgss$.

Let $\beta: \bGL_1\rightarrow \bC$ be 
a one-parameter subgroup of $\bC$ such that
$\langle\alpha,\beta\rangle>0$ for every root $\alpha$ of $\bC$
in the Lie algebra $\lieU$ of $\rU$.
Consider the ray $\yo(t)=\yo+t\beta$, $t\ge0$, contained in the apartment 
$\Apt(\bG,\bT,k)$, emanating from the point $\yo$ in the direction of $\beta$.
Then, one can verify the following (see also \cite[(6.7)]{MP2}):

\begin{enumerate}
\item
$\rM\cap \rK_{\yo(t)\vdash}=\rM_{\yo(t),\vec\stm^+}=\rM_{\yo,\vec\stm^+}$; \  
$\rU\cap\rK_{\yo(t)\vdash}=\rU_{\yo(t),\vec\stm}$; \ 
$\overline\rU\cap\rK_{\yo(t)\vdash}=\overline\rU_{\yo(t),\vec\stm^+}$.
\item 
If $t'\ge t$, then $\rU_{\yo(t'),\vec\stm}\supset \rU_{\yo(t),\vec\stm}$ and 
$\overline\rU_{\yo(t'),\vec\stm^+}\subset \overline\rU_{\yo(t),\vec\stm^+}$.
\item
$\rU^i_{\yo(\tm),\stm_{i-1}}=\tJ_{\yo(\tm)}^i\cap\rU$
and $\rU_{\yo(\tm),\vec\stm}=\rU^0_{\yo(\tm),0}\rU^1_{\yo(\tm),\stm_0}
\cdots\rU^d_{\yo(\tm),\stm_{d-1}}$. 
\item
Any compact subgroup of $\rU$ is contained in $\rU_{\yo(t),\vec\stm}$
for $\tm$ sufficiently large.
\item
There is a sequence $0=t_0<t_1<t_2<\cdots$ tending to $\infty$
so that $\rU_{\yo(t),\vec\stm},\,\rU_{\yo(t),\vec\stm^+},\,
\overline\rU_{\yo(t),\vec\stm}$ and $\overline\rU_{\yo(t),\vec\stm^+}$ 
are constant on the open intervals $t_{i-1}<t<t_i$ ($i\ge1$). 
Then, in fact we have that $\rU_{\yo(t),\vec\stm}$, 
$\overline\rU_{\yo(t),\vec\stm^+}$ are constant for $\tm_{i-1}\le t<\tm_i$,
and $\rU_{\yo(t),\vec\stm^+}$, $\overline\rU_{\yo(\tm),\vec\stm}$
are constant for $\tm_{i-1}<\tm\le\tm_i$.

Let $\xx_i=\yo+\tm_i\beta$. Fix $t'_i\in\bbR$ 
such that $t_{i-1}<t'_i<t_i$. Let $\yo_0=\xx_0$, and $\yo_i=\yo+t'_i\beta$. 
We observe that 
$\rU_{\xx_0,\vec\stm}=\rU_{\yo_0,\vec\stm}=\rU_{\yo_1,\vec\stm}$, 
$\overline\rU_{\xx_0,\vec\stm}=\overline\rU_{\yo_0,\vec\stm}
\supset\overline\rU_{\yo_1,\vec\stm}$, and
$\rU_{\xx_j,\vec\stm}=\rU_{\yo_{j+1},\vec\stm}$.
\item
$\rU_{\yo_i,\vec\stm^+}=\rU_{\yo_i,\vec\stm}$ and 
$\overline\rU_{\yo_i,\vec\stm^+}=\overline\rU_{\yo_i,\vec\stm}$, when $i\ge1$.

Proving the first equality, $\subset$ is obvious.
We may assume that $\gss$ is $k$-split without loss of generality. 
Let $\psi\in\Psi(\bG,\bT,k)$ such the gradient 
$\dot\psi$ of $\psi$ is in $\Phi(\bU^j,\bT,k)$, $j=0,\cdots,d$.
Let $\rU_\psi$ be the root subgroup associated to $\psi$. 
Suppose $\rU_\psi\subset\rU_{\yo_i,\vec\stm}$. 
Then, by the definition of $\yo_i$ in (5), 
$\rU_\psi\subset\rU_{\yo_i-\epsilon\beta,\vec\stm}$ for sufficiently 
small $\epsilon>0$. Hence, 
$\psi(\yo_i-\epsilon\beta)=\psi(\yo_i)-\epsilon\langle\dot\psi,\beta\rangle
\ge\stm_{j-1}$.
Since $\langle\dot\psi,\beta\rangle>0$, $\psi(\yo_i)>\stm_{j-1}$. 
Hence $\rU_\psi\subset\rU_{\yo_i,\vec\stm^+}$, and $\supset$ follows.
The second equality is similarly proved.
\item
$\tJ_{\xx_i}^\ell\cap\rU=\tJ_{\yo_{i+1}}^\ell\cap\rU
=\rU^\ell_{\yo_{i+1},\stm_{\ell-1}}$ and
$\tJ_{\xx_i}^\ell\cap\overline\rU=\tJ_{\yo_i}^\ell\cap\overline\rU
=\overline\rU^\ell_{\yo_i,\stm_{\ell-1}}$ ;
\newline
$\tJ_{\xx_i+}^\ell\cap\rU
=\rU^\ell_{\xx_i,\stm_{\ell-1}^+}\rU^{\ell-1}_{\xx_i,\stm_{\ell-1}}$ and 
$\tJ_{\xx_i+}^\ell\cap\overline\rU
=\overline\rU^\ell_{\xx_i,\stm_{\ell-1}^+}
\overline\rU^{\ell-1}_{\xx_i,\stm_{\ell-1}}$, \ 
$\ell=1,\cdots, d$.
\item
$\G^0_{\yo_i,0}=\overline\rU^0_{\yo_i,0}\rM^0_{\yo,0}\rU^0_{\yo_i,0}$,
$\G^0_{\yo_i,0^+}=\overline\rU^0_{\yo_i,0}\rM^0_{\yo,0^+}\rU^0_{\yo_i,0}$  and
$\G^0_{\xx_i,0^+}=\overline\rU^0_{\yo_{i+1},0}\rM^0_{\yo,0^+}\rU^0_{\yo_i,0}$.
\item
$\G^0_{\xx_i,0}$ contains both $\G^0_{\yo_i,0}$ and $\G^0_{\yo_{i+1},0}$. 
Suppose $\G^0_{\yo_i,0}\ne\G^0_{\yo_{i+1},0}$. Then the images of 
$\G^0_{\yo_i,0}$ and $\G^0_{\yo_{i+1},0}$ form the opposite parabolic subgroups
in $\G^0_{\xx_i,0}/\G^0_{\xx_i,0^+}$ with the common Levi factor 
given by the image of $\rM^0_{\yo,0}=\rM^0_{\yo_i,0}=\rM^0_{\yo_{i+1},0}$. 
\end{enumerate}

Let $v=v_0\in V_\pi^{\chi_\sgss}$ be a nonzero element.
Note that $\rU_{\xx_j,\vec\stm}=\rU_{\yo_{j+1},\vec\stm}$ by (5).
For $i=1,2,\cdots$, define 
\[
v_j=\int_{\rU_{\xx_j,\vec\stm}}\pi(n)v_{j-1}\,dn
=\int_{\rU_{\yo_{j+1},\vec\stm}}\pi(n)v_{j-1}\,dn.
\]
Inductively, we will show 

$(i)$ $v_j\ne0$;

$(ii)$ $v_j$ is a nonzero multiple of $\int_{\rU_{\xx_j,\vec\stm}}\pi(n)v\,dn$;

$(iii)$ $\rM_{\yo,\vec\stm^+}$ acts on $v_j$ 
as a character represented by $\gss$;

$(iv)$ for $j\ge0$, we have 
$Stab(v_{j})\supset\overline\rU_{\xx_{j+1},\vec\stm}
=\overline U_{\yo_{j+1},\vec\stm^+}$ and 
$Stab(v_{j})\supset\rU_{\xx_{j+1},\vec\stm^+}=\rU_{\yo_{j+1},\vec\stm^+}$;

$(v)$ $\G^0_{\xx_j,0^+}$ acts on $v_j$ 
as a character represented by $\gss$; 

$(vi)$ $v_j=\ 
\int_{\rU^0_{\xx_j,0}}\int_{\rU^1_{\xx_j,\stm_0}}
\cdots\int_{\rU^d_{\xx_j,\stm_{d-1}}}
\pi(u_0u_1\cdots u_d)v_{j-1}\,du_d\cdots du_0.$

\smallskip

The equality in $(ii)$ follows from Fubini's theorem and the bi-invariance of
the Haar measure on the unipotent subgroup $\rU$.

Since the compact group $\rM_{\yo,\vec\stm^+}=\rM_{\yo_{j+1},\vec\stm^+}$ 
normalizes $\rU_{\yo_{j+1},\vec\stm}$  
and acts on $v$ as a character $\chi_\gss$, $(iii)$ follows from $(ii)$. 

Proving $(iv)$, the second inclusion is obvious.
To prove the first, since $v_{0}\in V_\pi^{\chi_\sgss}$, we note that
$\rU_{\yo_0,\vec\stm}\cup \overline\rU_{\yo_0,\vec\stm^+}
\subset \ker(\chi_\gss|\rK_{\yo_0\vdash})=\ Stab(v_{0})$.
Since $\rU_{\yo_1,\vec\stm^+}\subset\rU_{\yo_1,\vec\stm}=\rU_{\yo_0,\vec\stm}$
and $\overline\rU_{\yo_1,\vec\stm^+}\subset\overline\rU_{\yo_0,\vec\stm^+}$,
the case $j=0$ follows. Now, suppose $j\ge1$.
We observe that 
$\overline\rU_{\xx_{j+1},\vec\stm}=\overline\rU_{\yo_{j+1},\vec\stm^+}\subset
\overline\rU_{\yo_j,\vec\stm^+}\subset Stab(v_{j-1})$ by 
induction hypotheses.  
Consider the character $\chi_\gss$ realized by $\gss$ on 
$\rK_{\yo_{j+1}\vdash}$.
Denote the kernel of $\chi_\gss$ in $\rK_{\yo_{j+1}\vdash}$ by $\rA$.
Then since $\rU_{\yo_{j+1},\vec\stm}, \ \overline\rU_{\yo_{j+1},\vec\stm^+}
\subset \rA$ and $\rK_{\yo_{j+1}\vdash}$ is decomposable 
with respect to $P=MU$, 
$\rA$ is also decomposable with respect to $P=MU$. Let 
$\rA_M:=\rA\cap\rM$. Then we have $\rA_M\subset Stab(v_{j-1})$ by $(iii)$, 
and $\rA=\rU_{\yo_{j+1},\vec\stm}\rA_M\overline\rU_{\yo_{j+1},\vec\stm^+}
=\overline\rU_{\yo_{j+1},\vec\stm^+}\rA_M\rU_{\yo_{j+1},\vec\stm}$.
Since $\overline\rU_{\yo_{j+1},\vec\stm^+}\rA_M
=\rA_M\overline\rU_{\yo_{j+1},\vec\stm^+}\subset Stab(v_{j-1})$, for
a proper normalization, we have
\[
v_{j}
=\int_{\rU_{\yo_{j+1},\vec\stm}}\pi(n)v_{j-1}\,dn
=\int_{\rU_{\yo_{j+1},\vec\stm}}\int_{\rA_M}
\int_{\overline\rU_{\yo_{j+1},\vec\stm^+}}\pi(nm\overline n)v_{j-1}\,
d\overline n\,dm\,dn=\int_{\rA}\pi(a)\,da.
\]
Hence $\overline\rU_{\yo_{j+1},\vec\stm^+}\subset\rA \subset Stab(v_j)$ 
and $(iv)$ follows.

For $(v)$, if $j=0$, it is obvious from the choice of $v_0$.
If $j\ge1$, since $\G^0_{\xx_j,0^+}
=\overline\rU^0_{\yo_{j+1},0}\rM^0_{\yo,0^+}\rU^0_{\yo_j,0}$,
and since $v_j$ is fixed by $\rU^0_{\yo_j,0}$ and
$\overline\rU^0_{\yo_{j+1},0}=\overline\rU^0_{\yo_{j+1},0^+}$,
combining these with $(iii)$, $(v)$ follows.

The equalities in $(vi)$ result from the decomposition
$\rU_{\xx_j,\vec\stm}
=\rU^0_{\xx_j,0}\rU^1_{\xx_j,s_0}\cdots\rU^d_{\xx_j,s_{d-1}}$.

To prove $(i)$, let $w_{d+1}:=v_{j-1}$ and 
$w_i:=\int_{\rU^i_{\xx_j,\stm_{i-1}}}\pi(u_i)w_{i+1}\,du_i$.
Then $v_j=w_0$.
As in the case of $(iii)$ and $(iv)$, we have 
for $\ell\le i$,
$\rM^\ell_{\yo,\stm^+_{\ell-1}}$ acts on $w_i$ via $\chi_\gss$ and 
each $w_i$ is fixed by $\overline\rU^\ell_{\xx_{j},\stm_{\ell-1}}$.
By induction hypothesis, $v_{j-1}=w_{d+1}\ne0$.
Suppose $w_{i+1}\ne0$ with $i\ne0$. Then 
consider the $\tJ^i_{\xx_j}$ stable vector space $W_i$ generated
by $\{\pi(g)w_{i+1}\mid g\in\tJ^i_{\xx_j}\}$. Note that $w_{i+1}$ is fixed by
$\overline\rU^i_{\xx_j,\stm_{i-1}}$ and $\rU^i_{\xx_j,\stm^+_{i-1}}$.
Moreover, 
from $\rU^{i-1}_{\xx_j,\stm_{i-1}}\subset\rU^{i-1}_{\xx_j,\stm_{i-2}^+}
\subset\rU_{\xx_j,\vec\stm^+}\subset Stab(v_{j-1})$, we have that
$\rU^{i-1}_{\xx_j,\stm_{i-1}}$ fixes $w_{i+1}$.
Hence, the representation of $\tJ^i_{\xx_j}$ on $W_i$ is a
Heisenberg representation of $\tJ^i_{\xx_j}$ twisted by
a character $\chi_\gamma$ of $\tJ^i_{\xx_j}$ represented by 
$\gamma=\gss_i+\gss_{i+1}+\cdots+\gss_d$
(recall each $\gss_i,\gss_{i+1},\cdots,\gss_d$
defines a character of $\G^{i}$). Then by Corollary \ref{cor: fact1},
$w_i\ne0$. Inductively, $w_1\ne0$. Since $w_1$ is fixed
under $\overline\rU^0_{\xx_{j},0}$ and $\rU^0_{\xx_{j},0^+}$ and
since $\overline\rU^0_{\xx_{j},0}$ and $\rU^0_{\xx_{j},0}$ project to 
the unipotent radicals of opposite parabolic subgroups in 
$\G^0_{\xx_{j},0}/\G^0_{\xx_{j},0^+}$, by Lemma \ref{lem: fact2},
$w_0=v_j=\int_{\rU^0_{\xx_{j},0}}\pi(u_0)w_1\,du_0$ is not zero. 

Since any open compact subgroup of $\rU$ is contained 
in $\rU_{\xx_{j},\vec\stm}$ for sufficiently large $j$, by $(i)$ and $(ii)$, 
we conclude that the image of $v$ under the Jacquet module map  
with respect to $\rU$ is nonzero. This is a contradiction to
the supercuspidality of $(\pi,V_\pi)$. 
Hence, we conclude that $\G^0_{\yo,0}$ is maximal.
\qed
\end{numbering}

\begin{numbering} {\it Proof of Proposition \ref{prop: cusp II}-(2).} \rm
Suppose $\cuspf$ is not cuspidal. Then there is $z\in\Bd(\bG^0,k)$
such that $\G^0_{z,0}\subsetneq\G^0_{\yo,0}$ and 
$\phi_\datum^{-1}\otimes V_\pi^{\chi_\sgss}$ has
$\G^0_{z,0^+}$ invariants. 
Let $\bT\subset\bG^0$ be an $E$-split maximal $k$-torus such that 
$\Apt(\bG^0,\bT,k)$ is of maximal dimension and $\yo,z\in\Apt(\bG^0,\bT,k)$.
Let $\bM$ be a $k$-Levi subgroup of $\bG$ with $\bT\subset\bM$
adapted for $\G^0_{z,0}$. Let $\bC$ be the center of $\bM$.
Fix a $k$-parabolic subgroup $\bP$ with $\bP=\bM\bU$.

Since $\phi_\datum^{-1}\otimes V_\pi^{\chi_\sgss}$ has $\G^0_{z,0^+}$
invariants, and since $\rM^0_{\yo,0}\rU^0_{\yo,0}$ projects to
a parabolic subgroup in $\G^0_{\yo,0}/\G^0_{\yo,0^+}$ sharing 
the Levi factor with the image of $\G^0_{z,0}$, 
we can deduce  from \cite{HCf} (see also \cite[(1.1)]{HL}) that 
$\phi_\datum^{-1}\otimes V_\pi^{\chi_\sgss}$ also has
$\rU^0_{\yo,0}$ invariants. Then, by Corollary \ref{cor: extend}-$(ii)$, 
$(\pi,V_\pi)$ also contains $(\rK_{\yo\vdash},\chi_\gss)$ where 
$\rK_{\yo\vdash}$ is defined with respect to $\bM\bU$ as in (\ref{num: Jv}).

Let $\beta: \bGL_1\rightarrow \bC$ be 
a one-parameter subgroup of $\bC$ such that
$\langle\alpha,\beta\rangle>0$ for every root $\alpha$ of $\bC$
in the Lie algebra $\lieU$ of $\rU$.
Take the ray $\yo(t)=\yo+t\beta$, $t\ge0$, contained in the apartment 
$\Apt(\bG,\bT,k)$, emanating from the point $\yo$ in the direction of $\beta$.
Now, the rest of the proof is similar to that of 
Proposition \ref{prop: cusp II}-(1).
\qed
\end{numbering}

\section{\bf Extending $\phi_\datum$ of $\rK^+_\datum$ to $\rK^d$}
\label{sec: extend}

Let $(\pi,V_\pi)$ be a fixed irreducible supercuspidal representation of
positive depth.
Suppose $\datum=(\vec\bG,\yo,\vec\rtm,\vec\phi)$ is
a good $\G$-datum of positive depth such that 
$(\rK^+_\datum,\phi_\datum)$ is contained in $(\pi,V_\pi)$. 
Define
\[
\rK:=\vec\bG(k)_{(0^+,\stm_0,\stm_1,\cdots,\stm_{d-1})}
=\G^0_{\yo,0^+}\G^1_{\yo,\stm_0}\cdots\G^d_{\yo,\stm_{d-1}}.
\]
Note that $\rK^d=\G^0_{[y]}\rK$ where $\rK^d$ is defined in 
(\ref{num: notation yu}).
To simplify notation, we denote $\rK^+_\datum$ by $\rK^+$.
Let $N:=\ker(\phi_\datum)$.

\begin{proposition}\label{prop: heisII}
$\rK/N$ is a Heisenberg $p$-group with center $\rK^+/N$.
\end{proposition}

To prove the above proposition, recall from 
(\ref{num: notation yu}) that 
\[
J^i_+=(\bG^{i-1},\bG^i)(k)_{\yo,(\rtm_{i-1},\stm^+_{i-1})},\qquad
J^i=(\bG^{i-1},\bG^i)(k)_{\yo,(\rtm_{i-1},\stm_{i-1})}.
\]
Define a pairing on $J^i/J^i_+$ given by 
$\langle a,b\rangle_i=\phi_\datum(aba^{-1}b^{-1})$.

Note that 
$\rK=\G^0_{\yo,0^+}\rJ^1\cdots\rJ^d$ \ and \ 
$\rK^+=\G^0_{\yo,0^+}\rJ^1_+\cdots\rJ^d_+$.

\begin{lemma}\label{lem: heisII} \ 
\begin{enumerate}
\item
The pairing $\langle\ ,\ \rangle_i$ is nondegenerate on $J^i/J^i_+$.
\item
$J^i/(N\cap J^i_+)$ is a Heisenberg $p$-group with center 
$J^i_+/(N\cap J^i_+)$.
\end{enumerate}
\end{lemma}

\proof
Let $\hat\phi_j$ be as in (\ref{num: construction Ktypes}) 
or (\ref{num: yusc}).
For $j>i-1$, since $aba^{-1}b^{-1}\in (G^j,G^j)$, 
$\hat\phi_{j}(aba^{-1}b^{-1})=1$.
For $j<i-1$, since $aba^{-1}b^{-1}\in\G^i_{\yo,\rtm_{i-1}}
\subset\G_{\yo,\rtm_j^+}\subset\ker(\hat\phi_j)$,
$\hat\phi_j(aba^{-1}b^{-1})=1$.
Hence $\langle a,b\rangle=\phi_\datum(aba^{-1}b^{-1})
=\hat\phi_{j-1}(aba^{-1}b^{-1})$. Then this pairing coincides with 
the one in Lemma 11.1 of \cite{yu}, and the above follows from that lemma.
\qed

\medskip

\noindent{\it Proof of Proposition \ref{prop: heisII}.}
Define a pairing on $\rK/\rK^+$ given by 
$\langle a,b\rangle=\phi_\datum(aba^{-1}b^{-1})$.
Note that $\langle\ ,\ \rangle$ on $J^i/J^i_+\subset\rK/\rK^+$ 
coincides with $\langle\ ,\ \rangle_i$.
We claim that if $i\ne j$, $\langle\ ,\ \rangle$ is trivial on 
$\rJ^i/\rJ^i_+\times \rJ^j/\rJ^j_+$. Without loss of generality, 
we may assume $i<j$. Let $a\in\rJ^i$ and $b\in\rJ^j$. If $\ell<j$, 
since $(\rJ^i,\rJ^j)\subset \rJ^j_+\subset\ker(\hat\phi_\ell)$, 
$\hat\phi_\ell(aba^{-1}b^{-1})=1$. 
If $\ell\ge j$, $aba^{-1}b^{-1}\in (\G^j,\G^j)\subset(\G^\ell,\G^\ell)$ and 
$\hat\phi_\ell(aba^{-1}b^{-1})=1$. 

Hence, we have
$\langle\ ,\ \rangle
=\langle\ ,\ \rangle_1\oplus\cdots\oplus\langle\ ,\ \rangle_d$,
and the pairing $\langle\ ,\ \rangle$ is nondegenerate.
Now, we can conclude $\rK/N$ is a Heisenberg $p$-group with 
center $\rK^+/N$.
\qed

\begin{numbering}\label{num: last ext1} \rm 
By Proposition \ref{prop: heisII} and the theory of Heisenberg groups,
there is a unique representation of $\rK$ extending the character
$\phi_\datum$. Let $\phi_\datum^\rh$ be the unique irreducible
representation of $\rK$ extending the character $\phi_\datum$, and 
let $V_{\phi_\datum^\rh}\subset V_\pi^{\phi_\datum}$ be 
the space of $\phi_\datum^\rh$. 

Let $\omega_\pi$ be the central character of $(\pi,V_\pi)$. 
Then, $Z_\G$ acts via $\omega_\pi$ on $V_{\phi_\datum^\rh}$.
Let $\phi_\datum^{\rh\pi}$ denote
the irreducible representation of $Z_\G\rK$ such that 
$\phi_\datum^{\rh\pi}=\phi_\datum^\rh$ on $K$ and
$\phi_\datum^{\rh\pi}=\omega_\pi\cdot 1$ on $Z_\G$.
Note that $(\pi,V_\pi)$ also contains $(Z_\G\rK,\phi_\datum^{\rh\pi})$.
Moreover, $\G^0_{[\yo]}$ stabilizes the representation 
$(\phi_\datum^\rh,V_{\phi_\datum^\rh})=
(\phi_\datum^{\rh\pi},V_{\phi_\datum^{\rh\pi}})$.
We fix an extension of $\phi_\datum^{\rh\pi}$ of $Z_\G\rK$
to $\G^0_{[\yo]}\rK$ on $V_{\phi_\datum^\rh}$, 
which we denote by $\tilde\phi_\datum^{\rh\pi}$.
\end{numbering}

\begin{remark} \label{rmk: kappa}\rm
Let $\kappa$ be as in (\ref{num: yusc}) associated to $\datum$. 
Then since $\kappa|\rK^+$ is
$\phi_\datum$-isotypic, $\kappa|\rK\simeq\phi_\datum^\rh$.
Hence we can choose the extension $\tilde\phi_\datum^{\rh\pi}$ to be $\kappa$.
\end{remark}

\begin{proposition}\label{prop: last ext}
Any irreducible representation of $\G^0_{[\yo]}\rK$ extending 
$(Z_\G\rK,\phi^{\rh\,\pi}_\datum)$ is of the form 
$\rho\otimes\tilde\phi_\datum^{\rh\pi}$
where $\rho$ is an irreducible representation of $\G^0_{[\yo]}\rK$ 
factoring through $Z_\G\rK$.
\end{proposition}

\proof 
Note that $\rho$ in the above induces an irreducible representation of 
$\G^0_{[\yo]}/\G^0_{\yo,0^+}$.
Then, since $\tilde\phi_\datum^{\rh\pi}|K$ is already irreducible,
the representation of the form $\rho\otimes\tilde\phi_\datum^{\rh\pi}$ 
is also irreducible.

Let $\Theta_1$ (resp. $\Theta_2$) be the character of 
$\lp\textrm{Ind}_{Z_\G\rK}^{\G^0_{[\yo]}\rK}1\rp\,
\otimes\tilde\phi_\datum^{\rh\pi}$
(resp. $\textrm{Ind}_{Z_\G\rK}^{\G^0_{[\yo]}\rK}\phi_\datum^{\rh\pi}$).
Now, since $[\G^0_{[\yo]}:\G^0_{\yo,0^+}Z_\G]$ is finite, we can explicitly
compute $\Theta_1$ and $\Theta_2$.
Since $\G^0_{[\yo]}\rK$ normalizes $Z_\G\rK$,
for $g, h\in\G^0_{[\yo]}\rK$, $^hg\in Z_\G\rK$  
if and only if $g\in Z_\G\rK$.
Hence,
\[
\Theta_1(g)=\Theta_2(g)=\left\{
\begin{array}{ll} 
\phi_\datum^{\rh\pi}(g)&\textrm{if }g\in Z_\G\rK, \\
0&\textrm{otherwise.}
\end{array}\right.
\]
Therefore, 
$\lp\textrm{Ind}_{Z_\G\rK}^{\G^0_{[\yo]}\rK}1\rp\,
\otimes\tilde\phi_\datum^{\rh\pi}$ 
and $\textrm{Ind}_{Z_\G\rK}^{\G^0_{[\yo]}\rK}\phi_\datum^{\rh\pi}$ 
are isomorphic, and the proposition follows from this.
\qed

\begin{corollary}\label{cor: last ext}
We keep the notation from Proposition \ref{prop: last ext}.
Let $(\pi,V_\pi)$ be as before.
Then there is a strongly good $\G$-datum of positive depth
$\datum=(\vec\bG,\yo,\vec\rtm,\vec\phi)$ and 
an irreducible representation $\rho$ of $\G^0_{[\yo]}\rK$ factoring 
through $Z_\G\rK$ such that $(\pi,V_\pi)$ contains
$(\G^0_{[\yo]}\rK,\rho\otimes\kappa)$, where $\kappa$ is as in 
(\ref{num: yusc}). Moreover, 
$c\textrm{-}\mathrm{Ind}_{\G^0_{[\yo]}}^{\G^0}\rho$ 
is an irreducible supercuspidal representation of $\G^0$
\end{corollary}

\proof 
The first statement follows from Corollary \ref{cor: tempered} and 
(\ref{num: last ext1})--(\ref{prop: last ext}).
The second follows from Proposition \ref{prop: cusp II} and \cite[(6.6)]{MP2}.
\qed

\section{\bf Conclusion: Exhaustion of supercuspidal representations}

\begin{theorem}\label{thm: Main} 
Suppose $\Hypk$, $\HypB$, $\HypGT$ and $\HypN$ are valid.
Then any irreducible supercuspidal representation $(\pi,V_\pi)$ of $\G$
arises through Yu's construction.
\end{theorem}

\proof 
Let $(\pi,V_\pi)$ be a given irreducible supercuspidal representation.
If the depth $\dpi(\pi)$ of $(\pi,V_\pi)$ is zero, 
this follows from Remark \ref{rmk: depth zero exhaustion}. 
Hence, we may assume $\dpi(\pi)>0$.
Then by Corollary \ref{cor: last ext}, 
there is a strongly good $\G$-datum of positive depth
$\datum=(\vec\bG,\yo,\vec\rtm,\vec\phi)$ and 
an irreducible representation $\rho$ of $\G^0_{[\yo]}\rK$ factoring 
through $Z_\G\rK$ such that $(\pi,V_\pi)$ contains
$(\G^0_{[\yo]}\rK,\rho\otimes\kappa)$. 
Moreover, $\bZ_{\bG^0}/\bZ_{\bG}$ is $k$-anisotropic 
by Proposition \ref{prop: cusp I}, and 
$c\textrm{-Ind}_{\G^0_{[\yo]}}^{\G^0}\rho$ is
an irreducible supercuspidal representation of $\G^0$ 
from Corollary \ref{cor: last ext}. Hence,
$\datum_Y=(\vec\G,\yo,\vec\rtm,\vec\phi,\rho)$ is a generic $\G$-datum.
Let 
$\pi_{\datum_Y}:=\,c\textrm{-Ind}_{\G^0_{[\yo]}\rK}^\G(\rho\otimes\kappa)$, 
the supercuspidal representation constructed in \cite{yu}. 
Since $\rho\otimes\kappa< \pi$, 
from the Frobenius reciprocity and the irreducibility of $(\pi,V_\pi)$,
we conclude $\pi\simeq\pi_{\datum_Y}$.
\qed

\begin{corollary}
Suppose $\Hypk$, $\HypB$, $\HypGT$ and $\HypN$ are valid.
Then, all supercuspidal representations of $\G$ are compactly induced from  
an irreducible representation of an open compact mod center group.
\end{corollary}

\def\noopsort#1{}
\providecommand{\bysame}{\leavevmode\hbox to3em{\hrulefill}\thinspace}

\section*{Selected notation}

\vskip .3cm

$\begin{array}{ll}
\S1.1&\bG^\der\\
   &\Bd(\bG,E)\\
   &\Apt(\bG,\bT,E) \\
\S1.2&\Bd(\G)\\
   &\lieG_\rtm\\
\S1.3&\Phi(\bG,\bT,E),\ \Psi(\bG,\bT,E)\\
   &\bX_\ast(\bT,E),\ \bX^\ast(\bT,E)\\
\S2.1&\vec\bG\\
\S2.2&\Apt(\bH,\bT,k)\\
\S2.5& \orb^H_X,\ \orb_X\\
   &\orb(H,X),\ \orb(X)\\
\S3.1& \G_{y,f},\ \lieG_{y,f} \\
   & f_{\vec\utm} \\
   & \vec\G_{y,\vec\utm},\ \vec\lieG_{y,\vec\utm} \\
   & \lieG^{i\,\perp},\ \lieG^i_\perp,\ \lieG^{i\,\perp}_{x,\rtm}\\
\S3.4&\HypB,\ \HypGT,\ \Hypk,\ \HypN\\
   &\bilinear(X,Y) \\
\S3.7&\lieG^{i\perp},\ \lieG^i_\perp\\
\S3.8& \rtm^\ast,\ \vec\rtm^\ast\\
   & \vec\lieG_{y,\vec\utm^\ast} \\
   & \CaL^{\,\ast} \\
\S3.10& \chi_\gamma\\
\S4.1&\dth(X)\\
\S5.1& \datum=(\vec\bG,\xo,\vec\rtm,\vec\phi)\\
   &\ell(\datum)\\
   &\bG^{-1},\ \bG^{d+1},\ \bZ^i, \\
\S5.3&\gss^i,\ \gss_\datum=\gss \\
\S5.8&\stm_i \\
   &\vec\stm(\epsilon),\ \vec\stm^+(\epsilon),\ \vec\stm(0),\ \vec\stm(0^+)\\ 
   &\vec\stm,\ \vec\stm^+\\
   &\rK_\datum^{i+},\ \rK_{\datum+}^i,\ \rK^+_\datum\\
   &\hat\phi_i\\
\S5.9&\phi_\datum,\ \chi_\gss \\
   &\gss^{-1} \\
\S5.10&\datum_x\\
\S5.11&\G^0_{[\xo]}\\
\S5.12&\CaL_{\xx,\epsilon}=\CaL_{\xx},\ \rL_{\xx,\epsilon}=\rL_{\xx}\\
   & \CaL_{\xx,\epsilon}^\sharp=\CaL_{\xx}^\sharp\\
   &\CaL_x^i,\ \CaL_x^{i \sharp}\\
\end{array}$ \hskip 2cm
$\begin{array}{ll}
\S5.16&\admG,\ \dualG,\ \temG\\
   &V_\pi^{(J,\sigma)}=V_\pi^\sigma\\
\S6&\charpi\\
   &\widehat{\slf}(\pi),\ \widehat f(X)\\
\S7.2&\lieG^0_\gss\\
\S7.4&\lieG_\gss\\
   &\sgsc\\
\S8.5&\admG_\gss(\epsilon),\ \admG_0\\
\S8.7&\temG_\gss(\epsilon)\\
   &\admG_\type,\ \temG_\type\\
\S8.8&\frB_d\\
   &\temG(\frt),\ \frt\in\frB_d\\
\S8.11&\bar{\temG_\gss}(\epsilon)\\
\S9.1&\distr(\lieG)\\
\S9.2&\distr^{i,\gss}_{x,-s},\ \distr^\gss\\
   &\Fcharpi\\
   &D_t,\ D_{-\stm_{d-1}}\\
\S11.1&\CaA\overset\mu=\CaA'\\
\S12.1&\datum_Y=(\vec\bG,\xo,\vec\rtm,\vec\phi,\rho)\\
   &\G^0_{[\xo]}\\
\S12.3&\rK^i,\ \rK^{i+}\\
   &J^i,\ J^i_+\\
\S12.4&\kappa,\ \rho_{\datum_Y}\\
\S13.1&\CaH(\G/\!/J,\sigma),\\
   &\supp(\CaH(\G/\!/J,\sigma))\\
\S13.2&\Phi(?,\bT,E)\\
\S13.9&\rK^i_{y\vdash},\ \tilde\rK^i_{y\vdash}\\
   &\rK_\vdash=\rK_{y\vdash}\\
   &\tilde\rK_\vdash=\tilde\rK_{y\vdash}\\
   &\overline U_{y,\vec\tm},\ M_{y,\vec\tm},\ U_{y,\vec\tm}\\
   &\overline U^i_{y,\vec\tm},\ M^i_{y,\vec\tm},\ U^i_{y,\vec\tm}\\
\S13.13& J_\ell,\ J_M,\ J_u\\
\S16.6&\tJ^i_{y+},\ \tJ^i_y\\
\S17.1&V_\pi^{\phi_\datum}\\
\S18&\rK,\ \rK^+\\
\S18.3&\phi_\datum^\rh,\ \phi_\datum^{\rh \pi},\ \tilde\phi_\datum^{\rh \pi}\\
& \\
\end{array}$

\end{document}